\documentclass[reqno]{amsart}
\usepackage{preprintStyle}
\pdfoutput=1 %

\usepackage{bm}
\usepackage{comment}
\usepackage{adjustbox} 

\usepackage{algorithmic}

\usepackage{amssymb}

\usepackage{bm}
\usepackage{comment}
\usepackage{adjustbox} 

\usepackage{mathtools}

\usepackage{xspace}

\usepackage{graphicx}
\usepackage{color}
\usepackage{circuitikz}
\usepackage{tikz}
\usetikzlibrary{calc,positioning,shapes}
\usetikzlibrary{patterns,decorations.pathmorphing,decorations.markings}
\usetikzlibrary{external}

\usepackage{pgfplots}
\pgfplotsset{compat=newest}
\usepackage[margin=10pt,font=small,labelfont=bf,labelsep=endash]{caption}
\usepackage{subcaption}

\usetikzlibrary{intersections, backgrounds}
\usetikzlibrary{circuits}
\usetikzlibrary{circuits.ee.IEC}
\usepackage{pgfplots}
\pgfplotsset{compat=newest}

\definecolor{color0}{rgb}{0.12156862745098,0.466666666666667,0.705882352941177}
\definecolor{color1}{rgb}{1,0.498039215686275,0.0549019607843137}
\definecolor{color2}{rgb}{0.172549019607843,0.627450980392157,0.172549019607843}
\definecolor{color3}{rgb}{0.83921568627451,0.152941176470588,0.156862745098039}
\definecolor{color4}{rgb}{0.580392156862745,0.403921568627451,0.741176470588235}
\definecolor{color5}{rgb}{0,0,0}

\definecolor{mycolor1}{rgb}{0.00000,0.44700,0.74100}
\definecolor{mycolor2}{rgb}{0.85000,0.32500,0.09800}
\definecolor{mycolor3}{rgb}{0.92900,0.69400,0.12500}
\definecolor{mycolor4}{rgb}{0.46600,0.67400,0.18800}
\definecolor{mycolor5}{rgb}{0.49400,0.18400,0.55600}

\newcommand{\lineWidth}{1.2pt}
\newcommand{\imageWidth}{2.0in}
\newcommand{\imageHeight}{1.8in}

\newcommand{\smoothFunctions}[3][]{\ifthenelse{\equal{#1}{}}{\mathcal{C}^{#2}}{\mathcal{C}_{#1}^{#2}}(#3)}

\newcommand{\dist}[1][]{\ifthenelse{\equal{#1}{}}{\mathbb{D}}{#1_{\mathbb{D}}}}

\newcommand{\T}{\top} 

\DeclareMathOperator{\diag}{diag}

\DeclareMathOperator{\sspan}{span}

\DeclareMathOperator*{\argmin}{arg\,min}

\DeclareMathOperator{\distance}{dist}
\newcommand {\UB}       {{\rm{UB}}}
\DeclareMathOperator*{\argmax}{arg\,max}
\newcommand{\dimPar}{p}

\newcommand{\calE}{\mathcal{E}}
\newcommand{\calF}{\mathcal{F}}

\newcommand{\calL}{\mathcal{L}}
\newcommand{\calM}{\mathcal{M}}

\newcommand{\calO}{\mathcal{O}}
\newcommand{\calP}{\mathcal{P}}

\newcommand{\calU}{\mathcal{U}}
\newcommand{\calV}{\mathcal{V}}
\newcommand{\calW}{\mathcal{W}}

\newcommand{\fv}{\ensuremath{\bm{v}}}

\newcommand{\ftheta}{\ensuremath{\bm{\theta}}}

\newcommand{\fLambda}{\ensuremath{\bm{\Lambda}}}

\newcommand{\zeroVec}{\mathbf{0}}

\newcommand{\abbr}[1]{\textsf{#1}\xspace}
\newcommand{\FOM}{\abbr{FOM}}
\newcommand{\MOR}{\abbr{MOR}}

\newcommand{\ROM}{\abbr{ROM}}

\newcommand{\SCM}{\abbr{SCM}}
\newcommand{\PDE}{\abbr{PDE}}
\newcommand{\PDEs}{\abbr{PDEs}}
\newcommand{\QSS}{\abbr{QSS}}
\newcommand{\RBM}{\abbr{RBM}}

\newcommand{\GM}{\abbr{GM}}
\newcommand{\EIGOPT}{\abbr{EigOpt}}

\newcommand{\ONB}{\abbr{ONB}}
\newcommand{\EIG}{\abbr{eig}}
\newcommand{\EIGS}{\abbr{eigs}}
\newcommand{\matlab}{{\sc Matlab}}

\newcommand{\fdim}{N} 
\newcommand{\rdim}{r} 

\newcommand{\N}{\ensuremath\mathbb{N}}

\newcommand{\R}{\ensuremath\mathbb{R}}
\newcommand{\C}{\ensuremath\mathbb{C}}

\newcommand{\IdMat}{\mathbf{I}}

\newcommand{\NullVector}{\mathbf{0}}

\newcommand{\Matrix}{\mathbf{A}}
\newcommand{\EigSpe}{\mathcal{W}}
\newcommand{\basisV}{\mathbf{V}}
\newcommand{\basisU}{\mathbf{U}}
\newcommand{\basisW}{\mathbf{W}}

\newcommand{\Proj}{\mathbf{P}}

\newcommand{\IndBasisV}[1]{\mathbf{V}_{\mkern-5mu #1}}
\newcommand{\IndBasisU}[1]{\mathbf{U}_{\mkern-5mu #1}}

\newcommand{\matB}{\mathbf{B}}

\newcommand{\matS}{\mathbf{S}}

\newcommand{\vek}{\mathbf{w}}
\newcommand{\Vek}{\mathbf{W}}
\newcommand{\uVec}{\mathbf{u}}
\newcommand{\eVec}{\mathbf{e}}

\newcommand{\yVec}{\mathbf{y}}
\newcommand{\zVec}{\mathbf{z}}

\newcommand{\ImgUnit}{\mathrm{i}}

\newcommand{\CompressedMat}[1]{\Matrix^{\mkern-6mu #1}}
\newcommand{\vekRed}[1]{\vek^{\mkern-3mu #1}}
\newcommand{\uVecRed}[1]{\uVec^{\mkern-3mu #1}}
\newcommand{\ResSubsp}[1]{\Res^{\mkern-3mu #1}}
\newcommand{\resSubsp}[1]{\res^{\mkern-1mu #1}}
\newcommand{\RedEigVal}[1]{\lambda^{\mkern-3mu #1}}
\newcommand{\PropRedEigVal}[1]{\tilde \lambda^{ #1}}

\newcommand{\indsub}{s}

\newcommand{\prmtr}{\bm \mu}
\newcommand{\prmtrs}{\mu}
\newcommand{\prmtrSet}{\calP}

\newcommand{\LB}{\mathrm{LB}}
\newcommand{\SLB}{\mathrm{SLB}}
\newcommand{\SUB}{\mathrm{SUB}}

\newcommand{\subspace}{\calV}
\newcommand{\ISubspace}[1]{\calV_{\mkern-2.5mu #1}}
\newcommand{\Res}{\mathbf{R}}
\newcommand{\res}{\mathbf{r}}

\newcommand{\TrainSet}{\Xi}

\usepackage{graphicx}
\usepackage{color}
\usepackage{circuitikz}
\usepackage{tikz}
\usetikzlibrary{calc,positioning,shapes}
\usetikzlibrary{patterns,decorations.pathmorphing,decorations.markings}
\usetikzlibrary{external}
\usepackage{pgfplots}
\pgfplotsset{compat=newest}
\usepackage[margin=10pt,font=small,labelfont=bf,labelsep=endash]{caption}
\usepackage{subcaption}

\usetikzlibrary{intersections, backgrounds}
\usetikzlibrary{circuits}
\usetikzlibrary{circuits.ee.IEC}
\usepackage{pgfplots}
\pgfplotsset{compat=newest}

\definecolor{color0}{rgb}{0.12156862745098,0.466666666666667,0.705882352941177}
\definecolor{color1}{rgb}{1,0.498039215686275,0.0549019607843137}
\definecolor{color2}{rgb}{0.172549019607843,0.627450980392157,0.172549019607843}
\definecolor{color3}{rgb}{0.83921568627451,0.152941176470588,0.156862745098039}
\definecolor{color4}{rgb}{0.580392156862745,0.403921568627451,0.741176470588235}
\definecolor{color5}{rgb}{0,0,0}

\definecolor{mycolor1}{rgb}{0.00000,0.44700,0.74100}
\definecolor{mycolor2}{rgb}{0.85000,0.32500,0.09800}
\definecolor{mycolor3}{rgb}{0.92900,0.69400,0.12500}
\definecolor{mycolor4}{rgb}{0.46600,0.67400,0.18800}
\definecolor{mycolor5}{rgb}{0.49400,0.18400,0.55600}

\title[\MOR for parametric Hermitian eigenproblems]{Certified Model Order Reduction for parametric Hermitian eigenproblems}
\author{Mattia Manucci${}^\star$, Benjamin Stamm${}^\dagger$  \and Zhuoyao Zeng${}^\dagger$}

\address{${}^{\star}$ Stuttgart Center for Simulation Science (SC SimTech), University of Stuttgart, Universit\"{a}tsstr.~32, 70569 Stuttgart, Germany}
\email{mattia.manucci@simtech.uni-stuttgart.de}
\address{${}^{\dagger}$ Institute of Applied Analysis and Numerical Simulation, University of Stuttgart, Pfaffenwaldring 57, 70569 Stuttgart, Germany}
\email{\{benjamin.stamm,zhuoyao.zeng\}@mathematik.uni-stuttgart.de}

\date{\today}

\begin{document}

\begin{abstract}
    This article deals with the efficient and certified numerical approximation of the smallest eigenvalue and the associated eigenspace of a large-scale parametric Hermitian matrix. 
    For this aim, we rely on projection-based model order reduction (MOR), i.e., we approximate the large-scale problem by projecting it onto a suitable subspace and reducing it to one of a much smaller dimension. 
    Such a subspace is constructed by means of weak greedy-type strategies. 
    After detailing the connections with the reduced basis method for source problems, we introduce a novel error estimate for the approximation  
    of the eigenspace associated with the smallest eigenvalue. 
    Since the difference between the second smallest and the smallest eigenvalue, the so-called spectral gap, is crucial for the reliability of the error estimate, we propose efficiently computable upper and lower bounds for higher eigenvalues and for the spectral gap, which enable the assembly of a subspace for the MOR approximation of the spectral gap. 
    Based on that, a second subspace is then generated for the MOR approximation of the eigenspace associated with the smallest eigenvalue. 
    We also provide efficiently computable conditions to ensure that the multiplicity of the smallest eigenvalue is fully captured in the reduced space. This work is motivated by a specific application: the repeated identifications of the states with minimal energy, the so-called ground states, of parametric quantum spin system models. 
\end{abstract}

\maketitle
{\footnotesize \textsc{Keywords:} Parametric eigenvalue problem, eigenspace approximation, model order reduction, greedy algorithm, certified error bound, spectral gap approximation, subspace procedure, quantum spin systems.} 

{\footnotesize \textsc{AMS subject classification: } Primary 41A63, 65F15; Secondary 65N25, 65F50, 65K10, 81P68.}


\section{Introduction}\label{sec:intro}
Given an analytic Hermitian matrix-valued function $\Matrix  :  \prmtrSet \rightarrow \C^{\fdim \times \fdim}$ defined on
a compact domain $\prmtrSet \subseteq  \R^{\dimPar}$, our aim is to obtain an accurate and efficient
approximation of the pair $(\lambda_{1}(\prmtr),\EigSpe_1(\prmtr))$ for all $\prmtr \in \prmtrSet$, where $\lambda_{1}(\prmtr)$ denotes the smallest eigenvalue of $\Matrix(\prmtr)$ and $\EigSpe_1(\prmtr)$ is the associated eigenspace. 
Note that, as $\prmtr$ varies in $\prmtrSet$, $\lambda_1(\prmtr)$ may not be simple, i.e. the algebraic multiplicity of $\lambda_1(\prmtr)$ could be larger than one, and so $\EigSpe_1(\prmtr)$ could be of dimension larger than one. 
For simplicity, we focus on the case $1 \ll \fdim < \infty$, although a generalization to the infinite-dimensional setting with compact self-adjoint operators is possible.

Our research is driven by the approximation problem that arises from the computation of the ground energy and the associated ground states of
parametric \emph{quantum spin systems} (\QSS)
models \cite{Herbst2022, Brehmer2023}. 
The study of quantum spin models, which dates back to the early days of quantum mechanics, remains a key focus in modern condensed matter physics. 
Indeed, as basic quantum many-body systems with inherently strong correlations, these models
often display interesting ground states that include complex ordering patterns, quantum disordered regimes, or topological order such as in quantum spin liquids \cite{DecLS11, HeiS22,SRFB2004}. 
In addition, these ground states serve in many cases as good approximations for the low-temperature behavior of real physical systems or compounds that these models aim to describe. 
Traditionally, several numerical instances must be performed in the Hamiltonian parameter space to map the phase diagram of the model. 
This is typically done by calculating relevant observables, such as the energy associated with the ground states of the model. Indeed, the smallest eigenvalue of the \QSS Hamiltonian corresponds to the minimal energy of the system, while its associated eigenspace corresponds to the states with minimal energy, the so-called ground states.
Computations may become easily prohibitive because of the exponential growth of the Hilbert space dimension with respect to the number of spins in a finite-size sample. 
Thus, solving the full-scale problem for \textit{many} different parameters using a standard eigensolver, such as the Lanczos method \cite{BaiDDRV00} or tensor methods \cite{Sch11}, may be computationally infeasible. 

More generally, the framework for approximating eigenvalues and the associated eigenspaces is of interest in all classes of parametric {partial differential equations} (\PDE) where one wants to study the spectral properties of the \PDE operator \cite{HWD17, GSH23}.

To address this problem, we present a certified projection-based \emph{model order reduction} (\MOR) framework for eigenproblems. 
Specifically, the approach we describe is divided into \emph{offline} and \emph{online} phases.
In the offline phase, a subspace $\subspace \subseteq  \C^{\fdim}$, also called a reduced space, of dimension $\rdim  \coloneqq \dim (\subspace) \ll\fdim$ is constructed at a computational cost that scales with $\fdim$. 
In the online phase, we compute the smallest eigenvalue $\RedEigVal{\subspace}_1(\prmtr)$ and the associated eigenspace $\EigSpe_1^\subspace(\prmtr)$ of the reduced Hermitian matrix-valued function $\CompressedMat{\basisV} \, : \, \prmtrSet \rightarrow \C^{\rdim\times \rdim}$ at $\prmtr \in \prmtrSet$, defined by
\begin{equation}\label{eqn:red:pro}
\Matrix^{\mkern-6mu\basisV}(\prmtr)		 	\;	\vcentcolon=	\;	
	\basisV^\ast \Matrix(\prmtr) \basisV,
\end{equation}
where $\basisV \in \C^{\fdim \times \rdim}$ is a matrix whose columns form an \emph{orthonormal basis} (\ONB) for $\subspace$. 
Once $\CompressedMat{\basisV} (\prmtr)$ is formed, the computational cost to evaluate $\RedEigVal{\subspace}_1(\prmtr)$ and its associated eigenspace is independent of $\fdim$ and scales with $\rdim\ll\fdim$.

As we want to approximate eigenspaces, it is crucial to state how the distance between two subspaces is measured. 
We follow the literature \cite{DavK70} using the standard Euclidean norm for vectors and the induced spectral norm for matrices, both denoted by $\|\cdot\|$. 
Consider two arbitrary subspaces $ \EigSpe,  \EigSpe' \subseteq  \C^\fdim$ and the associated matrices $ \basisW,  \basisW'$ such that their columns form an \ONB for $ \EigSpe$ and $ \EigSpe'$, respectively. 
One can define the orthogonal projectors onto the subspaces $ \EigSpe$ as 
\begin{equation} \label{eq:def projector}
\Proj^{ \EigSpe } \;\vcentcolon=\;  \basisW  \basisW^*,
\end{equation}
and similarly for $ \EigSpe'$ as $\Proj^{ \EigSpe'} \vcentcolon=  \basisW' ( \basisW')^*$. 
The distance between $ \EigSpe$ and $ \EigSpe'$ is then given as $\|\Proj^{ \EigSpe}- \Proj^{ \EigSpe'}\|$. 
Note that orthogonal projectors are independent of the choice of the basis \cite[Thm.~6.55]{Axl15}, and so $\|\Proj^{ \EigSpe}- \Proj^{ \EigSpe'}\|$ does not depend on the choice of the basis $ \basisW$ and $ \basisW'$. 
{Since $\EigSpe_1^\subspace(\prmtr) \subseteq \subspace \subseteq \C^\fdim$, we view $\EigSpe_1^\subspace(\prmtr)$ as a subspace of $\C^\fdim$ and can compare it with $\EigSpe_1(\prmtr)$ for the error analysis.}

By \emph{certified \MOR approximation} for the eigenpairs, we mean that for a {\emph{prescribed parameter subdomain}}  $\Xi \subseteq \prmtrSet$ and given arbitrary $\varepsilon_{\lambda}, \varepsilon_{\EigSpe}>0$, the constructed subspace $\subspace$ satisfies the following conditions:
\begin{equation}\label{eqn:goal}	
|\lambda_{1}(\prmtr)	- \RedEigVal{\subspace}_1	(\prmtr)| 
\;\le\; 
\varepsilon_{\lambda}
\quad  \text{and}\quad 
\|\Proj^{ \EigSpe_1(\prmtr)}- \Proj^{ \EigSpe_1^\subspace(\prmtr)}\|
\;\le\;
\varepsilon_{\EigSpe}
\qquad \text{for all }  \prmtr \in \Xi.
\end{equation}

The following result is crucial for understanding eigenspace approximations. 
\begin{theorem}(\cite[Sec.~2]{Wed83})\label{theo:err:mes}
    Consider two subspaces $ \EigSpe, \EigSpe' \subseteq \C^{\fdim}$ with associated \ONB matrices $ \basisW$, $\basisW'$, and orthogonal projectors $\Proj^{ \EigSpe}$, $\Proj^{ \EigSpe'}$ as defined in \eqref{eq:def projector}. Then
 \begin{enumerate}[label=\roman*)]
    \item  
    $\|\Proj^{ \calW}- \Proj^{ \calW'}\|\;=\;1$,  
        \hfill $\text{if} \quad \dim( \calW)>\dim( \calW')$;
    \item 	
    $\|\Proj^{ \calW}- \Proj^{ \calW'}\| \;=\; \|\Proj^{ \calW^{\perp}} \Proj^{ \calW'}\| \;=\; \|\Proj^{ \calW^{\perp}} {\basisW}' \|$,
        \hfill $\text{if} \quad \dim( \calW)=\dim( \calW')$;\label{thm:err:mes:item2}
    \item 
    $\|\Proj^{ \calW}- \Proj^{ \calW'}\| \;=\; \|\Proj^{ \calW^{\perp}} { \basisW }' \|\;=\; 1$,  
        \hfill $\text{if} \quad \dim( \calW)< \dim( \calW')$.\label{thm:err:mes:item3} 
\end{enumerate}
\end{theorem}
 From \Cref{theo:err:mes} it is clear that the necessary requirement to achieve \eqref{eqn:goal} for arbitrary $\varepsilon_{\EigSpe}$ is
 $\dim( \EigSpe_1(\prmtr) ) = \dim( \EigSpe_1^\subspace(\prmtr) )$ for all $\prmtr \in \Xi$, i.e. the dimension of the approximated eigenspace must be equal to the dimension of the exact eigenspace. 
 This is of fundamental relevance for the certified eigenspace approximation and will be discussed later in \Cref{sec:4.3}. 
In a framework where $\dim( \EigSpe ) = \dim( \EigSpe' )$ holds, $\|\Proj^{ \calW}- \Proj^{ \calW'}\|$ can be replaced by $\|\Proj^{ \calW^{\perp}} { \basisW}' \|$ due to \Cref{theo:err:mes}~\ref{thm:err:mes:item2}, {which makes drawing down the projection error $\|\Proj^{ \calW^{\perp}} { \basisW}' \|$ another necessary condition to achieve our goal.} {Therefore, provided that we can guarantee the exact dimension subspace recovery, we can arbitrarily switch between $\|\Proj^{ \calW}- \Proj^{ \calW'}\|$ and $\|\Proj^{ \calW^{\perp}} { \basisW}'\|$; this allow us to work with the error definition that is more convenient for our proposes. 
Following this consideration, we will first proceed by deriving an upper error bound for $\|\Proj^{ \calW^{\perp}} { \basisW}'\|$ (i.e. \Cref{thm:error bound EigVec}) and later present a result (i.e. \Cref{thm:exact:dim:rec}) that can be used to verify and impose the exact eigenspace dimension recovery.}

\subsection{Literature review and state-of-the-art}
Subspace approaches based on additional conditions on parameter dependencies have been proposed in \cite{MacMOPR00,PruRVP02}. In the context of 
eigenproblem optimization, subspace acceleration approaches have been discussed in 
\cite{AndS12, VliM12,KreV14,Meerbergen2017,KanMMM18,KreLV18,MehM24}. 

The subspace framework has several ideas shared with projection-based \MOR methods, which also start with the assumption that the solution set of the problem addressed can be well approximated by a small-dimensional subspace. 
Projection-based \MOR techniques can be distinguished into two main fields: the one related to the \emph{reduced basis method} (\RBM) \cite{HesRS16}, developed to address the fast approximation of parametric \PDE solutions, and the system-theoretic one, where the reduction of systems in control form is addressed; see \cite{BenGW15} for a detailed survey.

Within the context of \RBM, several approaches have been developed to approximate the smallest eigenvalue $\lambda_{1}(\prmtr)$ of a parametric Hermitian matrix $\Matrix(\prmtr)$ for $\prmtr \in {\prmtrSet}$, where $\Matrix(\prmtr)$ has the following affine-decomposition structure
\begin{equation}\label{eqn:affine}
	\Matrix(\prmtr)	\;	=	\;	\theta_1(\prmtr) \Matrix_1	+	\dots		+	\theta_Q(\prmtr) \Matrix_Q
\end{equation}
for given Hermitian matrices $\Matrix_1, \dots, \Matrix_Q \in {\mathbb C}^{\fdim\times \fdim}$ and real analytic scalar-valued functions $\theta_1, \dots , \theta_Q : {\R}^{\dimPar} \rightarrow {\R}$. 
The decomposition \eqref{eqn:affine} with $Q$ relatively small can be found in several important applications, such as all classes of linear {\PDE}s, and it is fundamental for the efficiency of projection-based \MOR, since it allows one to form the reduced matrix \eqref{eqn:red:pro} with a computational cost independent of $\fdim$ for all $\prmtr\in\prmtrSet$; see \cite[Sec.~3.3]{HesRS16}. The \emph{successive constraint method} (\SCM) \cite{HuyRSP07} is a well-known approach in the \RBM community to deal with the
approximation of the smallest eigenvalue of a parametric Hermitian matrix of type \eqref{eqn:affine}.
It is based on the construction of an upper bound $\lambda_1^{\UB}(\prmtr)$ and a lower bound $\lambda_1^{\LB}(\prmtr)$ for $\lambda_{1}(\prmtr)$, which are computationally cheaper to evaluate, i.e. they can be evaluated with few operations whose number is independent of $\fdim$. 
The accuracy of these approximations is iteratively improved via a 
greedy strategy, i.e., the parameter
\begin{equation}\label{eq:max_error2}
	\prmtr_{\mathrm{next}} \; = \;  \argmax_{\prmtr\in\Xi} \, (\lambda_1^{\UB}(\prmtr)-\lambda_1^{\LB}(\prmtr))
\end{equation}
is computed, where the maximization is performed over a discrete, finite set $\Xi\subseteq \prmtrSet$ chosen a priori.
Then $\lambda_1^{\UB}(\prmtr)$, $\lambda_1^{\LB}(\prmtr)$ are modified in such a way that they 
interpolate $\lambda_{1}(\prmtr)$ 
at $\prmtr = \prmtr_{\mathrm{next}}$. 
The \SCM is often slow to converge and lacks the Hermite interpolation 
property in the lower bound $\lambda_1^{\LB}(\prmtr)$, i.e., in general $\lambda_{1}(\prmtr_{\mathrm{next}})\neq \lambda_1^{\LB}(\prmtr_{\mathrm{next}})$ as well as 
$\nabla \lambda_{1}(\prmtr_{\mathrm{next}})\neq \nabla \lambda_1^{\LB}(\prmtr_{\mathrm{next}})$ when $\lambda_{1}(\prmtr)$ is differentiable at $\prmtr = \prmtr_{\mathrm{next}}$. 
To overcome this issue, Sirkovic and Kressner introduced in \cite{SirK16} a new lower bound $\lambda_1^{\SLB}(\prmtr,\subspace)$, for which the Hermite interpolation is shown to hold when $\lambda_{1}(\prmtr)$ is simple. 
Moreover, very recent work \cite{ManMG24} shows how to tackle the problem \eqref{eq:max_error2} in the whole continuum domain $\prmtrSet$ rather than in a discrete subset of it. 
Handling a continuum domain results in a larger computational time for the construction of the subspace, but it gives the advantage of selecting more suitably the parameters for the construction of the subspace, leading to subspaces with better approximation quality for a fixed accuracy; see \cite[Sec.~6.1]{ManMG24}.

In \cite{ManN15}, an approximation of $\lambda_{1}(\prmtr)$ is constructed with a heuristic strategy based on interpolant radial basis functions. 
In \cite{GugM23}, the authors have proposed a method based on the computation of 
\begin{equation*}
	\underline{\lambda} \; \vcentcolon= \;\min_{\prmtr \in \prmtrSet}\;\lambda_{\min}(\prmtr),
\end{equation*}
through gradient optimization. 
This is not suitable in our context, since it replaces $\lambda_{1}(\prmtr)$ with the uniform lower bound $\underline{\lambda}$ 
for all $\prmtr \in {\prmtrSet}$, while we need a reliable approximation for each instance of the parameter $\prmtr$ in the domain $\prmtrSet$. 
We also mention \cite{PraB24}, where the authors deal with parametric nonlinear eigenvalue problems utilizing contour integration-based strategies.

Fewer attempts have been made for the approximation of eigenspaces of parametric problems. 
A first step to tackle the problem, although not yet addressing the eigenspaces, was taken in \cite{MacMOPR00} where an output bound is presented for symmetric positive definite eigenvalue problems. 
{Approximation of eigenpairs of an infinite dimensional operator using sparse tensor grid method was introduced in \cite{AndS12}, in which a priori convergence of the method is also discussed. 
In \cite{GSH23}, a stochastic collocation method is employed to approximate isolated eigenspaces of parametric operators.}
A recent work \cite{AlgBB25} introduces a \MOR algorithm based on a sparse grid adaptive refinement, for the approximation of the eigensolutions of parametric problems arising from elliptic \PDEs. 
In \cite{Herbst2022,Brehmer2023} a \RBM framework is employed to deal with the fast approximation of the solution of a certain \QSS. 
The approximation space is constructed with a greedy-type algorithm \cite[Sec.~3.2.2]{HesRS16} which, however, uses an uncertified surrogate error estimate. 
In \cite{GarS24}, the authors consider a \RBM framework using derivatives of eigenvectors for a local convergence improvement.  

\subsection{Main contributions} 
We start with detailing the connections between \RBM for parametric source problems and the subspace framework for parametric Hermitian eigenproblems. 
In particular, with \Cref{example:no crossing}, we demonstrate that the difficulty in approximating parametric eigenproblems is linked to the existence of a large orthonormal system within the solution set intended for approximation. 
This is a shared feature with \RBM; see \cite[Sec.~5.1]{OhlR16}. 

Then, we introduce a novel error bound for the eigenspace {projection} error approximation; see \Cref{thm:error bound EigVec}. 
The effective evaluation of the error bound is related to the efficient and reliable approximation of the spectral gap, i.e. the difference between the second smallest eigenvalue (counted without multiplicities) and the smallest one of $\Matrix(\prmtr)$. 
This motivates us to find the computable and certified upper and lower bounds of the eigenvalues and of the spectral gap; see \Cref{prop:EigVal SLB}, \Cref{subsec:gap bounds}, and \Cref{alg:gap}, respectively. 

We extensively discuss how the proposed method behaves when $\lambda_1(\prmtr)$ is not simple in \Cref{sec:4.3} and we provide efficiently verifiable sufficient conditions to guaranty that the multiplicity of $\RedEigVal{\subspace}_1(\prmtr)$ coincides with that of $\lambda_1(\prmtr)$; see \Cref{thm:exact:dim:rec}. 

The subspace approximation provided by \Cref{alg:gap} is then used for the spectral gap approximation in a second greedy algorithm that constructs a suitable subspace for the approximation of the eigenspace using the a posteriori error estimate \eqref{eqn:new:err:est} and \Cref{prop1}; see \Cref{alg:sf}. To our knowledge, the method we propose is the only one in the literature able to efficiently construct approximation spaces that are able to capture the correct eigenspace dimension and that can be arbitrary accurate, thus ensuring approximation error control of the eigenspaces of parametric Hermitian matrices {on the prescribed training domain. Furthermore, for samples lying outside the training set, we provide an efficient and reliable criteria to efficiently identify potential discrepancies between the approximate eigenspace and the true eigenspace; see \Cref{rem:partial_certification}.} 

Finally, we present numerical simulations on artificial and \QSS examples, which validate the proposed certified projection-based \MOR framework.

\subsection{Organization of the manuscript} 
The remaining parts of this paper are organized as follows. 
In \Cref{sec2}, after detailing the problem setting and discussing some general assumptions, we recall some existing theoretical error bounds for eigenvalues and eigenvectors, and we introduce a novel theoretical bound for the eigenspace projection error. 
In \Cref{sec3}, we discuss practically computable upper and lower bounds for the eigenvalues and for the spectral gap, and we also introduce practical conditions to check whether the approximated eigenspaces fully capture the dimension of the original eigenspaces. 
In \Cref{sec4} the two-stage greedy strategy for the efficient and certified approximation of the eigenspaces is presented. 
\Cref{sec5} is to corroborate our theoretical findings with numerical results on synthetic examples and quantum spin models. 
In \Cref{conclusion}, we state our conclusions.

\subsection{General notations} \label{subsec:notation}
For any $n\in \N$, we denote by $\IdMat_n$ the identity matrix of size $n \times n$. 
$\|\cdot\|$ indicates the canonical Euclidean norm for vectors and the induced spectral norm for matrices. 

{The bold font is used for vectors and matrices, while the normal font is reserved for scalar quantities.}
For any $\prmtr\in\prmtrSet$, the eigenvalues of a Hermitian matrix $\Matrix(\prmtr) \in \C^{\fdim\times \fdim}$ will be denoted in ascending order and counting multiplicities by $\lambda_1(\prmtr) \leq \cdots \leq \lambda_\fdim(\prmtr)$, with the corresponding normalized eigenvectors $\vek_1(\prmtr), \dots, \vek_\fdim(\prmtr)$.

To deal with non-simple eigenvalues, we denote the number of distinct eigenvalues of $\Matrix(\prmtr)$ by $\tilde \fdim (\prmtr)$.
We also indicate by $ \tilde \lambda_1(\prmtr) < \cdots < \tilde \lambda_{\tilde \fdim (\prmtr)}(\prmtr)$ the \emph{different} eigenvalues of $\Matrix(\prmtr)$, each with multiplicity $m_k(\prmtr)$ for $k=1,\dots, \tilde \fdim (\prmtr)$, that is, $\fdim = \sum_{k=1}^{\tilde \fdim(\prmtr)} m_k(\prmtr)$ { and in particular $\lambda_1(\prmtr) = \dots = \lambda_{m_1(\prmtr)}(\prmtr) = \tilde \lambda_1(\prmtr) < \tilde \lambda_2(\prmtr)= \lambda_{m_1(\prmtr)+1}(\prmtr).$}  
The orthogonal projection onto the eigenspace $\tilde \EigSpe_k(\prmtr)$ associated with $ \tilde \lambda_k(\prmtr)$ is denoted by $ \Proj^{ \tilde \EigSpe_k(\prmtr)} \in \C^{\fdim \times \fdim}$ for $k=1, \dots , \tilde \fdim (\prmtr)$, so $\IdMat_\fdim = \sum_k \Proj^{ \tilde  \EigSpe_k(\prmtr)}$ by the spectral theorem. 
For simplicity, we also denote $\EigSpe_1(\prmtr) = \tilde \EigSpe_1(\prmtr)$. 

For an arbitrary matrix $\matB$, we represent with $\mathrm{Col}(\matB)$ the column space of $\matB$. 
Finally, for a given vector $\fv\in\C^{\fdim}$, $\diag(\fv)$ represents the square diagonal matrix with the elements of the vector $\fv$ on the main diagonal.

\section{Theoretical Foundations} \label{sec2}
 We start with formalizing the problem setting and highlighting connections to standard projection-based \MOR for parametric problems in \Cref{subsec:problem setting}. 
Then, we recall some existing theoretical error bounds for eigenvalues and eigenvectors in \Cref{subsec:EigVal EigVec bounds}.
In \Cref{subsec:bounds EigVec}, we propose and discuss a novel theoretical error bound for eigenspace approximations.

\subsection{Problem setting and general discussions}\label{subsec:problem setting}

We consider a parameter-dependent Hermitian matrix $\Matrix(\prmtr) \in \C^{\fdim \times \fdim}$ for $\prmtr \in \prmtrSet$, where $\prmtrSet \subseteq \R^p$ is compact and $\fdim \in \N$ is supposed to be large.  

Our goal is to efficiently approximate the smallest eigenpair $(\lambda_1(\prmtr),\EigSpe_1(\prmtr))$ for different values of~$\prmtr$. 
To describe $\EigSpe_1(\prmtr)$, which is characterized by its orthogonal projector $\Proj^{\EigSpe_1(\prmtr)}$, it is sufficient to know a complete set of orthonormal eigenvectors $\vek_1(\prmtr),\dots,\vek_{m_1(\prmtr)}(\prmtr)$ corresponding to $\lambda_1(\prmtr)$, see~\eqref{eq:def projector}. 
We note that the choice of such an orthonormal set of eigenvectors is not unique. However, since the spectral projections are independent of the specific choice of the eigenvectors basis, it is not restrictive to consider a particular representation of the eigenspace through an \ONB of eigenvectors.

Let us now introduce the corresponding solution set of the parametric eigenvalue problem 
\begin{align}\label{eqn:sol:set}
    \calM \;\vcentcolon=\; \left\{ \uVec  \in \C^\fdim \, \middle|  \,\,\|\uVec\|=1 \text{ and } \exists\,\prmtr\in\prmtrSet: \Matrix(\prmtr)\uVec = \lambda_1(\prmtr)\uVec \right\}
\end{align}
{and discuss its approximability.}
Note that the solution set $\calM$ is in general \textit{not} a manifold since the solution set could be locally not homeomorphic to a Euclidean space, as we show in \Cref{example:no manifold}.  

Following classical approximation theory for solution sets, the {(linear)} approximability degree of $\calM$ by an $n$-dimensional subspace of $\C^\fdim$ can be measured by
\begin{align}\label{eqn:KnW}
    d_n(\calM) \coloneqq \inf_{\substack{        \calL \subseteq \C^\fdim, \\ \dim(\calL)=n}} \sup_{\uVec \in \calM} \| \uVec- \Proj^\calL \uVec \|,  
\end{align}
i.e. the \textit{Kolmogorov $n$-width} of $\calM$. 
If $d_n(\calM)$ decays sufficiently fast as $n$ increases, then we say that \eqref{eqn:sol:set} is suitable to be approximated by a small-dimensional subspace.
However, as demonstrated by the upcoming \Cref{example:no crossing}, the solution set $\calM$ of a general parametric eigenvalue problem may have complicated behaviors and a slow decay of the associated Kolmogorov $n$-width. 

Next, we present two examples. 
The first one aims to stress that the solution set \eqref{eqn:sol:set} is in general not a manifold.
\begin{example} \label{example:no manifold}
Let $\eVec_1,\eVec_2,\eVec_3 \in \C^3$ be the canonical \ONB, $\prmtrSet \vcentcolon= [-2, \, 2] \subseteq \R$ and 
\begin{align*}
\Matrix : [-2, \, 2] \to \C^{3 \times 3}, \; \mu \mapsto 
\mathrm{diag}(\mu, \, \mu^2-2, \,  -\mu).
\end{align*}
Observe that $\lambda_1(\mu) $ has multiplicity $2$ for $\mu = \pm 1$, and that the solution set $\calM$ in the sense of \eqref{eqn:sol:set} is the unit circle of the $\eVec_1$-$\eVec_2$-plane unified with the unit circle of the $\eVec_2$-$\eVec_3$-plane. 
Hence, $\calM$ at $\eVec_2$ is locally not homeomorphic to any $\C^n$, and so by definition $\calM$ is not a manifold.
\end{example}

The second example aims to demonstrate that without restrictions or assumptions on the parametric eigenvalue problem, the approximation of the eigenspaces, or equivalently of the solution set $\calM$, using a subspace of small dimension could be of poor quality. 
\begin{example} \label{example:no crossing}
    Let $\prmtrSet= [-1,1] \subseteq \R$ and  $-1  \leq  x_1 < \dots < x_\fdim \leq 1 $ be pairwise different points. 
    For $k=1,\dots,\fdim$ let $L_k$ be the $k$-th Lagrange basis interpolation polynomial w.r.t. $\{x_k\}_{k=1}^\fdim$, in particular $L_j(x_k)=\delta_{j,k}$ for $j,k=1,\dots,\fdim$. 
    Note that each $L_k$ is a polynomial of degree $\fdim-1$ with zeros $\{x_j\}_{j=1,j\neq k}^\fdim$.
    Hence, we see that
    \begin{align}
        \sum_{j=1}^\fdim L_j^2(\mu) >0 \quad \text{ for all }\mu \in [-1,1]. \label{eq:sum of Lagrange polynomials squared}
    \end{align}
    Let $\{\eVec_k\}_{k=1}^\fdim$ be an arbitrary \ONB of $\C^\fdim$ and 
    \begin{align*}
        \vek : [-1,1] \rightarrow \C^\fdim, \; 
        \mu \mapsto \frac{1}{\sqrt{ \sum_{j=1}^\fdim L_j^2(\mu)}}  \sum_{k=1}^\fdim L_k(\mu)\eVec_k.
    \end{align*}
    Due to \eqref{eq:sum of Lagrange polynomials squared}, $\vek(\mu)$ is well-defined for every $\mu\in[-1,1]$, and it is immediate to verify that $\|\vek(\mu)\|=1$ for all $\mu\in[-1,1]$.
    Now, let us consider the parametric matrix 
    \begin{align*}
        \Matrix : 
        [-1,1]\rightarrow \C^{\fdim\times \fdim}, \, 
        \mu \mapsto - \vek(\mu) \, \vek(\mu) ^* = \frac{-1}{ \sum_{j=1}^\fdim L_j^2(\mu)}  \sum_{i,k=1}^\fdim L_i(\mu)L_k(\mu)\eVec_i \eVec_k^*.
    \end{align*}
    By construction, $\Matrix(\mu)$ is a real-analytic Hermitian matrix, which has a constant and simple smallest eigenvalue $\lambda_1(\mu)=-1$ with associated real analytic eigenvector function $\vek_1(\mu)=\vek(\mu)$.
    Since $\vek_1(x_k)=\eVec_k$, the solution set $\calM$ of this parametric matrix contains an \ONB of $\C^\fdim$. 
    Hence, we deduce with \cite[Rmk.~4.2]{GreU19} that the Kolmogorov $n$-width of $\calM$ decays sub-linearly, and so the approximation error w.r.t. $\calM$ using any subspace method decays poorly.
    In addition, there is another issue with this example, which leads to a poor approximation of $\calM$ by linear subspaces: By \cite[Thm.~3.7]{SirK16}, there exists $R>1$ s.t. approximating $\vek_1(\mu)$ using a subspace with parameters containing Chebyschev nodes has an error decay of order $R^{-J}$.  
    However, it turns out that the constant $R$ of this example is very close to $1$, e.g. $R\approx 1.003$ for $\fdim=14$.
    Thus, the approximation error of the eigenvector slowly decays despite the theoretical exponential error decay rate.
\end{example}
We should remark that in practical applications we often (but not always) observe good approximability of the solution set. 
{The previous example aims to provide an explicit construction for a Hermitian eigenvalue problem with poorly approximable solution set. }

{Methodologically, the approximation of $\calM$ with a subspace of small dimension is accomplished via Galerkin-projection methods. 
In more details, given a subspace $\subspace \subseteq \C^\fdim$ with $\dim(\subspace)\ll \fdim$ and a matrix $\basisV$ whose columns form an \ONB of $\subspace$, we consider the compressed matrix $\CompressedMat{\basisV} (\prmtr) \coloneqq \basisV^* \Matrix(\prmtr) \basisV $, and we calculate its smallest eigenvalue $\RedEigVal{\subspace}_1(\prmtr)$, which has multiplicity $m_1(\prmtr,\subspace)$, i.e. $\RedEigVal{\subspace}_1  = \ldots = \RedEigVal{\subspace}_{m_1(\prmtr,\subspace)}(\prmtr)$, as well as a complete set of orthonormal eigenvectors $\vekRed{\basisV}_1(\prmtr),\dots,\vekRed{\basisV}_{m_1(\prmtr,\subspace)}(\prmtr)$ corresponding to $\RedEigVal{\subspace}_1(\prmtr)$.
Then, we consider $( \RedEigVal{\subspace}_1(\prmtr), \EigSpe_1^\subspace(\prmtr) )$ as the approximation of $(\lambda_1(\prmtr),\EigSpe_1(\prmtr))$ w.r.t. the subspace $\subspace$, where $\EigSpe_1^\subspace(\prmtr)$ is a subspace of $\C^\fdim$, which has $(\basisV \vek_k^\basisV(\prmtr))_{k=1}^{m_1(\prmtr,\subspace)}$ as \ONB.}

{Practically, we take the following assumptions on our parametric matrix, and we will discuss their motivations or implications afterwards:} 
\begin{assumption}\label{general assumptions}
We assume that our problem satisfies the following properties.
\begin{enumerate}[label=\roman*)]
\item 
\label{assume:param sep} \textbf{Parameter separability}, i.e., {for a $Q \in \N $ with $Q \ll \fdim$ there exist $Q$} real-valued functions $\theta_1, \dots, \theta_Q : \prmtrSet \rightarrow \R$ and Hermitian matrices $\Matrix_1, \dots, \Matrix_Q \in \C^{\fdim\times \fdim}$, s.t. 
\begin{align}
\Matrix(\prmtr) = \sum_{q=1}^{Q} \theta_q(\prmtr) \, \Matrix_q \quad \text{ for all } \prmtr \in \prmtrSet. \label{eq:param sep}
\end{align}
\item 
\label{assume:analytic} \textbf{Real-analyticity in the parameter}, i.e., {both the real and the imaginary part of} each entry of the matrix-valued map $\Matrix : \prmtrSet \rightarrow \C^{\fdim \times \fdim}, \prmtr \mapsto \Matrix(\prmtr)$ are real-analytic in $\prmtr$.  
\item 
\label{assume:separated lower spectrum} \textbf{Separation of the lower part of the spectrum}, i.e., there exists $m_0 \in \N$ with $2\le m_0 \ll \fdim$, such that for all $\prmtr \in \prmtrSet$ we have $\tilde \lambda_{2} (\prmtr)  \leq \lambda_{m_0}(\prmtr) < \lambda_{m_0 +1}(\prmtr)${, where we recall $\tilde \lambda_{2} (\prmtr) = \lambda_{m_1(\prmtr)+1}(\prmtr)>\lambda_1(\prmtr)$ denotes the second smallest eigenvalue of $\Matrix(\prmtr)$ counted without multiplicities}.
\end{enumerate}
\end{assumption}

\Cref{general assumptions}~\ref{assume:param sep} is common and has both practical and theoretical utility. 
This condition, also referred to as affine decomposition, enables efficient computations, as we shall see in the next sections. 
We point out that many parametric matrices in applications, including the \QSS models and the matrices arising from the discretization of linear parametric \PDEs, have the form \eqref{eq:param sep}.

{\Cref{general assumptions} \ref{assume:analytic} is also common to take in this context.
It is equivalent to the functions  $\theta_1, \dots, \theta_Q$ being real-analytic, or requiring that the complex extension of the map $\Matrix$ is holomorphic. 
Let us discuss the parameter dependency of $\lambda_1(\prmtr)$ and $\EigSpe_1(\prmtr)$ based on this assumption. 
Since $\Matrix(\prmtr)$ is continuous in $\prmtr$, the smallest eigenvalue $\lambda_1(\prmtr) = \min_{\mathbf{x}\in \C^\fdim, \|\mathbf{x}\|=1} \mathbf{x}^* \Matrix(\prmtr) \, \mathbf{x}$ is also continuous. }
Moreover, if $\lambda_1(\prmtr)$ is simple, then it is also arbitrarily differentiable at $\prmtr \in\prmtrSet$ \cite{Kato76, Bau84}. 
In the case of a one-dimensional parameter space $\prmtrSet \subseteq \R$,  $\lambda_1(\mu)$ is piecewise real-analytic, where the non-analytic points are the points where eigenvalue crossings take place \cite[Sec.~2.5.7]{Kato76}.  
The regularity of $\EigSpe_1(\prmtr)$ is more nasty, see \cite{DoeE24} for more discussions.  
We note that, whenever the multiplicity of $\lambda_1(\prmtr)$ changes, the eigenspace changes discontinuously, as demonstrated by \Cref{example:no manifold}. 

As we also want to tackle the situation where $\dim(\EigSpe_1(\prmtr))>1$ for some parameter values, it turns out that we need an assumption on the degeneracy level of not only $\EigSpe_1(\prmtr) = \tilde \EigSpe_1(\prmtr)$ {(the eigenspace of the smallest eigenvalue)} but also of $\tilde \EigSpe_2(\prmtr)$ {(the eigenspace of $\tilde \lambda_2(\prmtr)$ the  second eigenvlue counted without multiplicities)}, even though the value $m_0$ is practically difficult to know a priori. 
Notice that \Cref{general assumptions}~\ref{assume:separated lower spectrum} implies $\dim (\EigSpe_1(\prmtr)) +\dim (\tilde \EigSpe_2(\prmtr)) \leq m_0 $ for all $\prmtr\in \prmtrSet$.

\subsection{Existing eigenvalue and eigenvector bounds} \label{subsec:EigVal EigVec bounds}

We recall and discuss known error bounds for the smallest eigenvalue and associated eigenvector approximations. 
Let us denote 
\begin{equation}\label{eqn:lifted:eigvec}
    \uVecRed{\basisV}_1  (\prmtr) 
    \;\vcentcolon=\; 
    \basisV \, \vekRed{\basisV}_1(\prmtr),
\end{equation}
where $  \vekRed{\basisV}_1 (\prmtr)$ refers to a normalized eigenvector associated to the smallest eigenvalue $\RedEigVal{\subspace}_1(\prmtr)$ of the compressed matrix $\CompressedMat{\basisV} (\prmtr)=\basisV^*\Matrix(\prmtr)\basisV$. 
The vector $\uVecRed{\basisV}_1  (\prmtr)$ can be seen as the approximated eigenvector $\vekRed{\basisV}_1(\prmtr) \in \C^{\rdim}$ with $\rdim\vcentcolon=\dim(\subspace)$, embedded into the full space, as it holds that 
\begin{align*}
    \RedEigVal{\subspace}_1(\prmtr) = \uVecRed{\basisV}_1 (\prmtr)^* \Matrix(\prmtr) \, \uVecRed{\basisV}_1 (\prmtr) =  \vekRed{\basisV}_1(\prmtr)^* \CompressedMat{\basisV} (\prmtr) \, \vekRed{\basisV}_1(\prmtr).
\end{align*}
We can also interpret $\uVecRed{\basisV}_1  (\prmtr)$ as an eigenvector of $\Proj^\subspace \Matrix(\prmtr) \Proj^\subspace$ w.r.t. the eigenvalue $\RedEigVal{\subspace}_1(\prmtr)$.
Also note that the choice of $\vekRed{\basisV}_1(\prmtr)$, and consequentially of $\uVecRed{\basisV}_1  (\prmtr)$, is in general not unique.
However, the choice of $\uVecRed{\basisV}_1  (\prmtr)$ does not affect the validity of the upcoming results. 

\begin{remark}\label{rem:no galerkin}
    Note that in general $\lambda_1(\prmtr) \neq \RedEigVal{\subspace}_1(\prmtr)$, and so Galerkin-orthogonality does not hold for parametric eigenvalue problems, i.e. 
\begin{align*}
    \left(\vek_1(\prmtr)-\uVecRed{\basisV}_1(\prmtr)\right)^* \Matrix(\prmtr) \, \uVec
    \;=\; 
    \left(\lambda_1(\prmtr) \vek_1(\prmtr)- \RedEigVal{\subspace}_1(\prmtr) \uVecRed{\basisV}_1(\prmtr) \right)^* \uVec
    \;\neq\;0, \quad \text{for }\uVec\in \subspace.
\end{align*}
This constitutes a notable difference between \RBM for eigenvalue problems with respect to \RBM for source problems \cite[Sec.~3]{HesRS16}. 
Indeed, the absence of Galerkin orthogonality makes conventional residual based error estimates no longer applicable in a straightforward manner, which is a well-known problem in the numerical analysis of eigenvalue problems \cite{Boffi_2010}.
\end{remark} 

To present the known error bounds, let us denote the eigenproblem residual w.r.t. $\subspace$ at $\prmtr \in \prmtrSet$ as
\begin{equation}\label{eqn:res:eig}
    \resSubsp{\subspace}(\prmtr)
    \;\vcentcolon=\; 
    \Matrix(\prmtr) \, \uVecRed{\basisV}_1(\prmtr) - \RedEigVal{\subspace}_1(\prmtr) \, \uVecRed{\basisV}_1 (\prmtr).
\end{equation}
Since we only concern about $\| \resSubsp{\subspace}(\prmtr)\|$, which is independent of the choice of $\basisV$, as $\basisV$ has orthonomal columns, we omit the dependency of $\basisV$ in the notation of $\resSubsp{\subspace}(\prmtr)$. 
Finally, let us also define the so-called \textit{spectral gap} of $\Matrix(\prmtr)$ by 
\begin{equation}
    \gamma (\prmtr) \;\vcentcolon=\; \tilde \lambda_2(\prmtr) - \lambda_1(\prmtr), \label{eq:gap full space}
\end{equation}
where we recall $\tilde \lambda_2(\prmtr)$ is the smallest eigenvalue of $\Matrix(\prmtr)$ such that $\tilde \lambda_2(\prmtr)>\lambda_1(\prmtr)$, that is $\tilde \lambda_2(\prmtr)=\lambda_{m_1(\prmtr)+1}(\prmtr)$ for $m_1(\prmtr)$ the algebraic multiplicity of $\lambda_1(\prmtr)$. 

\begin{lemma} 
\label{lemma:EigVal EigVec bound} 
{Let $\prmtr \in \prmtrSet$ be an arbitrary parameter and let $\sigma(\Matrix(\prmtr)) \subseteq \R$ denote the spectrum of $\Matrix(\prmtr)$. 
Then, the following properties hold true.}
\begin{enumerate}[label=\roman*)]
\item 
\label{sub-lemma:BF} (Theorem of Bauer-Fike \cite[Thm.~9.2]{Atk89}) Let $\hat \lambda(\prmtr) \vcentcolon= \argmin_{\lambda \in \sigma(\Matrix(\prmtr))} |\lambda-\RedEigVal{\subspace}_1(\prmtr)| $. 
Then it holds that 
\begin{align*}
    |\hat \lambda(\prmtr)-\RedEigVal{\subspace}_1(\prmtr)|\; \leq\; \| \resSubsp{\subspace}(\prmtr) \|.
\end{align*}
\item 
\label{sub-lemma:KT} (Theorem of Kato-Temple \cite[Thm.~11.7.1]{Atk89}) For 
\begin{equation*}
    \hat \lambda(\prmtr)
    \; \vcentcolon=\; 
    \argmin_{\lambda \in \sigma(\Matrix(\prmtr))} |\lambda - \RedEigVal{\subspace}_1(\prmtr)|
    \quad\text{and}\quad 
    \hat \gamma(\prmtr)
    \;\vcentcolon=\; 
    \min_{\lambda' \in \sigma(\Matrix(\prmtr))\backslash\{\hat \lambda(\prmtr)\}} |\lambda' - \RedEigVal{\subspace}_1(\prmtr) |, 
\end{equation*}
 it holds that   
\begin{align*}
     |\hat \lambda(\prmtr)-\RedEigVal{\subspace}_1(\prmtr)|\; \leq\; \frac{\| \resSubsp{\subspace}(\prmtr) \|^2}{\hat \gamma(\prmtr)}.
\end{align*}
\item 
\label{sub-lemma:EigVal by EigVec} The equality  
\begin{align} \label{eq:EigVal expression exact}
    |\RedEigVal{\subspace}_1(\prmtr)-\lambda_1(\prmtr)|
    =  (\vek_1(\prmtr) - \uVecRed{\basisV}_1 (\prmtr))^*(\Matrix(\prmtr) -  \lambda_1(\prmtr)\IdMat_\fdim )(\vek_1(\prmtr) - \uVecRed{\basisV}_1 (\prmtr))
\end{align}
holds true. 
By setting $C(\prmtr) \vcentcolon= \lambda_\fdim(\prmtr)-\lambda_1(\prmtr)$, we obtain  
\begin{align}
    |\RedEigVal{\subspace}_1(\prmtr)-\lambda_1(\prmtr)|
    \;\leq\; C(\prmtr) \| \vek_1(\prmtr) - \uVecRed{\basisV}_1 (\prmtr) \|^2. \label{ineq:EigVal by EigVec}
\end{align}
\item 
\label{sub-lemma:EigVec by EigVal} 
For $\Proj^{\EigSpe_{1}^{\perp}(\prmtr)}$ the orthogonal projector onto the orthogonal complement of the eigenspace $\calW_1(\prmtr)$, it holds that  
\begin{align} \label{eq:EigVec by EigVal}
    \| \Proj^{\EigSpe_{1}^{\perp}(\prmtr)} \, \uVecRed{\basisV}_1 (\prmtr) \|^2 \leq 
    \frac{\RedEigVal{\subspace}_1(\prmtr)-\lambda_1(\prmtr)}{\gamma(\prmtr)}. 
\end{align}
\end{enumerate}
\end{lemma}
\begin{proof}[Proof of \ref{sub-lemma:EigVal by EigVec} and \ref{sub-lemma:EigVec by EigVal}]
The equality \eqref{eq:EigVal expression exact} can be verified using 
\begin{equation*}
\lambda_1(\prmtr) = \vek_1(\prmtr)^*\Matrix(\prmtr)\,\vek_1(\prmtr),
\quad 
\RedEigVal{\subspace}_1(\prmtr) = \uVecRed{\basisV}_1 (\prmtr)^*\Matrix(\prmtr)\,\uVecRed{\basisV}_1 (\prmtr),
\end{equation*}
and then \eqref{ineq:EigVal by EigVec} follows directly. 
{To derive \eqref{eq:EigVec by EigVal}, we write $\uVecRed{\basisV}_1 (\prmtr) = \sum_{j=1}^\fdim c_j (\prmtr) \vek_j(\prmtr)$ as a linear combination of a basis of orthonormal eigenvectors of $\Matrix(\prmtr)$. 
By the definition of orthogonal projector we have $\| \Proj^{\EigSpe_{1}^{\perp}(\prmtr)} \, \uVecRed{\basisV}_1 (\prmtr) \|^2 = \sum_{j=m_1(\prmtr)+1}^\fdim |c_j(\prmtr)|^2$ , and observe that 
\begin{align*}
\lambda_1^{\subspace}(\prmtr) - \lambda_1 (\prmtr) 
 \;&=\; 
 \uVecRed{\basisV}_1 (\prmtr)^*\Matrix(\prmtr)\,\uVecRed{\basisV}_1 (\prmtr) - \lambda_1 (\prmtr) \uVecRed{\basisV}_1 (\prmtr)^*\,\uVecRed{\basisV}_1 (\prmtr) \\
 \;&=\; \sum_{j=m_1(\prmtr)+1}^\fdim (\lambda_j - \lambda_1) |c_j(\prmtr)|^2
 \geq 
  \gamma(\prmtr) \| \Proj^{\EigSpe_{1}^{\perp}(\prmtr)} \, \uVecRed{\basisV}_1 (\prmtr) \|^2.
\end{align*}} 
\end{proof}

Note that all the previous statements can be generalized for infinite-dimensional self-adjoint compact operators. 

The Bauer-Fike bound in \Cref{lemma:EigVal EigVec bound}~\ref{sub-lemma:BF} is computationally efficient to evaluate \cite{Herbst2022, BuhEOR14}. 
The drawbacks of this bound are twofold. 
For one thing, it only guarantees $\RedEigVal{\subspace}_1(\prmtr)$ to be close to some eigenvalue of $\Matrix(\prmtr)$, and in particular not necessarily to $\lambda_1(\prmtr)$. 
For another, it is asymptotically not sharp. 
The Kato-Temple bound (i.e. \Cref{lemma:EigVal EigVec bound}~\ref{sub-lemma:KT}) is an improvement of the Bauer-Fike bound, in the sense that with the additional information of $\hat\gamma(\prmtr)$, the distance of $\RedEigVal{\subspace}_1(\prmtr)$ from its second closest element in the spectrum of $\Matrix(\prmtr)$, the eigenvalue error is bounded by the squared residual divided by $\hat\gamma(\prmtr)$. 
However, $\hat\gamma(\prmtr)$ is in general not computable. 
Besides, the Kato-Temple bound, as another residual based bound, shares the same drawbacks as the Bauer-Fike bound, see also \Cref{rem:no galerkin}. 

The equality in \Cref{lemma:EigVal EigVec bound} \ref{sub-lemma:EigVal by EigVec} establishes the general relationship between the approximation of $\lambda_1(\prmtr)$ and that of $\vek_1(\prmtr)$.
It can be roughly interpreted as that, in case of convergence, the approximation of $\RedEigVal{\subspace}_1(\prmtr)$ converges quadratically faster than the approximation of $\vek_1(\prmtr)$. 
As noted in \Cref{lemma:EigVal EigVec bound}~\ref{sub-lemma:EigVec by EigVal}, deducing an error bound of $\vek_1(\prmtr)$ by that of $\lambda_1(\prmtr)$ requires additional information on the spectral gap \eqref{eq:gap full space}, which is positive by definition, but could be very small for some parameters.
When dealing with a small spectral gap, it is possible to achieve an accurate estimate of $\lambda_1(\prmtr)$; however, the error in the approximation of $\vek_1(\prmtr)$ may still be significant due to the orthogonality of the eigenvectors associated with distinct eigenvalues. Therefore, accurately estimating the spectral gap is crucial for reliable eigenvector approximations. Moreover, we also observe that the eigenvector error on the left-hand side of \Cref{lemma:EigVal EigVec bound}~\ref{sub-lemma:EigVec by EigVal} is squared.
This means that given an approximate eigenpair with eigenvalue approximation error of the order of the machine precision $\calO ({\varepsilon})$, we could only expect an eigenvector approximation error of order $ \calO (\sqrt{\varepsilon})$, incurring in the so-called square-root effect.

The discussed theoretical and numerical drawbacks motivate us to derive a novel error bound for the eigenspace approximation. 

\subsection{A novel error bound for eigenspaces} \label{subsec:bounds EigVec}
In this subsection, we introduce a {theoretical error bound for the projection error of the approximated eigenspace onto the true eigenspace}
and outline an approach to make it practically feasible.

\begin{theorem}[{Upper bound for the eigenspace projection error}] \label{thm:error bound EigVec}
Let $\prmtr \in \prmtrSet$ be an arbitrary parameter. 
Let $\EigSpe_1(\prmtr)$ be the eigenspace associated with the smallest eigenvalue $\lambda_1(\prmtr)$ of $\Matrix(\prmtr)$.
Let $\subspace \subseteq \C^\fdim$ be a subspace of $\dim(\subspace)=\rdim$, and let $\basisV$ be the matrix whose columns form an \ONB of $\subspace$. 
Consider the matrix $\CompressedMat{\basisV} (\prmtr) $ defined in \eqref{eqn:red:pro}, its smallest eigenvalue $\RedEigVal{\subspace}_1(\prmtr)$ with algebraic multiplicity denoted by $m_1(\prmtr,\subspace)\in\N$, and the eigenvectors $\vekRed{\basisV}_1(\prmtr),\dots,\vekRed{\basisV}_{m_1(\prmtr,\subspace)} (\prmtr)$ associated with $\RedEigVal{\subspace}_1(\prmtr)$.
Let $\EigSpe_1^\subspace(\prmtr) \vcentcolon= \sspan \big\{ \basisV \, \vekRed{\basisV}_k(\prmtr) \big\}_{k=1}^{m_1(\prmtr,\subspace)}$  and $\basisW^{\basisV}_1(\prmtr) \in \C^{\fdim\times m_1(\prmtr,\subspace)}$ be the matrix whose $k$-th column is equal to $\basisV \,\vekRed{\basisV}_k(\prmtr)$. 
Finally, let us define the residual matrix \begin{equation}\label{eqn:Res}
    \ResSubsp{\subspace}(\prmtr)
    \;\vcentcolon=\;
    \Matrix(\prmtr)\, \basisW^{\basisV}_1(\prmtr) - \RedEigVal{\subspace}_1(\prmtr) \, \basisW^{\basisV}_1(\prmtr).
\end{equation}
Then, it holds that
\begin{equation}\label{eqn:bound:eig:vec}
    \| \Proj^{\EigSpe^{\bot}_1(\prmtr)} \basisW^{\basisV}_1(\prmtr) \|
    \; \leq\; 
    \frac{1}{\gamma(\prmtr)}\left(\RedEigVal{\subspace}_1(\prmtr) - \lambda_1(\prmtr) 
    \;+ \;
    \|\ResSubsp{\subspace}(\prmtr)\|\right), 
\end{equation}
where $\gamma(\prmtr)$ is the spectral gap defined in \eqref{eq:gap full space}.
\end{theorem}
Note that \Cref{general assumptions}~\ref{assume:separated lower spectrum} implies $\Matrix(\prmtr)\neq c \, \IdMat_\fdim$ for all $c\in \R$ and for all $\prmtr\in\prmtrSet$, and so the spectral gap $\gamma(\prmtr)$ of $\Matrix(\prmtr)$ is assumed to be well-defined for all $\prmtr\in\prmtrSet$.
\begin{proof}
Recall that by the spectral theorem we have
\begin{equation*}
     \left(\Matrix(\prmtr)-\lambda_1(\prmtr)\IdMat_\fdim\right)\;=\;\sum_{i=1}^{\fdim}\left(\lambda_i(\prmtr)-\lambda_1(\prmtr)\right)\vek_i(\prmtr)\vek_i^{\ast}(\prmtr).
\end{equation*}
Besides, let us introduce the following matrix
\begin{equation}\label{eqn:pseudoinv}
     \left(\Matrix(\prmtr)-\lambda_1(\prmtr)\IdMat_\fdim\right)^{\dagger}\;\vcentcolon=\;\sum_{i=m_1(\prmtr)+1}^{\fdim}\frac{1}{\lambda_i(\prmtr)-\lambda_1(\prmtr)}\vek_i(\prmtr)\vek_i^{\ast}(\prmtr),
\end{equation}
which can be verified to be the Moore-Penrose pseudoinverse of $\Matrix(\prmtr)-\lambda_1(\prmtr)\IdMat_\fdim$, and it  holds that  
\begin{align} \label{eq:orth proj identity}
\begin{aligned}
    \left(\Matrix(\prmtr)-\lambda_1(\prmtr)\IdMat_\fdim\right)^{\dagger} \left(\Matrix(\prmtr)-\lambda_1(\prmtr)\IdMat_\fdim\right)
    \;=&\;
    \sum_{i=m_1(\prmtr)+1}^{\fdim}\vek_i(\prmtr)\vek_i^{\ast}(\prmtr)\\
    =&\;
    \IdMat_\fdim-\Vek_1(\prmtr)\Vek_1(\prmtr)^{\ast}
    \;=\;
    \Proj^{\EigSpe^{\bot}_1(\prmtr)}.
    \end{aligned}
\end{align}
From \eqref{eqn:pseudoinv}, we also see that
\begin{equation} \label{eq:pseudo inverse norm}
    \Big\|\left(\Matrix(\prmtr)-\lambda_1(\prmtr)\IdMat_\fdim\right)^{\dagger}\Big\|\;=\;\frac{1}{\lambda_{m_1(\prmtr)+1}-\lambda_1(\prmtr)}\;=\;\frac{1}{\gamma(\prmtr)}. 
\end{equation}
Now, applying \eqref{eq:orth proj identity} to $\basisW^{\basisV}_1(\prmtr)$, we find 
\begin{align} \label{eq: bound intermediate step}
\begin{aligned}
    \Proj^{\EigSpe^{\bot}_1(\prmtr)} \basisW^{\basisV}_1(\prmtr)
\,=&\, 
 \left(\Matrix(\prmtr)-\lambda_1(\prmtr)\IdMat_\fdim\right)^{\dagger}\left(\Matrix(\prmtr)-\lambda_1(\prmtr)\IdMat_\fdim+\RedEigVal{\subspace}_1\IdMat_\fdim-\RedEigVal{\subspace}_1\IdMat_\fdim\right)\basisW^{\basisV}_1(\prmtr) \\
 =&\,\left(\Matrix(\prmtr)-\lambda_1(\prmtr)\IdMat_\fdim\right)^{\dagger}\left( (\RedEigVal{\subspace}_1(\prmtr) - \lambda_1(\prmtr)) \basisW^{\basisV}_1(\prmtr) +  \ResSubsp{\subspace}(\prmtr) \right).
 \end{aligned}
\end{align}
Inequality \eqref{eqn:bound:eig:vec} follows from taking the norm of both sides of \eqref{eq: bound intermediate step}, then applying the submultiplicative property of the spectral norm, \eqref{eq:pseudo inverse norm}, $\| \basisW^{\basisV}_1(\prmtr)\| = 1$, and the triangle inequality. 
\end{proof}

To apply the result of \Cref{thm:error bound EigVec} as an a posteriori error estimator for a greedy-type algorithm, it is crucial that all terms on the right-hand side of \eqref{eqn:bound:eig:vec} can be computed efficiently, specifically with a computational cost that does not depend on $\fdim$. Let us discuss this further.
\begin{enumerate}
    \item 
    The residual norm $\|\ResSubsp{\subspace}(\prmtr)\|$ can be evaluated in an efficient and stable manner thanks to \Cref{general assumptions}~\ref{assume:param sep}; see \cite{BuhEOR14} for more details. 
    Note that in case $m_1(\prmtr,\subspace)=1$ we have $\|\ResSubsp{\subspace}(\prmtr)\|=\|\resSubsp{\subspace}(\prmtr)\|$ for $\resSubsp{\subspace}(\prmtr)$ in \eqref{eqn:res:eig}.
    \item 
    The evaluation of the smallest eigenvalue $\lambda_1(\prmtr)$ is replaced by an efficiently computable subspace lower bound $\lambda_1^{\rm{SLB}}(\prmtr, \calV)$ that we will recall in \Cref{subsec:SLB}, see \eqref{eq:defn_LB}.
    \item 
    We approximate the spectral gap \eqref{eq:gap full space} in an efficient way by means of an independent greedy-subspace procedure. This is detailed in the upcoming \Cref{subsec:gap bounds} and \Cref{subsec:stage1}.
    \item In \Cref{subsec:stage2}, we will combine the previous points and derive an efficiently computable a posteriori error bound for the greedy approximation of eigenspaces; see \Cref{prop1}. 
\end{enumerate}

\section{Practical error bounds for eigenvalues and spectral gap} \label{sec3}

Let us recall the snapshot-based assembly of a suitable subspace for the approximation of eigenpairs. 
To construct a subspace efficiently, we employ a greedy procedure, i.e. we construct a subspace with the exact eigenvectors of some parameters, which are selected based on an efficiently computable error indicator function that, under certain detailed conditions, we show to be an upper bound of the true error. 
More precisely, suppose we have already selected some parameters $\prmtr_1, \dots, \prmtr_J \in \prmtrSet$. 
For $1 \leq \ell(j) \ll \fdim$ with $j=1,\dots,J$, we set 
\begin{align} \label{eq:low_bound}
    \subspace 
    \,\vcentcolon=\, 
    \ISubspace{J} 
    \, \vcentcolon=\, 
    \sspan \left\{ \vek_1(\prmtr_1), \dots, \vek_{\ell(1)}(\prmtr_1),	\dots	,	\vek_1(\prmtr_J), \dots, \vek_{\ell(J)}(\prmtr_J) \right\},
\end{align}
where we recall that $\vek_k(\prmtr_j)$ denotes an eigenvector of $\Matrix(\prmtr_j)$ corresponding to its
$k$-th smallest eigenvalue counting multiplicity, denoted as $\lambda_k(\prmtr_j)$, for $j = 1, \dots , J$ and $k = 1, \dots , \ell(i)$. 
Note that the number $\ell(j)$ in the assembly of $\ISubspace{J}$ can vary depending on the selected parameter~$\prmtr_j$: As we will see in \Cref{sec4}, we have $\ell(j) = m_1(\prmtr_j) + m_2(\prmtr_j)$ for the approximation of the spectral gap and $\ell(j) = m_1(\prmtr_j)$ for the approximation of $\EigSpe_1(\prmtr)$. 
 
Let $\basisV \coloneqq \IndBasisV{J}$ be a matrix whose columns form an \ONB of $\ISubspace{J}$. 
For any $\prmtr \in \prmtrSet$, we consider the matrix $\CompressedMat{\basisV} (\prmtr) = \basisV^* \Matrix(\prmtr) \basisV$, whose smallest eigenvalue and corresponding eigenspace will be denoted by $\RedEigVal{\subspace}_1 (\prmtr)$ and $\EigSpe_1^\subspace(\prmtr)$, respectively. 
Then, we approximate $\lambda_1(\prmtr)$ by $\RedEigVal{\subspace}_1 (\prmtr)$ as well as $\EigSpe_1(\prmtr)$ by $\EigSpe_1^\subspace(\prmtr)$.

To select the next parameter $\prmtr_{J+1}$ for the enlargement of the subspace, we use an error indicator function $\Delta_{\ISubspace{J}}:\prmtr \to [0,\infty)$, which should be evaluated at a computational cost independent of $\fdim$, and we choose $\prmtr_{J+1} = \argmax_{\prmtr \in \TrainSet} \Delta_{\ISubspace{J}}(\prmtr)$. 
The algorithms terminate when $\Delta_{\ISubspace{J}}(\prmtr_{J+1})$ is smaller than a user-prescribed accuracy; this ensures that the constructed subspace gives certified approximations. For the algorithm to work, it is crucial that the decay of the error estimate $\Delta_{\ISubspace{J}}$ along the iterations $J$ reflects the same behavior as the Kolmogorov $n$-width \eqref{eqn:KnW} over $\TrainSet$; see \cite{BinCDDPW11} for further discussions.

We need practical upper and lower bounds for the eigenvalues and spectral gap to obtain such error indicator functions. 
Hence, we next recall and introduce such bounds. 
After discussions on an upper bound for $\lambda_k(\prmtr)$ in \Cref{subsec:SUB}, we recall the subspace lower bound for $\lambda_1(\prmtr)$ from \cite{SirK16} in \Cref{subsec:SLB}. Extending the idea in \cite{SirK16}, we introduce a subspace lower bound for $\lambda_k(\prmtr)$ in \Cref{subsec:SLB k}.
Based on the bounds for eigenvalues, we introduce error bounds for $\gamma(\prmtr)$ in \Cref{subsec:gap bounds}, and we discuss in \Cref{sec:4.3} on conditions to ensure that the subspace approximation captures the dimensions of the exact eigenspaces.

\subsection{Subspace upper bound for \texorpdfstring{$\lambda_k(\prmtr)$}{TEXT}} 
\label{subsec:SUB}
By the variational characterization of eigenvalues for Hermitian matrices, which is also known as the min-max principle, we obtain 
\begin{align} \label{eqn:EigVal SUB}
\begin{aligned}
    \lambda_k(\prmtr) \,
    = \min_{\substack{\calL\subseteq \C^\fdim,\\ \dim \calL = k}} \max_{\vek \in \calL, \|\vek\|=1 } {\vek^* \Matrix(\prmtr)\vek} \;
    \leq \;\min_{\substack{\calL\subseteq \subspace,\\ \dim \calL = k}} \max_{\vek \in \calL,  \|\vek\|=1  } {\vek^* \Matrix(\prmtr)\vek}
    \;=\; 
    \RedEigVal{\subspace}_k (\prmtr),
    \end{aligned}
\end{align}
for $k=1,\dots, \dim (\subspace)$.
Therefore, we consider $\lambda_k^\SUB (\prmtr, \subspace) \vcentcolon= \RedEigVal{\subspace}_k (\prmtr)$.

Note that, if the eigenvalue $\lambda_{k}(\widehat{\prmtr})$ is simple at some $\widehat{\prmtr} \in {\prmtrSet}$,  
and if an eigenvector of $\Matrix(\widehat{\prmtr})$ associated to $\lambda_{k}(\widehat{\prmtr})$ is contained in the 
subspace ${\subspace}$ and serves as the eigenvector of the $k$-th eigenvalue of $\Matrix^{\mkern-6mu\basisV} (\widehat \prmtr)$, then $\RedEigVal{\subspace}_{k}(\prmtr)$ satisfies the Hermite interpolation property \cite[Lem.~6]{KanMMM18}, i.e. 
\begin{equation}\label{eq:interpolate}
	\lambda_{k}(\widehat{\prmtr})	
    \;	=	\;	
    \RedEigVal{\subspace}_{k}(\widehat{\prmtr})
	\quad\quad			\text{and}			\quad\quad
	\nabla \lambda_{k}(\widehat{\prmtr})	
    \;	=	\;	
    \nabla \RedEigVal{\subspace}_{k}(\widehat{\prmtr}) \, .
\end{equation}
Also note that the left-hand equality above, i.e. 
$\lambda_{k}(\widehat{\prmtr})=\RedEigVal{\subspace}_{k}(\widehat{\prmtr})$,
still holds true, even when $\lambda_{k}(\widehat{\prmtr})$ is not a simple eigenvalue of $\Matrix(\widehat{\prmtr})$, 
but as long as one eigenvector associated to $\lambda_{k}(\widehat{\prmtr})$ is contained in the subspace ${\subspace}$ and this eigenvector corresponds to the $k$-th eigenvalue of $\Matrix^{\mkern-6mu\basisV} (\widehat \prmtr)$. 

\subsection{Subspace lower bound for \texorpdfstring{$\lambda_1(\prmtr)$}{TEXT}} \label{subsec:SLB} 
In this section, we recall the derivation of a practically computable lower bound for $\lambda_1(\prmtr)$ in \cite{SirK16}. 
Based on the eigenvalue perturbation result in \cite{LiL05} combined with the \SCM method of \cite{HuyRSP07}, this lower bound is rather involved. 
Let us first review the theoretical perturbation results of \cite{LiL05}, and then recall how to obtain a practical lower bound of $\lambda_1(\prmtr)$ in \Cref{sec:lin:pro} following \cite{SirK16} . 

Suppose that there is an iterative procedure, which we will detail in \Cref{sec4}, returning $\prmtr_1, \dots , \prmtr_J \in {\prmtrSet}$ and a subspace $\ISubspace{J}$ as in \eqref{eq:low_bound} after $J$ iterations.  
Recall that $\IndBasisV{J}$ is a matrix whose columns form an \ONB for $\ISubspace{J}$. 

Let $\prmtr \in \prmtrSet$ be an arbitrary parameter. 
We denote by $\RedEigVal{\ISubspace{J}}_k(\prmtr)$ and $\vekRed{\IndBasisV{J}}_k(\prmtr)$ the $k$-th smallest eigenvalue counting multiplicity and a corresponding eigenvector of the compressed matrix
$\CompressedMat{\IndBasisV{J}} (\prmtr) = \IndBasisV{J}^\ast \, \Matrix(\prmtr) \IndBasisV{J}$, respectively.
For an integer $\indsub \leq \dim(\subspace)$ satisfying $\lambda^{\ISubspace{J}}_s(\prmtr) < \lambda^{\ISubspace{J}}_{s+1}(\prmtr) $ in case $ s <\dim(\subspace)$, we also define
\begin{equation}\label{eq:6}
	\IndBasisU{J}(\prmtr,s)	\, \vcentcolon= \,	
	\left[
	\begin{array}{ccc}
		\IndBasisV{J} \vekRed{\IndBasisV{J}}_1(\prmtr)	&	\dots		&		\IndBasisV{J} \vekRed{\IndBasisV{J}}_{\indsub}(\prmtr)
	\end{array}
	\right]		\quad	\text{and} \quad
	{\calU}_J (\prmtr,s)	\, \vcentcolon= \,	\mathrm{Col}(\IndBasisU{J} (\prmtr,s)). 
\end{equation} 
Note that ${\calU}_J (\prmtr,s)$ is independent of the choice of $\IndBasisV{J}$ due to the requirement $\RedEigVal{\ISubspace{J}}_s(\prmtr) < \RedEigVal{\ISubspace{J}}_{s+1}(\prmtr) $ for $ s <\dim(\subspace)$, since ${\calU}_J (\prmtr,s)$ is equal to $\subspace$ intersected with the eigenspaces of the matrix $\Proj^{\ISubspace{J}} \Matrix(\prmtr) \Proj^{\ISubspace{J}}$ associated with the eigenvalues $\RedEigVal{\ISubspace{J}}_1(\prmtr),\dots, \RedEigVal{\ISubspace{J}}_s (\prmtr)$. 
Also note that the number $s$ is not fixed and can vary depending on the object we are interested in approximating; we give more insight about this in \Cref{thm:exact:dim:rec}.

Let ${\calU}_J^{\bot}(\prmtr,s)$ be the orthogonal complement of ${\calU}_J(\prmtr,s)$, and $\IndBasisU{J}^{\bot}(\prmtr,s)$ be a matrix whose columns form an \ONB for ${\calU}_J^{\bot}(\prmtr,s)$. 
Since the matrix 
$\left[\begin{array}{cc}
			\IndBasisU{J}(\prmtr,s) & \IndBasisU{J}^{\bot}(\prmtr,s)
		\end{array}\right]$ 
is by construction unitary, the eigenvalues of $\Matrix(\prmtr)$ coincide with the eigenvalues of
\[
	\begin{bmatrix}
			\IndBasisU{J}(\prmtr,s)^\ast \Matrix(\prmtr) \IndBasisU{J}(\prmtr,s)	   &	 \IndBasisU{J}(\prmtr,s)^\ast \Matrix(\prmtr) \IndBasisU{J}^{\bot}(\prmtr,s)		\\[.2em]
			\IndBasisU{J}^{\bot}(\prmtr,s)^\ast \Matrix(\prmtr) \IndBasisU{J}(\prmtr,s)	   &	 \IndBasisU{J}^{\bot} (\prmtr,s)^\ast \Matrix(\prmtr) \IndBasisU{J}^{\bot}(\prmtr,s)	
		\end{bmatrix}.
\]
In particular, the smallest eigenvalue of the matrix above is also $\lambda_1(\prmtr)$. 
Disregarding the off-diagonal blocks of this matrix, we get the block-diagonal matrix 
\begin{align} \label{eq:off diagonal matrix}
    \Hat{\Matrix} (\prmtr,\indsub)
    \;\vcentcolon=\; 
    \begin{bmatrix}
			\IndBasisU{J}(\prmtr,s)^\ast \Matrix(\prmtr) \IndBasisU{J}(\prmtr,s)	   &	\zeroVec 		\\
			 	\zeroVec  &	 \IndBasisU{J}^{\bot} (\prmtr,s)^\ast \Matrix(\prmtr) \IndBasisU{J}^{\bot}(\prmtr,s)	
		\end{bmatrix},
\end{align}
whose eigenvalues are denoted by $\Hat{\lambda}_{{i}}(\prmtr)$ for ${{i}}=1,\ldots,\fdim$. 
Due to a classical eigenvalue perturbation result \cite[Sec.~10.3]{parlett98}, it holds for $1\leq i \leq \fdim$ that 
\begin{align}  \label{ineq:classical EigVal perturb}
    \lambda_i (\prmtr) \; \geq \; \hat{\lambda}_i(\prmtr) - \rho^{(J)}(\prmtr,s)
\end{align}
with $\rho^{(J)}(\prmtr,s) \vcentcolon= \| \basisU^{\bot}_J(\prmtr,s)^* \Matrix(\prmtr) \IndBasisU{J}(\prmtr,s) \|$. 
The lower bound
\eqref{ineq:classical EigVal perturb} has been improved in \cite[Thm.~2]{LiL05} and reads as 
\begin{align} \label{ineq:Li Li bound}
    \lambda_i (\prmtr) \; \geq \; \hat{\lambda}_i(\prmtr) - \frac{2\rho^{(J)}(\prmtr,s)^2}{ \delta_i(\prmtr) + \sqrt{\delta_i(\prmtr)^2 + 4 \rho^{(J)}(\prmtr,s)^2}} 
\end{align}
with 
\begin{align}
    \delta_i (\prmtr) \;\vcentcolon=\; 
    \begin{cases}
        \distance \left(\hat{\lambda}_i(\prmtr), \sigma(\Matrix^{\IndBasisU{J}^{\bot}(\prmtr,s)}(\prmtr))\right), & \text{ if}\quad \hat{\lambda}_i(\prmtr) \in  \sigma(\Matrix^{\IndBasisU{J}(\prmtr,s)}(\prmtr)), \\[.8em]
        \distance \left(\hat{\lambda}_i(\prmtr), \sigma(\Matrix^{\IndBasisU{J}(\prmtr,s)}(\prmtr))\right), & \text{ if}\quad \hat{\lambda}_i(\prmtr) \in  \sigma(\Matrix^{\IndBasisU{J}^{\perp}(\prmtr,s)}(\prmtr)),
    \end{cases} \label{eq:delta_j}
\end{align}
where $\distance(x,M)\vcentcolon= \inf_{m\in M} |m-x|$ for $x\in \R$ and $M \subseteq \R$. 
The term $\rho^{(J)}(\prmtr,s)^2$ can be efficiently calculated by solving an eigenvalue problem of a small matrix of size $\indsub \times \indsub$, namely
\begin{align}\label{eq:rhomu}
\begin{aligned}
	\rho^{(J)}(\prmtr,s)^2	
    \;	&	=	\;	
    \| (\IdMat_\fdim - \Proj^{\calU_J(\prmtr,s)})
    \Matrix(\prmtr) \IndBasisU{J}(\prmtr,s) \|^2	
     \; =\;		
     \| \Matrix(\prmtr) \IndBasisU{J}(\prmtr,s) - \IndBasisU{J}(\prmtr,s) \fLambda^{{\calU}_J}(\prmtr) \|^2		\\
    &	=	\;
	\lambda_{\max}\big(\IndBasisU{J}(\prmtr,s)^\ast \Matrix(\prmtr)^2  \IndBasisU{J}(\prmtr,s)  - \fLambda^{{\calU}_J}(\prmtr,s)^2 \big)	\,	,	
\end{aligned}
\end{align}
where $\lambda_{\max}(\cdot)$ denotes the largest eigenvalue of a Hermitian matrix, and 
\begin{align*}
    \fLambda^{{\calU}_J}(\prmtr,s)
		\;	\vcentcolon=	\;
	\IndBasisU{J}(\prmtr,s)^\ast \Matrix(\prmtr) \IndBasisU{J}(\prmtr,s)
		\;	=	\;
	\diag\left(\RedEigVal{\ISubspace{J}}_1(\prmtr),\ldots,\RedEigVal{\ISubspace{J}}_\indsub(\prmtr)\right)\,	.
\end{align*}
Also observe that $\rho^{(J)}(\prmtr,s)$ vanishes, whenever columns of $\IndBasisU{J}(\prmtr,s)$ are eigenvectors of $\Matrix(\prmtr)$, and so in particular at $\prmtr_1, \dots , \prmtr_J$, see \cite[Lem.~2.1]{ManMG24}.
 
An important observation from \cite{SirK16} is that 
\begin{equation*}
    \hat{\lambda}_1(\prmtr) = \min \{ \RedEigVal{\ISubspace{J}}_1(\prmtr), \lambda_1^{\calU_J^\perp(\prmtr,s)}(\prmtr) \}\quad\quad\text{and}\quad\quad \delta_1 (\prmtr) = | \RedEigVal{\ISubspace{J}}_1(\prmtr) - \lambda_1^{\calU_J^{\bot}(\prmtr,s)}(\prmtr) |.
\end{equation*} 
Hence, \eqref{ineq:Li Li bound} yields the following lower bound for the 
smallest eigenvalue $\lambda_1(\prmtr)$ of $\Matrix(\prmtr)$, i.e. 
\begin{align} \label{eq:LB by f of sth}
    f^{(J)}(\lambda^{{\mathcal U}_J^{\bot}(\prmtr,s)}_1(\prmtr))	\;	\leq		\;		\lambda_1(\prmtr)
    \end{align}
for 
\begin{align}\label{eq:low_bound1}
 f^{(J)}(\eta)
		\;	\vcentcolon=		\;
		\min \left\{ \RedEigVal{\ISubspace{J}}_1 (\prmtr) , \eta \right\}
		-
		\frac{2\rho^{(J)}(\prmtr,s)^2}{\left| \RedEigVal{\ISubspace{J}}_1(\prmtr) - \eta \right|  
		+   \sqrt{	\left| \RedEigVal{\ISubspace{J}}_1 (\prmtr) - \eta \right|^2   +   4 \rho^{(J)}(\prmtr,s)^2}}  \:.
\end{align}

From a practical perspective, the problem of the bound \eqref{eq:LB by f of sth} is that obtaining $ \lambda_1^{\calU_J^{\bot}(\prmtr,s)}(\prmtr)$ involves the calculation of the smallest eigenvalue of the large matrix $\IndBasisU{J}^{\bot} (\prmtr,s)^\ast \Matrix(\prmtr) \IndBasisU{J}^{\bot}(\prmtr,s)	$, which is nearly as expensive as the computation of $\lambda_1(\prmtr)$. 
However, the function $f^{(J)}(\eta)$ in (\ref{eq:low_bound1}) increases monotonically \cite[Lem.~3.1]{SirK16}.
Hence, any $\eta^{(J)}(\prmtr,s) \leq \lambda^{{\mathcal U}_J^{\bot}(\prmtr,s)}_1(\prmtr)$ that can be cheaply computed
yields a practical lower bound
\begin{equation}\label{eq:ineff_lb}
	f^{(J)}(\eta^{(J)}(\prmtr,s))	\;	\leq		\;		\lambda_1(\prmtr).
\end{equation}
In the following, we recall from \cite{SirK16} how an $\eta^{(J)}(\prmtr,s)$
satisfying $\eta^{(J)}(\prmtr,s) \leq \lambda^{{\mathcal U}_J^{\bot}(\prmtr,s)}_1(\prmtr)$
can be efficiently obtained. 

\subsubsection{Determination of an $\eta^{(J)}(\prmtr,s)$ such that 
    $\eta^{(J)}(\prmtr,s) \leq \lambda^{{\mathcal U}_J^{\bot}(\prmtr,s)}_1(\prmtr)$}
    \label{sec:lin:pro}
To simplify notations, let us denote $\lambda_k^{(j)} \coloneqq \lambda_k{(\prmtr_j)}$ and $\vek_k^{(j)} \coloneqq \vek_k{(\prmtr_j)}$ for $j=1,\dots,J$. 
As shown in \cite[Lem.~3.2]{SirK16}, under the assumption $\fdim \geq 2 s$, it holds for any $\prmtr \in \prmtrSet$ and any $j = 1, \dots , J$ that   
\begin{align} \label{eq:beta_i}
    \lambda^{{\mathcal U}_J^{\bot}(\prmtr,s)}_1(\prmtr_j)
		\;\;	\geq		\;\;
		\lambda_1^{(j)} + \beta^{(j,J)}(\prmtr,s)
\end{align}
with
\begin{align} 
\begin{aligned}\label{eq:def beta, Lambda}
        \beta^{(j,J)}(\prmtr,s)	\;	\vcentcolon=&	\;
		\lambda_{1}
		\left(
		(\fLambda^{(j)} - \lambda^{(j)}_1 \IdMat_{\ell(j)})
		-
		( \basisW^{(j)} )^\ast \IndBasisU{J}(\prmtr,s) \IndBasisU{J}(\prmtr,s)^\ast \basisW^{(j)} (\fLambda^{(j)} - \lambda^{(j)}_{{\ell(j)}+1} \IdMat_{\ell(j)})
		\right)	\,	, \\ 
  \fLambda^{(j)}\; \vcentcolon=&\;
		\diag\left(\lambda^{(j)}_1,\ldots,\lambda^{(j)}_{{\ell(j)}}\right),		\qquad 
		\basisW^{(j)}\;	\vcentcolon=\;
		\left[
		\begin{array}{ccc}
			\vek^{(j)}_1	&	\dots		&	\vek^{(j)}_{\ell(j)}
		\end{array}
		\right]. 
\end{aligned}
\end{align}
Notice that the term $\beta^{(j,J)}(\prmtr,s)$ is independent of the choice of the basis $\IndBasisV{J}$ and of the vectors $\vekRed{\IndBasisV{J}}_1(\prmtr) , \dots , \vekRed{\IndBasisV{J}}_s(\prmtr)$: The only term possibly depending on $\IndBasisV{J}$ or on $\vekRed{\IndBasisV{J}}_1(\prmtr) , \dots , \vekRed{\IndBasisV{J}}_s(\prmtr)$ in the definition of $\beta^{(j,J)}(\prmtr,s)$ would be $\IndBasisU{J}(\prmtr,s)\IndBasisU{J}(\prmtr,s)^*$. 
However, since $\IndBasisU{J}(\prmtr,s)\IndBasisU{J}(\prmtr,s)^* = \Proj^{\calU_J (\prmtr,s)}$ and we have chosen $\indsub \leq \dim(\subspace)$ such that $\RedEigVal{\ISubspace{J}}_s(\prmtr) < \RedEigVal{\ISubspace{J}}_{s+1}(\prmtr) $ in the case $ s <\dim(\subspace)$, the space $\calU(\prmtr,s)$ is independent of the choice of the base, so is its orthogonal projector.

With \eqref{eq:beta_i}, we can then derive a lower bound of $\lambda^{{\mathcal U}_J^{\bot}(\prmtr,s)}_1(\prmtr)$
using the \SCM approach \cite{HuyRSP07} based on linear programming:  
For an arbitrary parameter $\widehat{\prmtr} \in {\prmtrSet}$, we firstly observe that \Cref{general assumptions}~\ref{assume:param sep} enables writing $\lambda^{{\calU}_J^{\bot}(\prmtr,s)}_1(\widehat{\prmtr})$ as 
\begin{align*}
    \lambda^{{\calU}_J^{\bot}(\prmtr,s)}_1(\widehat{\prmtr}) 
    = \min_{\zVec \in {\mathbb C}^{\fdim - \indsub}, \|\zVec\|=1} \;
    \sum_{q = 1}^Q \theta_q(\widehat{\prmtr})  
    \left( \zVec^\ast \IndBasisU{J}^{\bot}(\prmtr,s)^\ast \Matrix_q  \IndBasisU{J}^{\bot}(\prmtr,s) \zVec \right).
\end{align*}
Then, by setting $\ftheta(\widehat{\prmtr}) \vcentcolon= \left[ \begin{array}{ccc} \theta_1(\widehat{\prmtr}) & \dots & \theta_Q(\widehat{\prmtr}) \end{array} \right] \in \R^Q$, 
\begin{align*}
    {\mathcal Q}_J^{(s)}(\prmtr) : \{ \zVec \in \C^{\fdim-s} | \, \|\zVec\|=1 \}	\rightarrow {\mathbb R}^Q	\,	, \ 
		\zVec \mapsto  
        \left\{ \zVec^\ast \IndBasisU{J}^{\bot}(\prmtr,s)^\ast \Matrix_q  \IndBasisU{J}^{\bot}(\prmtr,s) \zVec \right\}_{q=1}^Q
\end{align*}
as well as ${\mathcal Y}_J(\prmtr,s) \vcentcolon =  \{ {\mathcal Q}_J^{(s)}(\prmtr)(\zVec) \; | \; \zVec \in \C^{\fdim-s} ,  \|\zVec\|=1  \} \subseteq \R^Q $, 
we can express $\lambda^{{\mathcal U}_J^{\bot}(\prmtr,s)}_1(\widehat{\prmtr})$ as the solution of a minimization problem, i.e. 
\begin{align*}
    \lambda^{{\mathcal U}_J^{\bot}(\prmtr,s)}_1(\widehat{\prmtr}) = 
    \min_{\yVec \in {\mathcal Y}_J(\prmtr,s)}\ftheta(\widehat{\prmtr})^\T  \yVec. 
\end{align*}
A lower bound of $\lambda^{{\mathcal U}_J^{\bot}(\prmtr,s)}_1(\widehat{\prmtr})$ is then given by optimizing on a suitable set containing ${\mathcal Y}_J(\prmtr,s)$. 
To find such a set, we observe that at the chosen parameters $\prmtr_i$ for $i = 1, \dots , J$, we have
$\, \lambda^{{\mathcal U}_J^{\bot}(\prmtr,s)}_1(\prmtr_i) =
   \min_{\yVec \in {\mathcal Y}_J(\prmtr,s)}\ftheta(\prmtr_i)^\T  \yVec$.
Moreover, it holds for any $ \yVec \in {\mathcal Y}_J(\prmtr,s)$ that   
\begin{align*}
    \ftheta(\prmtr_i)^\T  \yVec
			\;	\geq		\;
		\lambda^{{\mathcal U}_J^{\bot}(\prmtr,s)}_1(\prmtr_i)
		\;	\geq		\;
		\lambda^{(i)}_1 + \beta^{(i,J)}(\prmtr,s),
\end{align*}
where the second inequality is due to (\ref{eq:beta_i}).
In addition, each component $y_q$ of $\yVec=[y_1 \; \dots \; y_Q] $ satisfies by definition
$y_q =  \zVec^\ast \IndBasisU{J}^{\bot}(\prmtr,s)^\ast \Matrix_q  \IndBasisU{J}^{\bot}(\prmtr,s) \,\zVec $
for some $\zVec\in {\mathbb C}^{\fdim - \indsub}$ with $ \|\zVec\|=1$, so
$y_q$ belongs to the interval $[\lambda_{1}(\Matrix_q), \lambda_{\fdim}(\Matrix_q)]$. 
It follows that any $\yVec \in {\mathcal Y}_J(\prmtr,s)$ 
also belongs to the set 
\begin{equation}\label{eq:LP_polytope}
{\mathcal Y}_{J}^{\rm{LB}}(\prmtr,s)
		\;	\vcentcolon=	\;
		\left\{
		\yVec	\in 	{\mathcal B}	\;	\middle|	\;	\ftheta(\prmtr_i)^\T  \yVec \: \geq \: \lambda^{(i)}_1 + \beta^{(i,J)}(\prmtr,s)
		\;	,	\;	i = 1, \dots, J
		\right\}, 
\end{equation}
where
\begin{align*}
    {\mathcal B}	
    \; \vcentcolon=	\;	
		[\lambda_{1}(\Matrix_1), \lambda_{\fdim}(\Matrix_1)]	\times	\dots	 \times [\lambda_{1}(\Matrix_Q), \lambda_{\fdim}(\Matrix_Q)] \subseteq \R^Q.
\end{align*}
Hence, ${\mathcal Y}_J(\prmtr,s)	\subseteq		{\mathcal Y}_{J}^{\rm{LB}}(\prmtr,s)$, which leads to the relation
\begin{equation}\label{eq:LP}
     \eta_{\ast}^{(J)}(\prmtr,s)	
     \;	\vcentcolon=	\;
     \min \{ \ftheta(\prmtr)^\T \yVec	\;	|	\;	\yVec \in {\mathcal Y}_{J}^{\rm{LB}}(\prmtr,s)	\} 
     \; \leq \;
     \lambda^{{\mathcal U}_J^{\bot}(\prmtr,s)}_1(\prmtr) \qquad  \text{ for all }\prmtr\in\prmtrSet.
\end{equation}
The minimization problem in \eqref{eq:LP} is a linear programming problem with $2Q+J$ constraints. 
Since its feasible set ${\mathcal Y}_{J}^{\rm{LB}}(\prmtr)$ is compact, problem \eqref{eq:LP}
must have a minimizer. Finally, using $\eta^{(J)}_\ast(\prmtr,s) \leq \lambda^{{\mathcal U}_J^{\bot}(\prmtr,s)}_1(\prmtr)$,
the monotonicity of $f^{(J)}$, and the inequality (\ref{eq:LB by f of sth}), one can define the subspace lower bound for the smallest eigenvalue as
\begin{equation}\label{eq:defn_LB}
	\lambda^{\SLB}_{1}(\prmtr, \ISubspace{J}) \; \vcentcolon= \;
	f^{(J)}(\eta^{(J)}_\ast(\prmtr,s))
	\;
	\leq		
	\;
	\lambda_1(\prmtr) .
\end{equation}
For a list of useful properties of $\eta^{(J)}_\ast(\prmtr,s)$ and $\beta^{(j,J)}(\prmtr,s)$, we refer to \cite[Lem.~2.3]{ManMG24}. 
Similarly to the upper bound, the lower bound $\lambda^{\SLB}_{1}(\prmtr, \ISubspace{J})$ in \eqref{eq:defn_LB} also interpolates the eigenvalue function $\lambda_{1}(\prmtr)$ at the selected parameter points $\prmtr_i$ in the Hermite interpolation sense for $i = 1,\dots , J$, provided that $\lambda_{1}(\prmtr_i)$ is simple and differentiable at $\prmtr_i$; otherwise, only the Lagrange interpolation holds. 
For further details, we refer to \cite[Thm.~3.6]{SirK16} or \cite[Thm.~2.4]{ManMG24}.

Regarding the computational expense of evaluating $\lambda^{\SLB}_{1}(\prmtr, \ISubspace{J})$, we refer to the comprehensive discussion in \cite[Sec.~3.3]{SirK16}. 
It is important to note that the process of evaluating $\lambda^{\SLB}_{1}(\prmtr, \ISubspace{J})$ incurs a computational cost that is independent of $\fdim$, and this cost is generally much smaller compared to that required for the evaluation of $\lambda_1(\prmtr)$.
We also refer to \cite[Alg.~1]{ManMG24} for a compact algorithmic description of the procedure to obtain $\lambda_1^\SLB(\prmtr,\ISubspace{J})$. 

Besides the Hermite interpolation property at selected parameter points, the involved lower bound $\lambda^{\SLB}_{1}(\prmtr, \ISubspace{J})$ has another advantage, namely it can be extended to a general subspace lower bound to higher eigenvalues, as presented in the next section. 

\subsection{Subspace lower bound for \texorpdfstring{$\lambda_k(\prmtr)$}{TEXT}} \label{subsec:SLB k}
We now introduce a lower bound $\lambda_k^\SLB(\prmtr,\ISubspace{J})$ for $\lambda_k(\prmtr)$ that holds for an arbitrary $k\le\dim(\subspace_J)$. 
To do so, we exploit the ideas used to derive the lower bound $\lambda_1^\SLB(\prmtr,\ISubspace{J})$, defined in \eqref{eq:defn_LB}, for $\lambda_1(\prmtr)$.

\begin{theorem}[Subspace lower bound for eigenvalues] \label{prop:EigVal SLB}
Consider the parameters $\prmtr_1,\dots, \prmtr_J$ and the subspace $\ISubspace{J}$ as in \eqref{eq:low_bound}. 
Let $ s \in \{1,\ldots,\dim(\ISubspace{J}) \}$ with $\RedEigVal{\ISubspace{J}}_s (\prmtr) < \RedEigVal{\ISubspace{J}}_{s+1}(\prmtr) $ in case $s<\dim(\ISubspace{J})$, and assume $\fdim \geq 2 s$.
Then, for any $\prmtr \in \prmtrSet$ and any $k\in \{1,\dots, s\}$, we have 
$$
\lambda_k^\SLB(\prmtr, \ISubspace{J}) \leq \lambda_k(\prmtr)
$$
with  {$\lambda_k^\SLB(\prmtr,\ISubspace{J}) \coloneqq \min\{\RedEigVal{\ISubspace{J}}_k (\prmtr), \eta^{(J)}_*(\prmtr,s) \}$ in case $\rho^{(J)}(\prmtr,s) = g_{k,\prmtr}(\eta^{(J)}_*(\prmtr,s)) = 0$, and otherwise}
\begin{equation}
    \lambda_k^\SLB(\prmtr,\ISubspace{J}) \;
    \vcentcolon= \;
    \min\{\RedEigVal{\ISubspace{J}}_k (\prmtr), \eta^{(J)}_*(\prmtr,s) \} 
    - \frac{2\rho^{(J)}(\prmtr,s)^2}{g_{k,\prmtr}(\eta^{(J)}_*(\prmtr,s)) + \sqrt{g_{k,\prmtr}(\eta^{(J)}_*(\prmtr,s))^2 + 4 \rho^{(J)}(\prmtr,s)^2 }}, 
    \label{eq:def SLB lambda_k}
\end{equation}
where 
\begin{itemize}
    \item $\rho^{(J)}(\prmtr,s)$ is defined as in \eqref{eq:rhomu};
    \item $\eta^{(J)}_{\ast}(\prmtr,s)$ is defined as in \eqref{eq:LP};
    \item $g_{k,\prmtr} : \R \rightarrow \R$ is defined as 
\begin{equation}\label{eqn:geta}
    g_{k,\prmtr}(\eta)\;\vcentcolon=\; \min\{ |\eta -\RedEigVal{\ISubspace{J}}_j  (\prmtr)| : 1\leq j \leq k \}.
\end{equation}
\end{itemize} 
\end{theorem}

\begin{remark} \label{rem:rho_is_0}
{The involved fraction term in \eqref{eq:def SLB lambda_k} can be viewed as a modified residual. 
It could happen that $\rho^{(J)}(\prmtr,s)=0$ and $g_{k, \prmtr}(\eta_*^{(J)} (\prmtr,s))=0$, leading to a formal “0/0'' problem in the fraction. 
However, $\rho^{(J)}(\prmtr,s)=0$ means that $\RedEigVal{\ISubspace{J}}_1(\prmtr), \dots, \RedEigVal{\ISubspace{J}}_s(\prmtr)$ are eigenvalues of $\Matrix(\prmtr)$; from  $g_{k,\prmtr}(\eta^{(J)}_*(\prmtr,s))=0$ we deduce that $\RedEigVal{\ISubspace{J}}_1(\prmtr) \leq \eta^{(J)}_*(\prmtr,s) \leq \RedEigVal{\ISubspace{J}}_k (\prmtr)$. 
Comparing with the definition of $\lambda_1^\SLB(\prmtr,\ISubspace{J})$ from \eqref{eq:defn_LB}, those conditions imply that $\lambda_1(\prmtr)=\RedEigVal{\ISubspace{J}}_1(\prmtr) $. 
Using again the definition of $\eta_*^{(J)}(\prmtr, s)$ and the bound \eqref{ineq:Li Li bound}, we also deduce that $\min\{\RedEigVal{\ISubspace{J}}_k (\prmtr), \eta^{(J)}_*(\prmtr,s) \} \leq \lambda_k(\prmtr)$. 
For this reason we can define $\lambda_k^\SLB(\prmtr,\ISubspace{J}) \coloneqq \min\{\RedEigVal{\ISubspace{J}}_k (\prmtr), \eta^{(J)}_*(\prmtr,s) \}$ for the case $\rho^{(J)}(\prmtr,s) = 0$ and $g_{k, \prmtr}(\eta_*^{(J)} (\prmtr,s))=0$ in the theorem above.
Let us highlight that the only potential issue in this configuration is that one of $\{\lambda_t(\prmtr)\}_{t=1}^k$ could have another associated eigenvector orthogonal to $\ISubspace{J}$. We discuss and handle this with \Cref{thm:exact:dim:rec}.}
\end{remark}

Moreover, we note that \Cref{prop:EigVal SLB} is a generalization of the lower bound in \Cref{subsec:SLB} from \cite{SirK16}, as in case $k=1$, \eqref{eq:def SLB lambda_k} is exactly the bound $\lambda_1^\SLB(\prmtr,\ISubspace{J})$ from \eqref{eq:defn_LB}.
Also note that $\rho^{(J)}(\prmtr,s)$ and $\eta_*^{(J)}(\prmtr,s)$ remains the same for all $k$ and that only $g_{k,\prmtr}$ depends on $k$.

\begin{proof}
We firstly recall from \Cref{sec:lin:pro} that $\eta^{(J)}_*(\prmtr,s) \leq \lambda_1^{\calU_J^\perp(\prmtr,s)}(\prmtr)$ and that $f^{(J)}(\eta^{(J)}_*(\prmtr,s)) \leq \lambda_1(\prmtr)$ for $f^{(J)}$ defined in \eqref{eq:low_bound1}. 
To prove \eqref{eq:def SLB lambda_k}, we need to discuss the four possible cases depending on the ordering of $\RedEigVal{\ISubspace{J}}_1(\prmtr) = \lambda_1^{\calU_J(\prmtr,s)}(\prmtr),$ $ \RedEigVal{\ISubspace{J}}_k(\prmtr) = \lambda_k^{\calU_J(\prmtr,s)}(\prmtr) $, $\lambda_1^{\calU_J^\perp(\prmtr,s)}(\prmtr)$ and $\eta^{(J)}_*(\prmtr,s)$. 
\begin{itemize}
\item 
Suppose $\eta^{(J)}_*(\prmtr,s) \leq \RedEigVal{\ISubspace{J}}_1 (\prmtr)$. 
In this case, we have $\min\{ \RedEigVal{\ISubspace{J}}_k(\prmtr), \eta^{(J)}_*(\prmtr,s) \} = \eta^{(J)}_*(\prmtr,s)$ as well as $g_{k,\prmtr} (\eta^{(J)}_*(\prmtr,s)) = |\eta^{(J)}_*(\prmtr,s)- \RedEigVal{\ISubspace{J}}_1(\prmtr)| = g_{1,\prmtr} (\eta^{(J)}_*(\prmtr,s))$. 
Therefore, we observe in this case that  
\begin{align*}
    \lambda_k^\SLB(\prmtr, \ISubspace{J}) \;=\; \lambda_1^\SLB(\prmtr, \ISubspace{J}) \;\leq\; \lambda_1(\prmtr) \le \lambda_k(\prmtr),
\end{align*} 
where $\lambda_1^\SLB(\prmtr, \ISubspace{J})$ is defined as in \eqref{eq:defn_LB}.

\item 
Suppose $\RedEigVal{\ISubspace{J}}_1  \leq \eta^{(J)}_*(\prmtr,s) \leq \RedEigVal{\ISubspace{J}}_k(\prmtr)  \leq \lambda_1^{\calU_J^\perp(\prmtr,s)}(\prmtr)$. 
In this case, we see that the $k$-th eigenvalue of the matrix $\hat\Matrix(\prmtr,\indsub)$ in \eqref{eq:off diagonal matrix} is $\RedEigVal{\ISubspace{J}}_k(\prmtr)$, i.e. $\hat \lambda_k(\prmtr)=\RedEigVal{\ISubspace{J}}_k(\prmtr)$, and so for $\delta_k(\prmtr) $ in \eqref{eq:delta_j}, we have $\delta_k(\prmtr) = |\lambda_1^{\calU_J^\perp(\prmtr,s)}(\prmtr) - \RedEigVal{\ISubspace{J}}_k(\prmtr)|$. 
Applying the perturbation theory result \cite[Thm.~2]{LiL05}, we see that 
\begin{align*}
   \lambda_k(\prmtr)
   \; & \geq \; 
    f_k\big( \lambda_1^{\calU_J^\perp(\prmtr,s)}(\prmtr) \big) 
\end{align*}
for $f_k:\R \to \R$ defined by  
$$
f_k(\eta) 
\coloneqq 
\min\{\RedEigVal{\ISubspace{J}}_k (\prmtr), \eta \} - 
\frac{2\rho^{(J)}(\prmtr,s)^2}{ | \eta - \RedEigVal{\ISubspace{J}}_k(\prmtr)| + \sqrt{ | \eta - \RedEigVal{\ISubspace{J}}_k(\prmtr)|^2 + 4 \rho^{(J)}(\prmtr,s)^2 }}.
$$
Furthermore, we deduce
\begin{equation*}
   \lambda_k(\prmtr)
   \; \geq \;f_k\big( \lambda_1^{\calU_J^\perp(\prmtr,s)}(\prmtr) \big)\;\geq\;f_k(\eta^{(J)}_*(\prmtr,s))\;\geq\;\lambda_k^\SLB(\prmtr, \ISubspace{J}),
\end{equation*}
where for the second inequality we use the fact that $f_k$, similar to $f^{(J)}$ of \eqref{eq:low_bound1} with $\RedEigVal{\ISubspace{J}}_1(\prmtr)$ replaced by $\RedEigVal{\ISubspace{J}}_k(\prmtr)$, is a monotonic increasing function (by the same argument as for $f^{(J)}$, see \cite[Lem. 3.1]{SirK16}), and for the last inequality, we observe
\[
   |\eta_*^{(J)}(\prmtr,s) - \RedEigVal{\ISubspace{J}}_{k}(\prmtr)| 
   \;\geq\; 
   \min\{ |\eta_*^{(J)}(\prmtr,s) -\RedEigVal{\ISubspace{J}}_j(\prmtr)| : 1\leq j \leq k \}.
\]
\item 
Suppose $\RedEigVal{\ISubspace{J}}_1 (\prmtr) \leq \eta_*^{(J)}(\prmtr,s) \leq \lambda_1^{\calU_J^\perp(\prmtr,s)}(\prmtr) \leq \RedEigVal{\ISubspace{J}}_k(\prmtr)$. 
In this case, let us define $\tau \coloneqq \min \big\{i\leq \dim(\ISubspace{J}) \, : \, \lambda_{i}^{\ISubspace{J}} (\prmtr) \geq \lambda_1^{\calU_J^\perp(\prmtr,s)}(\prmtr) \big\} $. 
Clearly $ k \geq \tau$, so $\lambda_k(\prmtr) \geq \lambda_{\tau}(\prmtr)$. 
Also observe that the $\tau$-th eigenvalue of $\hat\Matrix(\prmtr,\indsub)$ from \eqref{eq:off diagonal matrix} is $\lambda_1^{\calU_J^\perp(\prmtr,s)}(\prmtr)$, i.e. $\hat \lambda_\tau(\prmtr)=\lambda_1^{\calU_J^\perp(\prmtr,s)}(\prmtr)$, and so for $\delta_\tau(\prmtr)$ in \eqref{eq:delta_j}, it holds that $\delta_\tau(\prmtr) = g_{k,\prmtr}\big( \lambda_1^{\calU_J^\perp(\prmtr,s)}(\prmtr) \big)$. 
Then, \eqref{ineq:Li Li bound} yields    
\begin{align} 
    \lambda_k(\prmtr)
   \; \geq\; 
   \lambda_\tau(\prmtr)
   \; \geq\; 
   \lambda_1^{\calU_J^\perp(\prmtr,s)}(\prmtr) - \frac{2\rho^{(J)}(\prmtr,s)^2}{ \delta_\tau(\prmtr) + \sqrt{ \delta_\tau(\prmtr)^2 + 4 \rho^{(J)}(\prmtr,s)^2 }}. 
   \label{eq: lambda_k_SLB case 3}
\end{align}
By \Cref{lem:aux func monoton increasing ell bumps}, the function  
\begin{align}
    \hat{f}:\R \rightarrow \R,\, \eta \mapsto  \eta \;-\; \frac{2\rho^{(J)}(\prmtr,s)^2}{ g_{k,\prmtr} (\eta) + \sqrt{ g_{k,\prmtr} (\eta)^2 + 4 \rho^{(J)}(\prmtr,s)^2 }} \label{eq:aux func concrete k}
\end{align}
is monotonically increasing, and we see that the right most term of \eqref{eq: lambda_k_SLB case 3} is precisely $ \hat f\big(\lambda_1^{\calU_J^\perp(\prmtr,s)}(\prmtr) \big)$. 
Thus, we can further estimate \eqref{eq: lambda_k_SLB case 3} from below and obtain 
\begin{align*}
    \lambda_k(\prmtr)
  \; & \geq\; 
  \hat f \big( \eta_*^{(J)}(\prmtr,s) \big)
  \; =\; 
  \lambda_k^\SLB(\prmtr, \ISubspace{J}),
\end{align*}
where the inequality follows from $\eta_*^{(J)}(\prmtr,s) \leq \lambda_1^{\calU_J^\perp(\prmtr,s)}(\prmtr)$ as well as the monotonicity of $\hat f$,  and the equality is derived using the assumption $\eta_*^{(J)}(\prmtr,s) \leq \RedEigVal{\ISubspace{J}}_k(\prmtr)$ as well as the definition of $\lambda_k^\SLB(\prmtr, \ISubspace{J})$. 

\item 
Suppose $ \RedEigVal{\ISubspace{J}}_k(\prmtr) \leq \eta_*^{(J)}(\prmtr,s)$. 
In this case, we observe  $\hat\lambda_k(\prmtr) = \RedEigVal{\ISubspace{J}}_k(\prmtr) $ and $g_{k,\prmtr} \big(\eta_*^{(J)}(\prmtr,s)\big) = \eta_*^{(J)}(\prmtr,s) - \RedEigVal{\ISubspace{J}}_k(\prmtr) \leq \lambda_1^{\calU_J^\perp(\prmtr,s)}(\prmtr) - \RedEigVal{\ISubspace{J}}_k(\prmtr) = \delta_k(\prmtr)$ for $\delta_k(\prmtr)$ from \eqref{eq:delta_j}. 
Thus, $\lambda_k^\SLB(\prmtr, \ISubspace{J})$ is less than the lower bound of $\lambda_k(\prmtr)$ from \eqref{ineq:Li Li bound}, and  we deduce $\lambda_k^\SLB(\prmtr, \ISubspace{J}) \leq \lambda_k(\prmtr)$. 
\end{itemize}
In all cases, we obtain $\lambda_k^\SLB(\prmtr, \ISubspace{J}) \leq \lambda_k(\prmtr)$ and the proof is completed.
\end{proof}

Since we now have a practically computable lower bound for $\lambda_k(\prmtr)$, let us comment on a theoretical convergence result regarding $\lambda_k^{\SLB}(\prmtr,\ISubspace{J})$. 

\begin{remark} 
    Suppose that $\prmtrSet \subseteq \R$ is a closed interval and $\prmtr_1,\dots, \prmtr_J$ are the Chebyshev nodes of~$\prmtrSet$. 
    Then, \cite[Thm.~3.7]{SirK16} says that one can find $J_0 \in \N$ and some constants $M>0,R>1$  such that $|\lambda_1(\prmtr)-\lambda_1^\SLB(\prmtr, \ISubspace{J} )|\leq M R^{-2J}$ for all $J\geq J_0$, provided that $\lambda_1(\prmtr)$ is simple for all $\prmtr\in \prmtrSet$. 
    Besides, if $\lambda_k(\prmtr)$ is also always simple, then we can use the same argument to deduce that there exists some $\tilde J_0 \in \N$ and some $\tilde M>0, \tilde R >1 $ such that $|\lambda_k(\prmtr)-\lambda_k^\SLB(\prmtr, \ISubspace{J} )|\leq \tilde M \tilde R^{-2J}$ for all $J\geq \tilde J_0$. 
    However, we should note that in general $R, \tilde R$ could be very close to $1$, and so the convergence could be still slow; see \Cref{example:no crossing}. 
\end{remark}

We can also use $\lambda_k^\SLB(\prmtr,\ISubspace{J})$ for the subspace approximation of higher eigenvalues. 
\begin{remark}
Let $k\geq 2$. 
Since now there are certified subspace upper and lower bounds for $\lambda_k(\prmtr)$, namely $\lambda^\SUB_k(\prmtr,\ISubspace{J}) = \RedEigVal{\ISubspace{J}}_k(\prmtr)$ and $\lambda_k^\SLB(\prmtr,\ISubspace{J})$ from \Cref{prop:EigVal SLB}, one can use a greedy procedure to assemble a subspace for the approximation of $\lambda_k(\prmtr)$, e.g. using
\begin{align*}
{\lambda^\SUB_k(\prmtr,\ISubspace{J})-\lambda_k^\SLB(\prmtr,\ISubspace{J})} 
\end{align*}
as the error estimator, if one is only interested in the $k$-th eigenvalue. Alternatively one can use
\begin{align*}
\max_{1\leq j\leq k} \left\{
\lambda^\SUB_j(\prmtr,\ISubspace{J})-\lambda_j^\SLB(\prmtr,\ISubspace{J})
\right\}
\end{align*}
as the error estimator, if interested in a uniform approximation of all the smallest $k$ eigenvalues. 
In any case, we emphasize that at the newly selected parameter $\prmtr_{J+1}$, one should always calculate at least $\vek_1(\prmtr_{J+1}),\dots,\vek_k(\prmtr_{J+1})$ and append all $k$ vectors to the spanning set of the next subspace. 
In other words, it is not sufficient only to append $\vek_k(\prmtr_{J+1})$ for subspace approximations of the $k$-th eigenvalue, see \cite[Sec.~5.3]{BHP24} for further discussions. 
\end{remark}

We conclude this subsection by showing that the Hermite interpolation property holds for $\lambda_k^\SLB(\prmtr,\ISubspace{J})$. 

\begin{lemma}[Interpolation properties of $\lambda_k^\SLB(\prmtr,\ISubspace{J})$]\label{lem:int:EigVal:LB}
Let us take the same assumptions as in \Cref{prop:EigVal SLB}, and we consider $\lambda_k^\SLB(\prmtr)$ as defined in \eqref{eq:def SLB lambda_k}. 
Suppose at $\prmtr_i$ we have $\ell(i)\geq s$.
Then it holds for $k=1,\ldots,\indsub$ that 
\begin{equation}\label{eqn:int:EigVal:k:LB}
        \lambda_k^{\SLB}({\prmtr_i},{\ISubspace{J}}) 
        \;=\; 
        \lambda_k(\prmtr_i).
    \end{equation}
Moreover, if additionally $\lambda_k(\prmtr_i)$ is simple at $\prmtr_i$, i.e. $\lambda_{k-1}(\prmtr_i) < \lambda_{k}(\prmtr_i) < \lambda_{k+1}(\prmtr_i) $, then
\begin{equation}\label{eqn:Her:int:EigVal:k:LB}
    \nabla \lambda_k^{\SLB}({\prmtr_i},{\ISubspace{J}})
    \;=\;
    \nabla\lambda_k(\prmtr_i).
\end{equation}
\end{lemma}
\begin{proof}
Recall that by construction \eqref{eq:low_bound} we have $\vek_1(\prmtr_i),\dots, \vek_{\ell(i)}(\prmtr_i) \in \ISubspace{J}$. 
By \cite[Lem.~2.3]{ManMG24} and the condition $\ell(i) \geq \indsub \geq k$, we know that $\eta^{(J)}_{\ast}(\prmtr_i, \indsub) \ge \lambda_{\ell(i)+1}(\prmtr_i) \geq \lambda_{k+1}(\prmtr_i)$. 
 Then, we apply \cite[Lem.~2.1]{ManMG24} to obtain $\rho^{(J)}(\prmtr_i,\indsub)=0$. 
In addition, the interpolation properties of $\RedEigVal{\ISubspace{J}}_k(\prmtr_i)$ (see \cite[Lem.~2.6]{KanMMM18}) yield $\RedEigVal{\ISubspace{J}}_{k}(\prmtr_i) = \lambda_{k}(\prmtr_i)  $.
 Thus, we obtain by definition \eqref{eq:def SLB lambda_k} 
 \begin{equation*}
     \lambda^{\SLB}_{k}(\prmtr_i,\ISubspace{J}) 
     \;=\;
     \min\{\lambda_{k}(\prmtr_i),\lambda_{\ell(i)+1}(\prmtr_i)\}
     \;=\;
     \lambda_{k}(\prmtr_i).
 \end{equation*}

 For the interpolation of the derivative at $\prmtr_i$, let us assume $\lambda_{k-1}(\prmtr_i) < \lambda_{k}(\prmtr_i) < \lambda_{k+1}(\prmtr_i) $. 
 In this case, the previous considerations still hold, i.e., we have $\eta^{(J)}_\ast(\prmtr_i,\indsub) = \lambda_{\indsub+1}(\prmtr_i)$ and $\lambda_k^{\mathrm{SLB}}(\prmtr_i, \ISubspace{J}) = \RedEigVal{\ISubspace{J}}_k(\prmtr_i)$. 
As $\eta^{(J)}_\ast(\prmtr,\indsub)$ and $\RedEigVal{\ISubspace{J}}_k(\prmtr)$ 
are continuous functions of $\prmtr$, we deduce $\eta^{(J)}_\ast(\prmtr,\indsub) > \RedEigVal{\ISubspace{J}}_k(\prmtr)$ for all $\prmtr$ in a neighborhood of $\prmtr_i$, and therefore for such $\prmtr$ it holds
\begin{equation}\label{eq:LB_nearmui}
	\lambda^{\SLB}_{k}(\prmtr,\ISubspace{J})
    \;	=	\;
	\RedEigVal{\ISubspace{J}}_k(\prmtr)	-
	\frac{2\rho^{(J)}(\prmtr,\indsub)^2}{\left| \RedEigVal{\ISubspace{J}}_k(\prmtr) - \eta^{(J)}_\ast(\prmtr,\indsub) \right|  
		+   \sqrt{	\left| \RedEigVal{\ISubspace{J}}_k(\prmtr) - \eta^{(J)}_\ast(\prmtr,\indsub) \right|^2   +   4 \rho^{(J)}(\prmtr,\indsub)^2}}.
\end{equation}
By assumption $\lambda_{k}(\prmtr_i)$ is a simple eigenvalue of $\Matrix(\prmtr_i)$, so is 
$\RedEigVal{\ISubspace{J}}_k(\prmtr_i)$ as an eigenvalue of $\CompressedMat{\IndBasisV{J}} (\prmtr_i) = \IndBasisV{J}^\ast \Matrix(\prmtr_i) \IndBasisV{J}$, 
implying both $\lambda_{k}(\prmtr)$ and $\RedEigVal{\ISubspace{J}}_k(\prmtr)$ are differentiable at $\prmtr_i$.
In this case, as $\rho^{(J)}(\prmtr_i,\indsub) = 0$ due to \cite[Lem.~2.1]{ManMG24}, differentiating both sides of \eqref{eq:LB_nearmui} 
at $\prmtr = \prmtr_i$ gives rise to
$\nabla\lambda^{\SLB}_{k}(\prmtr_i,\ISubspace{J}) = \nabla \RedEigVal{\ISubspace{J}}_k(\prmtr_i)$.
By \cite[Lem.~2.6]{KanMMM18} we also have 
$\nabla \RedEigVal{\ISubspace{J}}_k(\prmtr_i) = \nabla \lambda_{k}(\prmtr_i)$. This concludes the proof.
\end{proof}

\subsection{Upper and lower bounds for \texorpdfstring{$\gamma(\prmtr)$}{TEXT}}
\label{subsec:gap bounds}
Taking advantage of the subspaces upper and lower bounds for eigenvalues, we introduce upper and lower bounds for the spectral gap as defined in \eqref{eq:gap full space}.
Once again suppose we have a subspace $\ISubspace{J} \subseteq \C^{\fdim}$ as in \eqref{eq:low_bound} with $ \dim(\ISubspace{J}) > m_1(\prmtr)$ for all $\prmtr\in\prmtrSet$. 
Note that $m_1(\prmtr)<m_0$ due to \Cref{general assumptions}~\ref{assume:separated lower spectrum}. 

If $m_1(\prmtr) \leq m_1(\prmtr,\ISubspace{J})$, then we can apply \eqref{eqn:EigVal SUB} as well as \eqref{eqn:gamma:LB} to derive 
\begin{align} \label{eq:gap SUB}
\gamma(\prmtr)
\,=\, 
\lambda_{m_1(\prmtr)+1}(\prmtr)-\lambda_1(\prmtr)
\,\leq\,
\RedEigVal{\ISubspace{J}}_{m_1(\prmtr,\ISubspace{J})+1}(\prmtr)- \lambda_1^{\SLB}(\prmtr,\ISubspace{J})
\,=\vcentcolon\,
\gamma^{\rm{SUB}}(\prmtr,{\ISubspace{J}})
\end{align}
for any $\prmtr\in\prmtrSet$. 
On the other hand, if $m_1(\prmtr) \geq m_1(\prmtr,\ISubspace{J})$, then \Cref{prop:EigVal SLB} and \eqref{eqn:EigVal SUB} yield
\begin{equation}\label{eqn:gamma:LB}
    \gamma(\prmtr) 
    \;=\;
    \lambda_{m_1(\prmtr)+1}(\prmtr)-\lambda_{1}(\prmtr)
    \;\geq\;
    \lambda^{\SLB}_{m_1(\prmtr,\subspace_J)+1}(\prmtr, \ISubspace{J})- \RedEigVal{\ISubspace{J}}_{1}(\prmtr)
    \;=\vcentcolon\;
    \gamma^{\SLB}(\prmtr,{\ISubspace{J}})
\end{equation}
for any $\prmtr\in\prmtrSet$. 
In other words, at least one of $\gamma^{\rm{SUB}}(\prmtr,{\subspace})$, $\gamma^{\SLB}(\prmtr,{\ISubspace{J}})$ is a rigorous upper or lower bound for the spectral gap. In particular, in the case $m_1(\prmtr) = m_1(\prmtr, \ISubspace{J})$, both $\gamma^{\rm{SUB}}(\prmtr,{\subspace})$ and $\gamma^{\SLB}(\prmtr,{\ISubspace{J}})$ become, respectively, a rigorous upper and lower bound for $\gamma(\prmtr)$. 
Hence, it is of interest to study under which conditions we can guarantee that $m_1(\prmtr) = m_1(\prmtr, \ISubspace{J}) $.  
This will be detailed in the next subsection. We conclude by mentioning that the Hermite interpolation properties hold for both $\gamma^{\SLB}(\prmtr,{\ISubspace{J}})$ and $\gamma^{\SUB}(\prmtr,{\ISubspace{J}})$ as a direct consequence of \Cref{lem:int:EigVal:LB} and \cite[Lem.~6]{KanMMM18}.

\subsection{Exact recovery of the eigenspaces dimension}\label{sec:4.3}
We begin with a lemma that provides an upper bound of $\eta_{\ast}^{(J)}(\prmtr,s)$. 
This result is needed to justify a theorem that specifies a computable sufficient condition to ensure the exact recovery of the dimensions of the eigenspaces. 

\begin{lemma}\label{lemma:eta} Let $k_0, t \in \N$ with $t \leq k_0$, $\prmtr\in\prmtrSet$, and let us define $\alpha(\prmtr, t) \vcentcolon= \sum_{i=1}^t m_i(\prmtr) $ as well as $\alpha(\prmtr, k_0) \vcentcolon= \sum_{i=1}^{k_0} m_i(\prmtr)$. 
Suppose a subspace $\ISubspace{J}$ is constructed w.r.t. parameters $\prmtr_1, \dots, \prmtr_J \in \TrainSet$ and $ \ell(j) = \alpha(\prmtr_j, k_0) $ for $j=1,\ldots,J$. 
Assume $s=\sum_{i=1}^t m_i(\prmtr , \ISubspace{J}) \le \dim(\ISubspace{J})$ and consider $\eta_{\ast}^{(J)}(\prmtr,s)$ as defined in \eqref{eq:LP}. 
Then it holds that 
    \begin{equation}\label{eqn:lemma:eta}
        \eta_{\ast}^{(J)}(\prmtr,s) 
        \;\le\;
        \lambda_{\alpha(\prmtr, t)+1}(\prmtr).
    \end{equation}
\end{lemma}
\begin{proof} 
By construction, see \Cref{sec:lin:pro}, it holds that $  \eta_{\ast}^{(J)}(\prmtr,s) \le \lambda^{{\mathcal U}_J^{\bot}(\prmtr,s)}_1(\prmtr) $.
{We recall that $\lambda^{{\mathcal U}_J^{\bot}(\prmtr,s)}_1(\prmtr)$ is the orthogonal complement of }$ \calU_J(\prmtr,s) = \sspan\{ \tilde \EigSpe_i^{\ISubspace{J}}(\prmtr) \}_{i=1}^t$, {where each $\tilde \EigSpe_i^{\ISubspace{J}}(\prmtr)$ is} the eigenspace, embedded in $\C^\fdim$,  of $\CompressedMat{\IndBasisV{J}} (\prmtr)$ associated to {eigenvalues $  \PropRedEigVal{\ISubspace{J}}_1 (\prmtr) < \dots < \PropRedEigVal{\ISubspace{J}}_t (\prmtr)$ counted without multiplicities}. 

If we can show that $\lambda^{{\calU}_J^{\bot}(\prmtr,s)}_1(\prmtr) \le \lambda_{\alpha(\prmtr, t)+1}(\prmtr)$, then \eqref{eqn:lemma:eta} follows immediately. 
Note that $\lambda_{\alpha(\prmtr, t)+1}(\prmtr) = \lambda_{\alpha(\prmtr, t+1)}(\prmtr) = \tilde \lambda_{t+1}(\prmtr) $, i.e. the $t$-th proper eigenvalue of $\Matrix(\prmtr)$ not counted in multiplicities; see \Cref{subsec:notation}. To proceed, let us discuss two exclusive cases. 
\begin{itemize}
\item 
Suppose $\calU^\perp_J(\prmtr,s) \perp \sspan\{ \tilde \EigSpe_i(\prmtr) \}_{i=1}^{t+1} $ with $\tilde \EigSpe_i(\prmtr)$ being the eigenspace associated to $\tilde \lambda_i (\prmtr)$ for $i=1,\dots t+1$. 
In other words, we have 
$$
\sspan\{ \tilde \EigSpe_1(\prmtr), \dots, \tilde \EigSpe_{t+1}(\prmtr) \} \subseteq \calU_J(\prmtr,s) = \sspan\{ \tilde \EigSpe_1^{\ISubspace{J}}(\prmtr), \dots, \tilde \EigSpe_t^{\ISubspace{J}}(\prmtr) \}.
$$ 
By the min-max principle and the interpolation property \cite[Lem~2.3]{KanMMM18}, we get
\begin{align*}
    \PropRedEigVal{\ISubspace{J}}_i = \tilde \lambda_i^{\calU_J(\prmtr,s)}(\prmtr) = \tilde \lambda_{i}(\prmtr),
    \quad \text{with}\quad m_i(\prmtr)=m_i(\prmtr,\calU_J(\prmtr,s))\quad\text{for}\quad i=1,\ldots,t+1. 
\end{align*}
However, by construction $\calU_J(\prmtr,s)$ only contains eigenspaces of $t$ different eigenvalues that are $\PropRedEigVal{\ISubspace{J}}_1 , \ldots, \PropRedEigVal{\ISubspace{J}}_t$. Therefore, the assumption $\calU^\perp_J(\prmtr,s) \perp \sspan\{ \tilde \EigSpe_i(\prmtr) \}_{i=1}^{t+1} $ leads to a contradiction. 
\item 
Now, consider $\calU^\perp_J(\prmtr,s) \not \perp \sspan\{ \tilde \EigSpe_i(\prmtr) \}_{i=1}^{t+1} $. 
Then, there exists $ \yVec_0 \in \sspan\{ \tilde \EigSpe_i(\prmtr) \}_{i=1}^{k+1}$ such that $\NullVector \neq \yVec = \Proj^{\calU_J^\perp(\prmtr,s)}\yVec_0 \in \calU_J^\perp(\prmtr,s)$, and we obtain by the min-max principle that 
\begin{equation}
\lambda_1^{ \calU_J^\perp(\prmtr,s) } (\prmtr) 
\; = 
\min_{ \uVec \in \calU_J^\perp(\prmtr,s) \backslash \{ \NullVector \}  } \frac {\uVec^* \Matrix(\prmtr) \uVec } {\uVec^*\uVec} \;
\leq \;
\frac{\yVec^* \Matrix(\prmtr) \yVec}  {\yVec^*\yVec} \;
\leq\; \lambda_{\alpha(\prmtr, t+1)}(\prmtr). 
\end{equation}
\end{itemize}
Therefore, $\lambda^{{\calU}_J^{\bot}(\prmtr,s)}_1(\prmtr) \le \lambda_{\alpha(\prmtr, t+1)}(\prmtr) = \lambda_{\alpha(\prmtr, t)+1}(\prmtr)$ holds. 
\end{proof} 
We are now ready to state a sufficient condition which ensures that the dimension of the eigenspaces is fully captured by the \ROM approximation. 

\begin{theorem} \label{thm:exact:dim:rec}
Let $\prmtr \in \prmtrSet$, $k_0\in\N$ such that $\sum_{i=1}^{k_0} m_i(\prmtr,\ISubspace{J}) \le \dim(\ISubspace{J})$, where $\ISubspace{J}$ is a subspace constructed w.r.t. parameters $\prmtr_1, \dots, \prmtr_J \in \TrainSet$ and $ \ell(j) = \alpha(\prmtr_j,k_0) =\sum_{i=1}^{k_0} m_i(\prmtr_j)  $ for $j=1,\ldots,J$.
For $1 \leq t \leq k_0$ we also set $s(t) \coloneqq \sum_{i=1}^t m_i(\prmtr,\ISubspace{J}) $ and $\alpha(\prmtr, t) = \sum_{i=1}^t m_i(\prmtr)$.
Suppose for a natural number $k \leq k_0$ we have
\begin{equation}\label{eqn:hyp}
    \eta_{\ast}^{(J)}(\prmtr,s(t))
    >
    \RedEigVal{\ISubspace{J}}_{s(t)}(\prmtr)+\varepsilon^{(J)}(\prmtr,s(t))
    \quad \text{ for all } t=1,\dots,k,
\end{equation} 
where {$\varepsilon^{(J)}(\prmtr,s(t)) \vcentcolon= 0$ in case $g_{s(t),\prmtr}(\eta^{(J)}_*(\prmtr,s(t))) = \rho^{(J)}(\prmtr,s(t)) = 0$, and otherwise } 
\begin{equation*}
 \varepsilon^{(J)}(\prmtr,s(t))
 \;\vcentcolon=\; 
 \frac{2\rho^{(J)}(\prmtr,s(t))^2}{g_{s(t),\prmtr}(\eta^{(J)}_*(\prmtr,s(t))) + \sqrt{g_{s(t),\prmtr}(\eta^{(J)}_*(\prmtr,s(t)))^2 + 4 \rho^{(J)}(\prmtr,s(t))^2 }} 
\end{equation*}
for $\rho^{(J)}(\prmtr,s(t))$ defined as in \eqref{eq:rhomu}, $\eta^{(J)}_{\ast}(\prmtr,s(t))$ defined as in \eqref{eq:LP}, and $g_{s(t),\prmtr}(\eta)$ defined as in \eqref{eqn:geta}. 
Then, it holds that $
m_t(\prmtr,\ISubspace{J}) =  m_t(\prmtr)
$ for all $1 \leq t \leq k$. 
\end{theorem}
{We emphasise that the key condition \eqref{eqn:hyp} is a \emph{strict} inequality to be checked. 
Recall from \Cref{rem:rho_is_0} that the special case 
\begin{equation*}
g_{s(t),\prmtr}(\eta^{(J)}_*(\prmtr,s(t)))=0\quad \text{and} \quad \rho^{(J)}(\prmtr,s(t))=0
\end{equation*}
indicates a possible mismatch between the eigenspace dimension of the \FOM and \ROM.
This is confirmed here: Indeed by $\varepsilon^{(J)}(\prmtr,s(t))=0$ and $\eta_{\ast}^{(J)}(\prmtr,s(t))
\leq
\RedEigVal{\ISubspace{J}}_{s(t)}(\prmtr)$ due to $g_{s(t),\prmtr}(\eta^{(J)}_*(\prmtr,s(t)))=0$, we see that \eqref{eqn:hyp} does not hold.} 

\begin{proof}
By \eqref{eq:LP} in \Cref{sec:lin:pro} it holds $\eta^{(J)}_{\ast}(\prmtr,s(t)) \le \lambda^{ \calU_J^{\perp}(\prmtr) }_1(\prmtr)$ with $ \calU_J(\prmtr, s(t)) = \sspan\{ \tilde \EigSpe_i^{\ISubspace{J}}(\prmtr) \}_{i=1}^t$ and $\tilde \EigSpe_i^{\ISubspace{J}}(\prmtr)$ is the eigenspace of $\CompressedMat{\IndBasisV{J}} (\prmtr)$ associated with $\PropRedEigVal{\ISubspace{J}}_i(\prmtr)$ and lifted to the full dimensional space{, see also the beginning of last lemma's proof for a reminder of the notations.}

To prove the desired equalities, let us begin with showing two auxiliary properties that hold for any $s(t)$ with $1\leq t \leq k$: 
\begin{enumerate}[label=(\roman*)]
\item\label{aux property: sep EigVal}
We have that $ \lambda_{s(t)}(\prmtr) < \lambda_{s(t)+1}(\prmtr) $: Indeed,
since $\eta^{(J)}_{\ast}(\prmtr,s(t))> \RedEigVal{\ISubspace{J}}_{s(t)}(\prmtr)$ by assumption 
 and  $\calU_J(\prmtr)$ is a subspace of dimension $s(t)$, we obtain from \cite[Thm.~2]{LiL05} and \Cref{lem:aux func monoton increasing ell bumps} that $\eta^{(J)}_{\ast}(\prmtr,s(t))-\varepsilon^{(J)}(\prmtr,s(t)) \le \lambda_{s(t)+1}(\prmtr)$. 
 Therefore we get
 \begin{equation*} 
     \lambda_{s(t)}(\prmtr)
     \;\le\; 
     \RedEigVal{\ISubspace{J}}_{s(t)}(\prmtr)
     \;<\;
     \eta^{(J)}_{\ast}(\prmtr,s(t))-\varepsilon^{(J)}(\prmtr, s(t))
     \;\le\;
     \lambda_{s(t)+1}(\prmtr). 
 \end{equation*} 
\item\label{aux property: alpha >= s}
We have that $ \alpha(\prmtr, t) \geq  s(t)$: 
Suppose by contradiction $ \alpha(\prmtr, t)  < s(t)$, then it holds that $\RedEigVal{\ISubspace{J}}_{\alpha(\prmtr, t)+1}(\prmtr) = \RedEigVal{\ISubspace{J}}_{\alpha(\prmtr, t) }(\prmtr)$.
Using the condition \eqref{eqn:hyp} and \Cref{lemma:eta}, we obtain
 \begin{equation*}
   \RedEigVal{\ISubspace{J}}_ {\alpha(\prmtr, t) +1}(\prmtr)
   \;=\;
   \RedEigVal{\ISubspace{J}}_{\alpha(\prmtr, t) }(\prmtr)
   \;<\;
   \eta_{\ast}^{(J)}(\prmtr,s(t)) 
  \;\le\;
   \lambda_{\alpha(\prmtr, t) +1}(\prmtr),
 \end{equation*}
 which contradicts the min-max principle. 
\end{enumerate}

Now, we prove the case $k = t = 1$, i.e., $s(1) = m_1(\prmtr, \ISubspace{J})=\alpha(\prmtr,1) = m_1(\prmtr)$. 
By \ref{aux property: sep EigVal}, \ref{aux property: alpha >= s} and the definition of $\lambda_{m_1(\prmtr)}(\prmtr)$, we obtain that $m_1(\prmtr) \leq s(1) = m_1(\prmtr,\ISubspace{J}) \leq \alpha(\prmtr,1) = m_1(\prmtr)$, and so $m_1(\prmtr) = m_1(\prmtr,\ISubspace{J})$.

The proof of the case $k > 1$ follows by a direct inductive argument on $t=1,\dots, k$. We show how the argument works for the case $k = 2$: 
Applying \ref{aux property: sep EigVal} and \ref{aux property: alpha >= s} to $t=1$, we see that $m_1(\prmtr) = m_1(\prmtr,\ISubspace{J})$, as presented before. 
Then, \ref{aux property: alpha >= s} for $t =2$ yields $m_1(\prmtr,\ISubspace{J}) + m_2(\prmtr,\ISubspace{J}) = s(2) \leq \alpha(\prmtr, 2) = m_1(\prmtr) + m_2(\prmtr) $. 
Using the facts $\lambda_{s(2)}(\prmtr) < \lambda_{s(2)+1}(\prmtr) $ by \ref{aux property: sep EigVal}, $m_1(\prmtr) = m_1(\prmtr,\ISubspace{J})$ as well as $m_2(\prmtr,\ISubspace{J})>1$, we deduce that $m_1(\prmtr) + m_2(\prmtr) \leq s(2) = m_1(\prmtr,\ISubspace{J}) + m_2(\prmtr,\ISubspace{J}) \leq m_1(\prmtr) + m_2(\prmtr) $; by $m_1(\prmtr) = m_1(\prmtr,\ISubspace{J})$ for $t=1$, we deduce that $m_2(\prmtr,\ISubspace{J}) = m_2(\prmtr)$ holds. 
\end{proof}

\section{A two stage greedy strategy} \label{sec4} 
We introduce a practical method to efficiently and reliably approximate $\lambda_1(\prmtr)$ and $\EigSpe_1(\prmtr)$. 
As discussed in \Cref{subsec:EigVal EigVec bounds}~and~\ref{subsec:bounds EigVec}, the approximation of $\EigSpe_1(\prmtr)$ requires knowledge of the spectral gap $\gamma(\prmtr)$.
Thus, we propose a two-stage approach: 
\begin{enumerate}
    \item 
    In the first stage, we construct a subspace $\subspace_{\gamma}$ by means of a greedy algorithm to approximate the spectral gap over the parametric domain, see \Cref{subsec:stage1}. 
    \item 
    In the second stage, we utilize $\subspace_{\gamma}$ to generate a different subspace $\subspace_{\EigSpe_1}$ also by means of a greedy algorithm. The subspace $\subspace_{\EigSpe_1}$ serves for the approximation of the eigenspace associated with the smallest eigenvalue in the parametric domain; see \Cref{subsec:stage2}.
\end{enumerate}

\subsection{Approximation of the spectral gap} \label{subsec:stage1}
Given a subspace $\subspace_{\gamma} \subseteq \C^\fdim$ and a matrix $\IndBasisV{\gamma}$  whose columns form an \ONB of $\subspace_{\gamma}$, we set 
\begin{align*}
    \gamma^{\subspace_{\gamma}} (\prmtr) 
    \;\vcentcolon=\; 
    \PropRedEigVal{\subspace_{\gamma}}_2(\prmtr) -   \RedEigVal{\subspace_{\gamma}}_1(\prmtr),
\end{align*}
where $\PropRedEigVal{\subspace_{\gamma}}_2$ is the smallest eigenvalue of $\CompressedMat{\IndBasisV{\gamma}} (\prmtr) \vcentcolon= \IndBasisV{\gamma}^*\Matrix(\prmtr)\IndBasisV{\gamma}$ such that $\tilde \lambda^{\subspace_\gamma}_2(\prmtr) \;>\; \RedEigVal{\subspace_{\gamma}}_1(\prmtr)$, i.e. $\gamma^{\subspace_{\gamma}}$ is the spectral gap of the matrix $\CompressedMat{\IndBasisV{\gamma}} (\prmtr) $.
We consider $\gamma^{\subspace_{\gamma}} (\prmtr) $ as the approximation of $\gamma(\prmtr)$ within the subspace $\subspace_\gamma$.
To assemble a suitable subspace $\subspace_\gamma$, we make use of the upper and lower bounds from \Cref{subsec:gap bounds}, 
and we define for an arbitrary subsapce $\subspace$ the error indicator function of the spectral gap by  
\begin{align}
    \Gamma^{\subspace }(\prmtr)\;\vcentcolon=\;\frac{\gamma^{\SUB}(\prmtr,\subspace)-\gamma^{\SLB}(\prmtr,\subspace )}{\gamma^{\subspace }(\prmtr)}. \label{eqn:spect:gap:err:est}
\end{align}
At the $J$-th iteration of the greedy algorithm, we then solve the maximization problem
\begin{equation} \label{eqn: max problem gap error}
    \max_{\prmtr\in\TrainSet}\, \Gamma^{\ISubspace{J}}(\prmtr)
\end{equation} 
and choose the next parameter $\prmtr_{J+1}$ as one satisfying \eqref{eqn: max problem gap error}. 

As mentioned in \Cref{subsec:gap bounds}, to provide a rigorous certified approximation of the spectral gap, the subspace $\subspace_{\gamma}$ provided by the greedy algorithm must satisfy $m_1(\prmtr,\subspace_{\gamma}) = m_1(\prmtr) $ for all $\prmtr\in\TrainSet$.
A sufficient condition for this equality comes from an application of \Cref{thm:exact:dim:rec} for $k=1$ (with $k_0=2$ for the assembly of $\subspace_\gamma$). 
Hence, we should check and enforce the condition \eqref{eqn:hyp} for $ k = 1 $ 
after the convergence of the greedy algorithm that employs \eqref{eqn:spect:gap:err:est} as the error estimator{, i.e., if the inequality \eqref{eqn:hyp} for $k=1$ is not satisfied at some $\prmtr \in \TrainSet$, we additionally select this point as an interpolation point}. 

Moreover, it is possible that  $\RedEigVal{\subspace}_{m_1(\prmtr)+1}(\prmtr)$ is not well-defined, i.e. $\CompressedMat{\basisV}  (\prmtr) = c \cdot  \IdMat_{\dim (\subspace)} $ for some $c\in \R$. 
In such cases, we should add $\prmtr$ as the next interpolation point, so that $\gamma(\prmtr)$ is captured by the subspace approximation.

In summary, we propose \Cref{alg:gap} for a certified reduced basis approximation of the spectral gap for a general parametric Hermitian matrix for which \Cref{general assumptions}~\ref{assume:param sep} holds.
\begin{algorithm}[t]
\begin{algorithmic}[1]
	\REQUIRE
Real analytic scalar functions $\theta_q(\prmtr): \R^Q\to \R$, Hermitian matrices $\Matrix_q \in {\mathbb C}^{\fdim \times \fdim}$ for $q=1,\dots,Q$ such that $\Matrix(\prmtr) = \sum_{q=1}^Q \theta_q(\prmtr)\Matrix_q$, a compact domain $\prmtrSet\subseteq {\R}^{\dimPar}$ and a set $\TrainSet \subseteq \prmtrSet$, a termination tolerance $\varepsilon_{\gamma}>0$.
	\ENSURE
    A reduced matrix-valued function  $\CompressedMat{\IndBasisV{\gamma}} (\prmtr)$ and the subspace ${\subspace_{\gamma}} = \mathrm{Col}(\IndBasisV{\gamma})$ such that 
	$\max_{\prmtr \in \TrainSet} \Gamma^{\subspace_{\gamma}}(\prmtr) \leq \varepsilon_{\gamma}$. 
	\STATE
    Compute $\lambda_{\min}(\Matrix_q)$ and $\lambda_{\max}(\Matrix_q)$ for $q = 1, \dots , Q$.
    \STATE 
    Set $J \gets 1$, choose an initial point $\prmtr_1 \in \TrainSet$, and set $S \gets \{\prmtr_1\}$. 
	\STATE
    Compute $\lambda_{k}(\prmtr_{1})$, $\vek_{k}(\prmtr_{1})$
	for $k = 1, \dots , \ell(1) \vcentcolon= m_1(\prmtr_{1})+m_2(\prmtr_{1})$, and $\lambda_{\ell(1)+1}(\prmtr_{1})$.
	\STATE 
    Set $\IndBasisV{1} \gets 
		\mathrm{orth}
		\left(
        \begin{bmatrix}
            \vek_1(\prmtr_{1})  &  \dots  & \vek_{\ell(1)}(\prmtr_{1})
        \end{bmatrix}
		\right)$ and 
    $\subspace_1 \gets \mathrm{Col}(\IndBasisV{1})$.
\WHILE{$\max_{\prmtr \in \TrainSet   } \Gamma^{\ISubspace{J}}(\prmtr) > \varepsilon_{\gamma}$}
\STATE Choose $\prmtr_{J+1} \gets \arg \max_{\prmtr \in \TrainSet   } \Gamma^{\ISubspace{J}}(\prmtr) $, and append $\prmtr_{J+1}$ to $S$. 
    \STATE \label{alg Gap: compute snapshot}
        Compute $ \left(\lambda_{k}(\prmtr_{J+1})  \vek_{k}(\prmtr_{J+1}) \right)_{k=1}^{\ell(J+1)}$ and $\lambda_{\ell(J+1)+1}(\prmtr_{J+1})$ for $ \ell(J+1)\vcentcolon=m_1(\prmtr_{J+1})+m_2(\prmtr_{J+1})$.
        \STATE \label{alg Gap: update space}
        Set $\IndBasisV{J+1} \gets \mathrm{orth}
		\left(
        \begin{bmatrix}
        \IndBasisV{J} & 
            \vek_1(\prmtr_{J+1})  &  \dots  & \vek_{\ell(J+1)}(\prmtr_{J+1})
        \end{bmatrix}
		\right)$ and 
    $\subspace_{J+1} \gets \mathrm{Col}(\IndBasisV{J+1})$. 
\STATE \label{alg Gap. update J} Set $J \gets J+1$. 
\ENDWHILE
    \WHILE{ $\TrainSet \backslash S \neq \emptyset $ }\label{alg gap: second loop}
    \STATE Choose $\prmtr \in \TrainSet \backslash S $, and append $\prmtr$ to $S$.
    \IF{condition \eqref{eqn:hyp} for $k = 1$ is \textit{not} satisfied \textbf{or} $\RedEigVal{\subspace}_{m_1(\prmtr)+1}(\prmtr)$ is \textit{not} well-defined  }
    \STATE Set $\prmtr_{J+1} \gets \prmtr$, and execute Step \ref{alg Gap: compute snapshot},\ref{alg Gap: update space} and \ref{alg Gap. update J}. 
 \ENDIF
    \ENDWHILE \label{alg gap: end second loop}
    \STATE \textbf{Terminate} with $\IndBasisV{\gamma}=\IndBasisV{J}$, $\CompressedMat{\IndBasisV{\gamma}} (\prmtr) = \IndBasisV{J}^\ast \Matrix(\prmtr) \IndBasisV{J}$ and $\ISubspace{\gamma}=\ISubspace{J}$.
\end{algorithmic}
\caption{A certified greedy strategy for the approximation of $\gamma(\prmtr)$ over $\TrainSet\subseteq\prmtrSet$}
\label{alg:gap}
\end{algorithm}

\subsection{Approximation of the eigenspace associated to the smallest eigenvalue} \label{subsec:stage2}
For given spaces $\subspace, \subspace_\gamma  \subseteq \C^\fdim$ and a number $\varepsilon_\gamma \in (0,1)$, {we define }
\begin{align}\label{eqn:new:err:est}
    \Delta_{\subspace_{\gamma}, \varepsilon_\gamma} (\prmtr, \subspace) 
   \; \vcentcolon= \;
    \frac{1}{ {\left(1-\varepsilon_{\gamma}\right)} \gamma^{ \subspace_{\gamma}}(\prmtr)}
    \left(
    \lambda_1^\SUB(\prmtr,\subspace) - \lambda_1^\SLB(\prmtr,\subspace) 
    \;+ \;
    \|  \ResSubsp{\subspace}(\prmtr)\|
    \right), 
\end{align}
where $\lambda_1^\SUB(\prmtr,\subspace)$ and $\lambda_1^\SLB(\prmtr,\subspace)$ are as introduced in \Cref{subsec:SUB}~and~\ref{subsec:SLB}, respectively.
{This expression is modified from  \Cref{thm:error bound EigVec} and the role of the new factor $(1-\varepsilon_\gamma)^{-1}$ will be clear later in \eqref{eqn:target}. }
The greedy algorithm for the eigenspace approximation is based on \eqref{eqn:new:err:est} as a posteriori error estimator, where the inputs $\subspace_{\gamma},\varepsilon_\gamma$ should result from \Cref{alg:gap} and the same set of parameters $\TrainSet \subseteq \prmtrSet$ must be used.
In the $J$-th step, the greedy algorithm solves the global maximization problem
\begin{equation}\label{eqn:eig:vec:err:est}
    \max_{\prmtr\in\TrainSet}\,  \Delta_{\subspace_{\gamma}, \varepsilon_\gamma} (\prmtr, \ISubspace{J}),
\end{equation}
and selects the maximizer of \eqref{eqn:eig:vec:err:est} as the next interpolation point. 
The following proposition clarifies why \eqref{eqn:new:err:est} is suitable as {a certified a posteriori error bound for the projection error.} 

\begin{proposition}\label{prop1}
Let $\subspace_{\gamma}$ be the subspace obtained from \Cref{alg:gap} with parameter set $\TrainSet \subseteq \prmtrSet$ and error accuracy $\varepsilon_\gamma \in (0,1)$. 
Then, for any subspace $\subspace$ and for any $\prmtr \in \TrainSet$, it holds that 
\begin{equation}\label{eqn:prop:stat}
       \| \Proj^{\EigSpe^{\bot}_1(\prmtr)} \Vek_1^{\basisV}(\prmtr) \|\;\le\; \Delta_{\subspace_{\gamma}, \varepsilon_\gamma} (\prmtr, \subspace),
\end{equation}
where $\basisV$ is matrix whose columns form an \ONB of $\subspace$, and $\Vek_1^{\basisV}(\prmtr)$ is a matrix whose columns form an \ONB of $\EigSpe_1^{\subspace}(\prmtr)$. 
\end{proposition}
\begin{proof}
    By \Cref{thm:error bound EigVec}, it is sufficient to show that $ \Delta_{\subspace_{\gamma}, \varepsilon_\gamma} (\prmtr, \subspace) $ is greater than the right-hand side of \eqref{eqn:bound:eig:vec}. 
    Since $\lambda_1^{\SUB}(\prmtr,\subspace)-\lambda_1(\prmtr)\le\lambda^{\SUB}_1(\prmtr,\subspace)-\lambda^{\SLB}_1(\prmtr,\subspace)$, we only need to show that \begin{equation}\label{eqn:target}
    \gamma(\prmtr)\;\ge\;\gamma^{\subspace_{\gamma}}(1-\varepsilon_\gamma).
    \end{equation}
Because of \Cref{alg:gap}, we have
\begin{align*}
    \begin{aligned}
        \gamma^{\SUB}(\prmtr,\subspace_{\gamma}) - \gamma(\prmtr)
        \;\le\;
        \gamma^{\SUB}(\prmtr,
    \subspace_{\gamma})-\gamma^{\SLB}(\prmtr,\subspace_{\gamma})\;\le\;\varepsilon_\gamma\left(\tilde \lambda^{\SUB}_2(\prmtr,\subspace_{\gamma})- \lambda^{\SUB}_1(\prmtr,\subspace_{\gamma})\right),
    \end{aligned}
\end{align*}
from where it follows that 
\begin{equation}\label{eqn:int}
    \gamma(\prmtr)\;\ge\; \lambda_2^{\SUB}(\prmtr,\subspace_{\gamma})(1-\varepsilon_\gamma)+\varepsilon_\gamma \lambda_1^{\SUB}(\prmtr,\subspace_{\gamma})-\lambda_1^{\SLB}(\prmtr,\subspace_{\gamma}).
\end{equation}
Observing that $-\lambda_1^{\SUB}(\prmtr,\subspace_{\gamma})\le-\lambda_1^{\SLB}(\prmtr,\subspace_{\gamma})$, and substituting it into \eqref{eqn:int}, we get \eqref{eqn:target}. 
\end{proof}
{By \Cref{theo:err:mes} and the discussions afterwards, in order to obtain a certified \ROM of the eigenspace in the sense of \eqref{eqn:goal},}
we need to make further adjustments to the greedy algorithm taking \eqref{eqn:new:err:est} as an error surrogate, since it could still happen that $\dim (\EigSpe_1^\subspace(\prmtr)) \neq \dim (\EigSpe_1(\prmtr)) $ for some $\prmtr \in \TrainSet$. 
Hence, we leverage \Cref{thm:exact:dim:rec} for $k=1$ and enforce condition \eqref{eqn:hyp} for all $\prmtr \in \TrainSet${, i.e., if the inequality \eqref{eqn:hyp} for $k=1$ is not satisfied at some $\prmtr \in \TrainSet$, we additionally select this point as an interpolation point}. 
Overall, we propose \Cref{alg:sf} for the certified reduced basis approximation of the eigenspace associated with the smallest eigenvalue. 
\begin{algorithm}[t]
\begin{algorithmic}[1]
	\REQUIRE
    The real analytic scalar functions $\theta_q(\prmtr): \R^Q\to \R$, Hermitian matrices $\Matrix_q \in {\mathbb C}^{\fdim \times \fdim}$ for $q=1,\dots,Q$ such that $\Matrix(\prmtr) = \sum_{q=1}^Q \theta_q(\prmtr)\Matrix_q$, compact domain $\prmtrSet\subseteq {\R}^{\dimPar}$ and a set $\TrainSet \subseteq \prmtrSet$ {from \Cref{alg:gap}},
    a termination tolerance $\varepsilon_{\EigSpe}>0$, the subspace $\subspace_{\gamma}$ and $\varepsilon_\gamma\in(0,1)$ provided by \Cref{alg:gap}.
	\ENSURE A reduced matrix-valued function  $\CompressedMat{\IndBasisV{\EigSpe}}(\prmtr)$ and the subspace ${\subspace_{\EigSpe}} = \mathrm{Col}(\IndBasisV{\EigSpe})$ s.t.
	$\max_{\prmtr \in \TrainSet}\Delta_{\subspace_{\gamma}, \varepsilon_\gamma}(\prmtr, \subspace_{\EigSpe})  \leq \varepsilon_{\EigSpe}$.  
	\STATE
    Compute $\lambda_{\min}(\Matrix_q)$ and $\lambda_{\max}(\Matrix_q)$ for $q = 1, \dots , Q$.
    \STATE 
    Set $J=1$, choose an initial point $\prmtr_1 \in \TrainSet$, and set $S \gets \{\prmtr_1\}$. 
	\STATE 
    Compute $\lambda_{k}(\prmtr_{1})$, $\vek_{k}(\prmtr_{1})$
	for $k = 1, \dots , \ell(1) \vcentcolon= m_1(\prmtr_{1})$, and $\lambda_{\ell(1)+1}(\prmtr_{1})$.
	\STATE 
    Set $\IndBasisV{1} \gets 
		\mathrm{orth}
		\left(
        \begin{bmatrix}
            \vek_1(\prmtr_{1})  &  \dots  & \vek_{\ell(1)}(\prmtr_{1})
        \end{bmatrix}
		\right)$ and 
    $\subspace_1 \gets \mathrm{Col}(\IndBasisV{1})$.
    \WHILE{$\max_{\prmtr \in \TrainSet } \Delta_{\subspace_{\gamma}, \varepsilon_\gamma}(\prmtr,\ISubspace{J}) > \varepsilon_\EigSpe$} \label{alg EigSpace: first loop}
    \STATE Choose $\prmtr_{J+1} \gets \arg \max_{\prmtr \in \TrainSet } \Delta_{\subspace_{\gamma}, \varepsilon_\gamma}(\prmtr,\ISubspace{J})$, and append $\prmtr_{J+1}$ to $S$.  
    \STATE \label{alg EigSpace: compute snap shot}
        Compute $\lambda_{k}(\prmtr_{J+1})$, $\vek_{k}(\prmtr_{J+1})$ for $k = 1, \dots , \ell(J+1) \vcentcolon= m_1(\prmtr_{J+1})$, and $\lambda_{\ell(J+1)+1}(\prmtr_{J+1})$.
        \STATE \label{alg EigSpace: update space}
        Set $\IndBasisV{J+1} \gets \mathrm{orth}
		\left(
        \begin{bmatrix}
        \IndBasisV{J} & 
            \vek_1(\prmtr_{J+1})  &  \dots  & \vek_{\ell(J+1)}(\prmtr_{J+1})
        \end{bmatrix}
		\right)$ and 
    $\subspace_{J+1} \gets \mathrm{Col}(\IndBasisV{J+1})$. 
    \STATE \label{alg EigSpace: update J}
    Set $J \gets J+1$
    \ENDWHILE \label{alg EigSpace: end first loop}
    \WHILE{ $\TrainSet \backslash S \neq \emptyset $ } \label{alg: second loop}
    \STATE Choose $\prmtr \in \TrainSet \backslash S $, and append $\prmtr$ to $S$.
    \IF{condition \eqref{eqn:hyp} for $k=1$ is \textit{not} satisfied }
    \STATE Set $\prmtr_{J+1} \gets \prmtr$, and execute Step \ref{alg EigSpace: compute snap shot}, \ref{alg EigSpace: update space} and \ref{alg EigSpace: update J}. 
    \ENDIF
    \ENDWHILE \label{alg: end second loop}
    \STATE \textbf{Terminate} with $\CompressedMat{\IndBasisV{\EigSpe}}(\prmtr) = \IndBasisV{J}^\ast \Matrix(\prmtr) \IndBasisV{J}$ and 
				$\subspace_{\EigSpe}=\ISubspace{J}$.
\end{algorithmic}
\caption{A certified greedy strategy for the approximation of $\EigSpe_1(\prmtr)$ over $\TrainSet\subseteq\prmtrSet$}
\label{alg:sf}
\end{algorithm}

{
\begin{remark}\label{rem:partial_certification}
By construction, \Cref{alg:sf} provides certified approximation of $\EigSpe_1^{\subspace_{\EigSpe}}(\prmtr)$ for $\prmtr\in \TrainSet$ in the sense of \eqref{eqn:goal}. 
When evaluating the resulting surrogate model  $\CompressedMat{\IndBasisV{\EigSpe}}(\prmtr)$ at $\prmtr\in \prmtrSet \backslash \TrainSet$ in the online phase, a \emph{conditional certification} of $\EigSpe_1^{\subspace_{\EigSpe}}(\prmtr)$ can also be obtained, provided that $\CompressedMat{\IndBasisV{\gamma}} (\prmtr)$, $\subspace_\gamma$ and $\varepsilon_\gamma$ from \Cref{alg:gap} is kept available in the online phase: 
If condition \eqref{eqn:hyp} for $k=1$ is satisfied at $\prmtr\in \prmtrSet \backslash \TrainSet$ w.r.t. \emph{both $\subspace_\gamma$ and $\subspace_\EigSpe$}, we know by \Cref{thm:exact:dim:rec} that $\dim\big(\EigSpe_1(\prmtr)\big) =  \dim\big(\EigSpe_1^{\subspace_{\gamma}}(\prmtr)\big) =  \dim\big(\EigSpe_1^{\subspace_{\EigSpe}}(\prmtr)\big)$. 
Suppose $\varepsilon_{\prmtr} \coloneqq \gamma^{\SUB}(\prmtr,\subspace_\gamma) - \gamma^{\SLB}(\prmtr,\subspace_\gamma) <1$, see \eqref{eq:gap SUB} and \eqref{eqn:gamma:LB}. 
Then, $\Delta_{\subspace_\gamma, \varepsilon_{\prmtr}}(\prmtr, \subspace_{\EigSpe}) $ from  \eqref{eqn:new:err:est} becomes a certified a-posteriori upper bound of $\|\Proj_{\EigSpe_1(\prmtr)} - \Proj_{\EigSpe_1^{\subspace_{\EigSpe}}(\prmtr)}\|$ due to \Cref{theo:err:mes} and  \Cref{prop1}. 
If, however, this condition is violated, a potential dimension mismatch may arise, in which case the user might consider performing an additional computation on the parameter $\prmtr$ at the \FOM level to update the subspaces $\subspace_\gamma, \subspace_\EigSpe$. 
\end{remark}
}

We also remark that there are parallels between the two-stage greedy selection strategies discussed here and the greedy method utilized in \RBM for parametric source problems:
As described in \cite[Sec.~3.2.2]{HesRS16}, developing a reduced space for parametric source problems can also employ a greedy strategy. 
Evaluating error estimates efficiently often requires approximating the smallest eigenvalue (or singular value) associated with the parametric matrix, which is usually achieved through another greedy algorithm (see \cite{HuyRSP07,SirK16,ManMG24}) as mentioned in the introduction. 
 In parametric eigenvalue problems, this process corresponds to \Cref{alg:gap} to approximate the spectral gap, while \Cref{alg:sf} is analogous to the weak-greedy algorithm for the parametric source problem.

\section{Numerical Experiments}\label{sec5}
This section presents numerical results to validate the two-stage greedy strategy. We consider both
randomly generated examples and parametrized \QSS. All computations are performed using \matlab~2023a on a laptop with an Apple M2 Pro processor.  {In cases where no different specification is given, the eigenvalues and corresponding eigenvectors of dense matrices are calculated using the \EIG function in \matlab. 
For sparse matrices, we utilize the \EIGS function with its default termination tolerance set at $10^{-14}$. 
To suitably upgrade the subspace $\calV$ accounting for finite machine precision errors, we adopt the following criteria to decide if two eigenvalues are coalescent: 
If both the eigenvalues are below $10^{-14}$ then we say that they are coalescent and equal to $0$; otherwise, we evaluate their relative distance and if this is below $10^{-8}$ we consider them as coalescent.}

For the numerical error analysis w.r.t. a given subspace $\subspace$, let us introduce the spectral gap error function
\begin{equation}\label{eqn:gap:err}
    \calE^{\subspace}(\prmtr)\;\vcentcolon=\;\frac{ |\gamma^{\subspace} (\prmtr)-\gamma(\prmtr)|}{\gamma^{\subspace} (\prmtr)}.
\end{equation}
We define the eigenvalue surrogate error as
\begin{equation}\label{eqn:eig:surr}
    H(\prmtr,\subspace)\;\vcentcolon=\;\lambda_1^{\SUB}(\prmtr,\subspace)-\lambda_1^{\SLB}(\prmtr,\subspace),
\end{equation}
and we note that, given $\subspace_\gamma$ and $\varepsilon_\gamma$, we can write
\begin{equation}\label{eqn:app:err:est}
\Delta_{\subspace_{\gamma}, \varepsilon_\gamma} (\prmtr, \subspace) 
\;\le\;
\frac{1}{(1-\varepsilon_{\gamma})} 
\left(
\frac{H(\prmtr,\subspace)} {\gamma^{\subspace_{\gamma}}(\prmtr)} 
\;+\; 
\frac{ \| \ResSubsp{\subspace} (\prmtr)\|}{ \gamma^{\subspace_{\gamma}}(\prmtr)}
\right),
\end{equation}
for $\Delta_{\subspace_{\gamma}, \varepsilon_\gamma} (\prmtr, \subspace)$ defined in \eqref{eqn:new:err:est}. 
Notice that the multiplication by the constant $(1-\varepsilon_{\gamma})^{-1}$ on the right-hand side of \eqref{eqn:app:err:est} does not affect the parameter points selected by \Cref{alg:sf}. 
Moreover, if $\varepsilon_{\gamma}$ is small enough, say $\varepsilon_{\gamma}\;\le\;10^{-2}$, then we can consider $(1-\varepsilon_{\gamma})^{-1} \approx 1$.

Usually greedy algorithms are performed w.r.t. a discrete parameter subset $\TrainSet \subseteq \prmtrSet$ \cite{HuyRSP07,SirK16}.  
However, when the dimension of the parameter set $\dimPar$ for $\prmtrSet \subseteq \R^\dimPar$ is small, such as $\dimPar = 1$ or $\dimPar = 2$, and when the optimized function is differentiable in the domain, it is also possible to perform greedy algorithms w.r.t. the entire set ${\prmtrSet}$ \cite{KanMMM18,ManMG24}. 
In such situations, finding a global solution to the optimization problems described in \eqref{eqn:spect:gap:err:est} and \eqref{eqn:eig:vec:err:est} for $\TrainSet=\prmtrSet$ can be performed through
\EIGOPT~\cite{MenYK14}, which takes advantage of the Lipschitz continuity in the eigenvalue functions with respect to the parameters. 
On the contrary, if $\dimPar>2$, or if the function to be optimized lacks differentiability at many points, it is preferable to use a discrete subset $\TrainSet$ of $\prmtrSet$. 
In such cases, the target function is evaluated at each point in $\TrainSet$ to solve the optimization problem. 
In this section, we will also test the applicability of $\EIGOPT$ for the first parametric eigenvalue problem, where the eigenspace of the smallest eigenvalue has always dimension $1$, see \Cref{subsec:toy example}. 

\vspace{0.2cm}
\noindent\fbox{%
    \parbox{0.98\textwidth}{%
        The code and data used to generate the subsequent results are accessible via
		\begin{center}
			\url{https://doi.org/10.5281/zenodo.15106369}
		\end{center}
		under MIT Common License.
    }%
}\\[.2em]

\subsection{Randomly generated Hermitian matrices} \label{subsec:toy example}
We consider a parametric matrix of the form
\begin{equation}\label{eqn:ex2:aff}
	\Matrix(\prmtrs) = \prmtrs^2 \Matrix_1 + \prmtrs \, \Matrix_2, \quad\quad \prmtrs \in \prmtrSet\vcentcolon=[10^{-1}, 10] ,
\end{equation}
where $\Matrix_1$, $\Matrix_2 \in \R^{\fdim \times \fdim}$ with $\fdim=2\times 10^3$ are randomly generated dense Hermitian matrices. 
The matrices $\Matrix_1, \Matrix_2$ are selected such that $\Matrix({\prmtrs})$ has no eigenvalue crossing in the parameter domain, and so we can apply \EIGOPT for the optimization for this example. 
Our objective is to find a subspace $\subspace_{\EigSpe}$ such that the eigenpair $(\lambda^{\subspace}_{1}(\prmtrs), \EigSpe_1^\subspace(\prmtrs))$ of the projected problem approximates the eigenpair $(\lambda_{1}(\prmtrs), \EigSpe_1(\prmtrs))$ of $\Matrix(\prmtrs)$ with an error not larger than $\varepsilon_{\EigSpe} = 10^{-8}$ in the continuum domain $\prmtrSet$ employing \EIGOPT.
Regarding the spectral gap error function $\calE^{\subspace_{\gamma}}(\prmtrs)$ defined in \eqref{eqn:gap:err}, we set the exit tolerance of \Cref{alg:gap} to $\varepsilon_\gamma=10^{-8}$. 
We first address the performance of the greedy strategy for approximating the spectral gap through \Cref{fig1}.
\begin{figure}[t]	
	\centering{
	\begin{subfigure}[t]{0.32\textwidth}
%
\tikzexternaldisable
\begin{tikzpicture}

\begin{axis}[%
	width=0.7*\imageWidth,
	height=\imageHeight,
	scale only axis,
	scaled ticks=false,
	grid=both,
	grid style={line width=.1pt, draw=gray!10},
	major grid style={line width=.2pt,draw=gray!50},
	axis lines*=left,
	axis line style={line width=\lineWidth},
xmode=log,
xmin=1e-1,
xmax=10,
xlabel style={font=\color{white!15!black}},
xlabel={$\mu$},,
ymode=log,
ymin=1e-2,
ymax=1e2,
yminorticks=true,
ylabel style={font=\color{white!15!black}},
	axis background/.style={fill=white},
	legend style={%
		legend cell align=left, 
		align=left, 
		font=\tiny,
		draw=white!15!black,
		at={(0.70,0.02)},
		anchor=south,},
]

\addplot [color = mycolor1, line width=\lineWidth] table [x=x1, y=y1, col sep=comma]{Img_Tex/DataCSV/Example_1_1.csv};
    \addlegendentry{$ \gamma(\prmtrs)$}
    
    \addplot [color = mycolor2,dashed, line width=\lineWidth] table [x=x1, y=y2, col sep=comma]{Img_Tex/DataCSV/Example_1_1.csv};
    \addlegendentry{$\gamma^{\subspace_\gamma}(\prmtrs)$}

        \addplot [color = mycolor3, line width=\lineWidth,only marks, mark=asterisk, mark options={solid, mycolor3}] table [x=x1, y=y1, col sep=comma]{Img_Tex/DataCSV/Example_1_1_bis.csv};
    \addlegendentry{$(\prmtrs_j,\gamma(\prmtrs_j)))$}

\end{axis}
\end{tikzpicture}%
     \subcaption{Plots of $\gamma(\prmtrs)$ and its approximation $\gamma^{\subspace_{\gamma}}(\prmtrs)$ over $\prmtrSet$. The interpolation points $(\prmtrs_j,\gamma(\prmtrs_j))$ selected by \Cref{alg:gap} are also displayed.}
     \label{fig1:a}
	\end{subfigure}
	\hspace{0.01cm}
 \begin{subfigure}[t]{0.32\textwidth}
%
\tikzexternaldisable
\begin{tikzpicture}
\begin{axis}[%
	width=0.7*\imageWidth,
	height=\imageHeight,
	scale only axis,
	scaled ticks=false,
	grid=both,
	grid style={line width=.1pt, draw=gray!10},
	major grid style={line width=.2pt,draw=gray!50},
	axis lines*=left,
	axis line style={line width=\lineWidth},
xmode=log,
xmin=1e-1,
xmax=10,
xlabel style={font=\color{white!15!black}},
xlabel={$\mu$},
ymode=log,
ymin=1e-14,
ymax=1e-10,
yminorticks=true,
ylabel style={font=\color{white!15!black}},
	axis background/.style={fill=white},
	legend style={%
		legend cell align=left, 
		align=left, 
		font=\tiny,
		draw=white!15!black,
		at={(0.90,0.98)},
		anchor=north east,},
]

\addplot [color = mycolor1, line width=\lineWidth] table [x=x1, y=y3, col sep=comma]{Img_Tex/DataCSV/Example_1_1.csv};
\addlegendentry{$\calE^{\subspace_\gamma}(\prmtrs)$ (see \eqref{eqn:gap:err})}
\end{axis}
\end{tikzpicture}%
     \subcaption{Spectral gap approximation error over $\prmtrSet$.}
     \label{fig1:b}
	\end{subfigure}
    \hspace{0.01cm}
 \begin{subfigure}[t]{0.32\textwidth}
\tikzexternaldisable
\begin{tikzpicture}

\begin{axis}[%
	width=0.7*\imageWidth,
	height=\imageHeight,
	scale only axis,
	scaled ticks=false,
	grid=both,
	grid style={line width=.1pt, draw=gray!10},
	major grid style={line width=.2pt,draw=gray!50},
	axis lines*=left,
	axis line style={line width=\lineWidth},
xmin=0,
xmax=32,
xlabel style={font=\color{white!15!black}},
xlabel={$j$},
ymode=log,
ymin=1e-9,
ymax=1e5,
yminorticks=true,
ylabel style={font=\color{white!15!black}},
	axis background/.style={fill=white},
	legend style={%
		legend cell align=left, 
		align=left, 
		font=\tiny,
		draw=white!15!black,
		at={(1.05,0.99)},
		anchor=north east,},
]

\addplot [color=mycolor1, line width=\lineWidth, mark=o, mark options={solid, mycolor1}]
  table [x=x1, y=y1, col sep=comma]{Img_Tex/DataCSV/Example_1_2.csv};

\addlegendentry{$\Gamma^{\subspace_{j-1}}(\prmtrs_j)$ (see \eqref{eqn:spect:gap:err:est})}
\end{axis}
\end{tikzpicture}%
     \subcaption{Decay of $\max_{\prmtrs\in\prmtrSet}\Gamma^{\subspace_{j-1}}(\prmtrs) $ with respect to the iteration counter $j$ of \Cref{alg:gap}.}
     \label{fig1:c}
	\end{subfigure}}
  \caption{Approximation of the spectral gap of a dense parametric Hermitian matrix $\Matrix(\mu)\in\R^{\fdim\times \fdim}$ as in \eqref{eqn:ex2:aff} with $\fdim=2\cdot 10^3$; the resulting reduced space $\subspace_{\gamma}$ for the gap approximation has dimension $\rdim={62}$.}
	\label{fig1}
\end{figure}

\Cref{fig1:a} illustrates the spectral gap of the full and reduced problems together with the selected interpolation points. \Cref{fig1:b} shows that, as expected, the approximation error is always below the prescribed accuracy $\varepsilon_\gamma=10^{-8}$. 
It can be observed that a subspace with dimension $\rdim={62}$ effectively approximates the spectral gap of a problem with dimension $\fdim=2\times10^3$, achieving a precision of at least $8$ decimal digits. {The average computational time to evaluate the spectral gap for the \FOM is 0.43 second while for the \ROM is $7.7\cdot10^{-4}$, therefore achieving a speed-up around $558$ times.}
Finally, \Cref{fig1:c} refers to the decay of the surrogate error with respect to the iterations of \Cref{alg:gap}. 
\begin{figure}[t]
	\centering{
	  \begin{subfigure}[t]{.48\linewidth}
%
\tikzexternaldisable
\begin{tikzpicture}

\begin{axis}[%
	width=\imageWidth,
	height=\imageHeight,
	scale only axis,
	scaled ticks=false,
	grid=both,
	grid style={line width=.1pt, draw=gray!10},
	major grid style={line width=.2pt,draw=gray!50},
	axis lines*=left,
	axis line style={line width=\lineWidth},
xmode=log,
xmin=1e-1,
xmax=10,
xlabel style={font=\color{white!15!black}},
xlabel={$\mu$},
ymode=log,
ymin=1e-14,
ymax=1e-8,
yminorticks=true,
ylabel style={font=\color{white!15!black}},
	axis background/.style={fill=white},
	legend style={%
		legend cell align=left, 
		align=left, 
		font=\tiny,
		draw=white!15!black,
		at={(0.28,0.02)},
		anchor=south west,},
]

\addplot [color=mycolor1, line width=\lineWidth]
  table [x=x1, y=y1, col sep=comma]{Img_Tex/DataCSV/Example_1_4.csv};
  
  \addlegendentry{$ \| \Proj^{\EigSpe^{\bot}_1(\prmtrs)} \Vek_1^\basisV(\prmtrs) \|$}
  
  \addplot [color=mycolor2, dashed, line width=\lineWidth]
  table [x=x1, y=y2, col sep=comma]{Img_Tex/DataCSV/Example_1_4.csv};
  \addlegendentry{$\Delta_{\subspace_\gamma,\varepsilon_\gamma}(\prmtrs,\subspace_{\EigSpe})$}

\end{axis}
\end{tikzpicture}%
    		\subcaption{Approximation error and error estimate \eqref{eqn:new:err:est} over $\prmtrSet$.}
   		\label{fig2:a}
   	\end{subfigure}
    \hfill
	\begin{subfigure}[t]{.48\linewidth}
\tikzexternaldisable
\begin{tikzpicture}

\begin{axis}[%
	width=\imageWidth,
	height=\imageHeight,
	scale only axis,
	scaled ticks=false,
	grid=both,
	grid style={line width=.1pt, draw=gray!10},
	major grid style={line width=.2pt,draw=gray!50},
	axis lines*=left,
	axis line style={line width=\lineWidth},
xmin=0,
xmax=40,
xlabel style={font=\color{white!15!black}},
xlabel={$j$},
ymode=log,
ymin=1e-17,
ymax=5e13,
yminorticks=true,
ylabel style={font=\color{white!15!black}},
	axis background/.style={fill=white},
	legend style={%
		legend cell align=left, 
		align=left, 
		font=\tiny,
		draw=white!15!black,
		at={(1.0,1.0)},
		anchor=north east,},
]

\addplot [color=mycolor1, line width=1.5*\lineWidth]
 table [x=x1, y=y1, col sep=comma]{Img_Tex/DataCSV/Example_1_3.csv};
 \addlegendentry{$\Delta_{\subspace_\gamma,\varepsilon_\gamma}(\prmtrs_{j},\subspace_{j-1})$ (see \eqref{eqn:new:err:est})}
 
 \addplot [color=mycolor2, dashed, line width=1.5*\lineWidth]
 table [x=x1, y=y2, col sep=comma]{Img_Tex/DataCSV/Example_1_3.csv};
 
 \addlegendentry{$\|\Res^{\subspace_{j-1}}(\prmtrs_{j})\|/\gamma^{\subspace_\gamma}(\prmtrs_{j})$}

 \addplot [color=mycolor3, dashdotted, line width=1.5*\lineWidth]
 table [x=x1, y=y3, col sep=comma]{Img_Tex/DataCSV/Example_1_3.csv};

\addlegendentry{$H(\prmtrs_{j},\subspace_{j-1})/\gamma^{\subspace_\gamma}(\prmtrs_{j})$ (see \eqref{eqn:eig:surr})}

\addplot [color=mycolor4, dotted, line width=1.5*\lineWidth]
 table [x=x1, y=y4, col sep=comma]{Img_Tex/DataCSV/Example_1_3.csv};

\addlegendentry{$\| \Proj^{\EigSpe_1(\mu_j)} -\Proj^{\EigSpe^{\subspace_{j-1}}_1(\mu_j)}
    \|$}
\end{axis}
\end{tikzpicture}%
   		\subcaption{Decay of the estimator \eqref{eqn:new:err:est} with its residual and eigenvalue error components.}
   		\label{fig2:b}
   	\end{subfigure}
	}
	
	\caption{Approximation of the eigenspace associated to the smallest eigenvalue of a dense parametric Hermitian matrix $\Matrix(\mu)\in\R^{\fdim\times \fdim}$ as in \eqref{eqn:ex2:aff} with $\fdim=2\cdot 10^3$; the resulting reduced space $\subspace_{\EigSpe}$ for the eigenspace approximation has dimension $\rdim=39$.}
	\label{fig2}
\end{figure}

We now consider the subspace $\subspace_{\gamma}$ returned by \Cref{alg:gap} for $\varepsilon_{\gamma}=10^{-8}$, set $\varepsilon_{\EigSpe}=10^{-8}$ and run \Cref{alg:sf}. 
Since the spectral gap of $\Matrix(\prmtrs)$ is well defined and positive, as we can see from \Cref{fig1:a}, $\vek_1(\prmtrs)$ is differentiable. 
Thus, we still use \EIGOPT for the approximation of the eigenspace in this example, and we obtain a reduced space $\subspace_{\EigSpe}$ of dimension $\rdim = 39$. {Again, the evaluation of the eigenvalue and associated eigenspace for the \ROM is about $1555$ times faster than using the \FOM.}

\Cref{fig2:a} depicts the approximation error for the eigenspace, which is correctly below the target accuracy and the a posteriori error estimate. 
\Cref{fig2:b} demonstrates the behavior of the eigenspace surrogate error \eqref{eqn:new:err:est} together with the error calculated for the selected interpolation parameter, the eigenvalue surrogate error \eqref{eqn:eig:surr} and the residual norm scaled by the spectral gap approximation. 
We observe that during the first iterations the influence of the eigenvalue error is similar to that of the residual error. 
However, over the iterations, its impact diminishes considerably, leaving the surrogate error primarily affected by the residual error. This behavior can be explained with \Cref{lemma:EigVal EigVec bound}~\ref{sub-lemma:KT}, i.e. whenever the approximated eigenvalue is close enough to the exact eigenvalue, then their difference decays at least quadratically fast, up to the scaling of the spectral gap, with respect to the residual.

\subsection{Quantum spin systems examples} \label{subsec:QSS}
Let us introduce an auxiliary definition to describe the \QSS examples. 
For $L,j\in \N$ with $j \leq L$ and for any matrix $\matS \in \C^{m \times m}$, we set  
\begin{align*}
    \matS^{(L,j)} \vcentcolon= \left( \bigotimes_{k=1}^{j-1} \IdMat_m \right) \bigotimes \matS \bigotimes \left( \bigotimes_{k=1}^{L-j} \IdMat_m \right) \in \C^{m^L \times m^L},
\end{align*}
where $\bigotimes$ denotes the Kronecker product of matrices, and we use the convention $\bigotimes_{k=1}^{0} \IdMat_m \vcentcolon= 1$.

\subsubsection{The xxz-chain model}\label{subsec:QSS:xxz}
The first \QSS we address is the xxz chain model with open boundary condition \cite{Herbst2022, SRFB2004}. 
For $\prmtr = [\, \prmtrs_1,\prmtrs_2 \,]^\top \in \prmtrSet \coloneqq [-1, \, 2.5]\times[0, \, 3.5]$, the system Hamiltonian has the following affine decomposition structure 
\begin{equation}\label{eqn:QSS:xxz:aff}
	\Matrix(\prmtr) = \Matrix_1+\prmtrs_1 \Matrix_2 - \prmtrs_2 \Matrix_3, 
\end{equation}
where 
\begin{align*}
 \Matrix_1 \vcentcolon= \frac{1}{4} \sum_{j=1}^{L-1} \left( \matS_{x}^{(L,j)} \matS_{x}^{(L,j+1)} + \matS_{y}^{(L,j)} \matS_{y}^{(L,j+1)} \right), 
 \quad 
 \Matrix_2 \vcentcolon= \frac{1}{4} \sum_{j=1}^{L-1}  \matS_{z}^{(L,j)} \matS_{z}^{(L,j+1)} , 
 \quad 
 \Matrix_3 \vcentcolon= \frac{1}{2} \sum_{j=1}^{L} \matS_{z}^{(L,j)} , 
 \quad 
\end{align*}
and the matrices $\matS_{x}, \matS_{y}, \matS_{z}$ are the \textit{spin-1/2-matrices} given by 
\begin{align*}
    \matS_{x} \vcentcolon= 
    \left[
    \begin{array}{cc}
        0 & 1 \\
        1 & 0
    \end{array}
    \right],    \qquad 
    \matS_{y} \vcentcolon= 
    \left[
    \begin{array}{cc}
        0 & -\ImgUnit \\
        \ImgUnit & 0
    \end{array}
    \right],    \qquad 
    \matS_{z} \vcentcolon= 
    \left[
    \begin{array}{cc}
        1 & 0 \\
        0 & -1
    \end{array}
    \right]. 
\end{align*}
For the numerical tests, we set $L=14$, which yields $\fdim=16384$. 

We run \Cref{alg:gap} with $\varepsilon_{\gamma}=10^{-8}$ over a grid $\Xi \subseteq \prmtrSet$ consisting of $35\times35$ Chebyshev-spaced parameter points{. Specifically, $33$ points in each direction are spaced as Chebyshev points of the first kind to which we added the extremal nodes such that the grid includes points from the boundary of $\prmtrSet$}; results are displayed in \Cref{fig3}. 
\begin{figure}[t]
	\centering{
	  \begin{subfigure}[t]{.48\linewidth}
             \input{Img_Tex/Example_2_Plot_1.tex}
    		\subcaption{Plot of the approximated spectral gap $ \gamma^{\subspace_\gamma}(\prmtr)$ over a discrete grid $\Xi$ consisting of $35\times 35$ Chebyshev nodes. }
   		\label{fig3:a}
   	\end{subfigure}
	\hfill
	\begin{subfigure}[t]{.48\linewidth}
        \input{Img_Tex/Example_2_Plot_2.tex}
   		\subcaption{Spectral gap approximation error $\calE^{\subspace_{\gamma}}(\prmtr)$ \eqref{eqn:gap:err} over $\TrainSet$; $\log_{10}$ scale is used. Red crosses are the selected interpolation points by \Cref{alg:gap}.}
   		\label{fig3:b}
   	\end{subfigure}
	}
	
	\caption{\QSS example, xxz chain model: 
    Approximation of the spectral gap over a discrete domain consisting of $35\times 35$ Chebyschev-nodes for $\Matrix(\prmtr)\in\R^{\fdim\times \fdim}$ sparse matrix as in \eqref{eqn:QSS:xxz:aff} with $\fdim=2^{14}=16384$.
    The resulting reduced space $\subspace_{\gamma}$ has dimension $\rdim={265}$.}
	\label{fig3}
\end{figure}
The generated subspace $\subspace_{\gamma}$ has dimension $\rdim={265}$.  
{The average computational time to evaluate the spectral gap for the \FOM is 0.17 seconds while for the \ROM is $1.1\cdot10^{-2}$, therefore achieving a speed-up of a factor of $15$. Note that, for this example, we approximate the smallest eigenvalues and associated eigenvectors of the \ROM using \EIGS.}
The plot on the left, i.e. \Cref{fig3:a}, demonstrates the approximated spectral gap $\gamma^{\subspace_\gamma}(\prmtr)$ in the parameter domain. 
We notice that $\hat \prmtr=[\,-1, \, 0\,]^\top$ is a parameter point of high degeneracy level, where the algebraic multiplicity of the smallest eigenvalue is $m_1(\hat\prmtr)=15$, and this point is successfully captured by our algorithm. 
We can see from \Cref{fig3:b} that the error \eqref{eqn:gap:err} is always smaller than the target accuracy. 
We also observe that the interpolation points are mainly clustered in the lower diagonal region of the rectangular domain.
\begin{figure}[t]
	\centering{
	  \begin{subfigure}[t]{.48\linewidth}
          {\input{Img_Tex/Example_2_Plot_3.tex}}
    		\subcaption{Approximation error $\| \Proj^{\EigSpe^{\bot}_1(\prmtr)} \Vek_1^\basisV(\prmtr) \|$ over a discrete grid $\Xi${; $\log_{10}$ scale is used}. Red crosses are the selected interpolation points by \Cref{alg:sf}.}
   		\label{fig5:a}
   	\end{subfigure}
\hfill
	\begin{subfigure}[t]{.48\linewidth}
		{
%
\tikzexternaldisable
\begin{tikzpicture}

\begin{axis}[%
	width=\imageWidth,
	height=\imageHeight,
	at={(1.011in,0.642in)},
	scale only axis,
	scaled ticks=false,
	grid=both,
	grid style={line width=.1pt, draw=gray!10},
	major grid style={line width=.2pt,draw=gray!50},
	axis lines*=left,
	axis line style={line width=\lineWidth},
xmin=0,
xmax=119,
xlabel style={font=\color{white!15!black}},
xlabel={$j$},
ymode=log,
ymin=1e-15,
ymax=5e14,
yminorticks=true,
ylabel style={font=\color{white!15!black}},
	axis background/.style={fill=white},
	legend style={%
		legend cell align=left, 
		align=left, 
		font=\tiny,
		draw=white!15!black,
		at={(1.0,1.0)},
		anchor=north east,},
]

\addplot [color=mycolor1, line width=1.5*\lineWidth]
 table [x=x1, y=y1, col sep=comma]{Img_Tex/DataCSV/Example_2_4.csv};
 \addlegendentry{$\Delta_{\subspace_\gamma,\varepsilon_\gamma}(\prmtr_{j},\subspace_{j-1})$ (see \eqref{eqn:new:err:est})}
 \addplot [color=mycolor2, dashed, line width=1.5*\lineWidth]
 table [x=x1, y=y2, col sep=comma]{Img_Tex/DataCSV/Example_2_4.csv};
 \addlegendentry{$\|\Res^{\subspace_{j-1}}(\prmtr_{j})\|/\gamma^{\subspace_\gamma}(\prmtr_{j})$}
 \addplot [color=mycolor3, dashdotted, line width=1.5*\lineWidth]
 table [x=x1, y=y3, col sep=comma]{Img_Tex/DataCSV/Example_2_4.csv};

\addlegendentry{$H(\prmtr_{j},\subspace_{j-1})/\gamma^{\subspace_\gamma}(\prmtr_{j})$ (see \eqref{eqn:eig:surr})}

\addplot [color=mycolor4, dashdotted, line width=1.5*\lineWidth]
 table [x=x1, y=y4, col sep=comma]{Img_Tex/DataCSV/Example_2_4.csv};

 \addlegendentry{$\| \Proj^{\EigSpe_1(\prmtr_j)} -\Proj^{\EigSpe^{\subspace_{j-1}}_1(\prmtr_j)}\|$}
\end{axis}
\end{tikzpicture}%
}
   		\subcaption{Decay of the estimator \eqref{eqn:new:err:est} with its residual and eigenvalue error components {along the iterations}.}
   		\label{fig5:b}
   	\end{subfigure}
	}
	
	\caption{\QSS example, xxz chain model: 
    Approximation of the eigenspace over a discrete domain consisting of $35\times 35$ Chebyschev-nodes for $\Matrix(\prmtr)\in\R^{\fdim\times \fdim}$ sparse matrix as in \eqref{eqn:QSS:xxz:aff} with $\fdim=2^{14}= 16384$. 
    The resulting reduced space $\subspace_{\EigSpe}$ has dimension  $\rdim={134}$.   }
	\label{fig5}
\end{figure}

Then, we use the same discrete grid $\Xi$ to run \Cref{alg:sf} for the eigenspace approximation with $\varepsilon_{\EigSpe}=10^{-8}$. We obtain a subspace of dimension $\rdim={134}$; the error is always below the prescribed accuracy, as we can observe from \Cref{fig5:a}, which also displays the selected interpolation points. {Moreover, the evaluation of the eigenvalue and associated eigenspace for the \ROM is around $36$ times faster than using the \FOM.}
\Cref{fig5:b} depicts the behavior of the computed error, the error estimator, and its components along the iterations of \Cref{alg:sf}.
Similarly to the last example, the error estimator is asymptotically dominated by the contribution due to its residual component. 
We observe that {around the $90$-th iteration} the error estimate is smaller than the computed error. 
This is not in contradiction with our theory: 
The selected point at that iteration is $\prmtr=[\,-1,\,0\,]^{\T}$, i.e. the point where the associated eigenspace has dimension larger than one, and since the approximation space used before interpolating does not capture the correct dimension, it follows from \Cref{theo:err:mes} that our error estimate is an upper bound of $\|\Proj^{\EigSpe^{\bot}_1(\mu_j)}\basisW_1^{\basisV_{j-1}}(\prmtr_{j}) \|$ but not an upper bound of the error $\|\Proj^{\EigSpe_1(\mu_j)} -\Proj^{\EigSpe^{\subspace_{j-1}}_1(\mu_j)}\|$. 
Let us stress that, after the convergence of \Cref{alg:sf}, we are always sure that our approximation captures the correct eigenspace dimension due to the strategy of enforcing the conditions provided by \Cref{thm:exact:dim:rec} and thus we can guarantee $\|\Proj^{\EigSpe^{\bot}_1(\mu)}\basisW_1^{\basisV}(\prmtr) \|=\|\Proj^{\EigSpe_1(\mu_j)} -\Proj^{\EigSpe^{\subspace_{\EigSpe}}_1(\mu)}\|$ for all $\prmtr\in\TrainSet$.

In summary, for the xxz-chain model, our method yields small reduced spaces for the approximation of the spectral gap and of the eigenspace associated with the smallest eigenvalue. These spaces allow for high-accuracy approximation with guaranteed error control.

\subsubsection{The bilinear-biquadratic spin-1 chain with uniaxial
single-ion anisotropy}\label{subsec:QSS:blbq}

{In the second \QSS example, we move to the realm of spin-1 chains and test our method on a particularly challenging example known to approach a gapless phase in the thermodynamic limit~\cite{Brehmer2023, DecLS11}, namely the bilinear-biquadratic model with a uniaxial single-ion anisotropy on a chain with open boundaries. 
As a consequence of this physical complexity, the Kolmogorov $n$-width of the ground state spaces decays slowly, presenting a severe challenge for model reduction. 
This example serves to validate our method under non-ideal conditions.}
For $\prmtr = [\prmtrs_1,\prmtrs_2]^\top \in \prmtrSet \coloneqq [-\pi, \, \pi]\times[-2, \, 3]$, the system Hamiltonian has the following affine decomposition structure  
\begin{equation}\label{eqn:QSS:bblq:aff}
	\Matrix(\prmtr) = \cos(\prmtrs_1)\Matrix_1+\sin(\prmtrs_1) \Matrix_2 + \prmtrs_2 \Matrix_3, 
\end{equation}
where 
\begin{align*}
    \Matrix_1 \vcentcolon= \hspace*{-1mm} \sum_{a\in \{x,y,z\} } \sum_{j=1}^{L-1} \matS_a^{(L,j)}\matS_a^{(L,j+1)},
    \quad 
    \Matrix_2 \vcentcolon= \hspace*{-1mm} \sum_{a,b\in \{x,y,z\} } \sum_{j=1}^{L-1} \left( \matS_a^{(L,j)} \matS_b^{(L,j+1)}  \right)^2 ,
    \quad 
    \Matrix_3 \vcentcolon= \sum_{j=1}^{L} \left( \matS_z^{(L,j)} \right)^2,
\end{align*}
and the matrices $\matS_{x},\matS_{y},\matS_{z}$ are now the \textit{spin-1-matrices} given by
\begin{align*}
    \matS_{x} \;\vcentcolon=\; 
    \frac{1}{\sqrt{2}}
    \left[
    \begin{array}{ccc}
        0 & 1 & 0 \\
        1 & 0 & 1 \\
        0 & 1 & 0
    \end{array}
    \right],    \qquad 
    \matS_{y}\; \vcentcolon=\; 
    \frac{-\ImgUnit}{\sqrt{2}}
    \left[
    \begin{array}{ccc}
        0  & 1  & 0 \\
        -1 & 0  & 1 \\
        0  & -1 & 0
    \end{array}
    \right],    \qquad 
    \matS_{z} \;\vcentcolon=\; 
    \left[
    \begin{array}{ccc}
        1 & 0 & 0 \\
        0 & 0 & 0 \\
        0 & 0 & -1
    \end{array}
    \right]. 
\end{align*} 
The comprehensive physics inherent in this model, characterized by several gapped and critical phases, has been extensively explored through both theoretical analysis and numerical simulations, such as exhaustive DMRG calculations. 
A notable reference is the ground-state phase diagram detailed in \cite{DecLS11}, which also summarizes previous investigations. For our numerical tests, we set $L=10$, which results in $\fdim=59049$. 

\begin{figure}[t]	
	\centering{
	\begin{subfigure}[t]{0.48\textwidth}
	    \input{Img_Tex/Example_3_Plot_1_tris.tex}
     \subcaption{{Heat map} of the approximated spectral gap $ \gamma^{\subspace_\gamma}(\prmtr)$ over a discrete grid $\Xi$ consisting of $50\times 50$ Chebyshev nodes.}
     \label{fig6:a}
	\end{subfigure}
\hfill
 \begin{subfigure}[t]{0.48\textwidth}
	    \input{Img_Tex/Example_3_Plot_1.tex}
     \subcaption{Spectral gap approximation error $\calE^{\subspace_{\gamma}}(\prmtr)$ \eqref{eqn:gap:err} over~$\Xi$.}
     \label{fig6:b}
	\end{subfigure}
  \caption{\QSS example, bblq chain model:
  Approximation of the spectral gap for 
  $\Matrix(\prmtr)\in\R^{\fdim\times \fdim}$ sparse matrix as in \eqref{eqn:QSS:bblq:aff} with $\fdim=59049${; $\log_{10}$ scale is used}.  
  The resulting reduced space $\subspace_{\gamma}$ has dimension $\rdim={859}$. {Red crosses in both the plots indicate the selected interpolation points by \Cref{alg:gap}.}
	}
	\label{fig6}
    }
\end{figure}

We run \Cref{alg:gap} with $\varepsilon_{\gamma}=10^{-2}$ over a grid $\TrainSet \subseteq \prmtrSet$ consisting of $50\times50$ Chebyshev-spaced parameter points {with the same criteria applied for the xxz-chain model}; results are displayed in \Cref{fig6} and, in particular, \Cref{fig6:a} depicts the {heat map of the} approximated spectral gap in the domain $\Xi${, with the red crosses being the selected interpolation points; we observe that most of the points picked by the greedy algorithm are where the gap is rather small.}
The algorithm returns a subspace $\subspace_{\gamma}$ of dimension $\rdim={859}$, which rigorously approximates the spectral gap on $\TrainSet$ under the prescribed accuracy; see \Cref{fig6:b}. 
Despite the relatively large reduced dimension, the computation of the spectral gap for the reduced matrix of dimension {$859$ was $11.6$ times faster, on average,} than the calculation of the spectral gap for the full size matrix. 
{Specifically, we measured over the entire domain $\Xi$ an average computational time of $0.36$ seconds to compute the spectral gap of the the reduced matrix; in contrast, we found an average computational time of $4.2$ seconds to compute the one of the sparse full size matrix.
This measured speed-up, while modest compared to the other examples, is a result of the certification process: The rigorous error tolerance $\varepsilon_\gamma = 10^{-2}$ for the challenging blbq-chain model necessitates a reduced model of size $\rdim=859$, which accurately reflects the problem's complexity and slow Kolmogorov $n$-width decay. 
The primary goal of this step is the reliable certification of the spectral gap, which then enables an efficient construction of the final eigenspace approximation.}

Using $\subspace_{\gamma}$, we run \Cref{alg:sf} for the eigenspace approximation, where we consider the same domain $\TrainSet$ and target error accuracy $\varepsilon_{\EigSpe}=10^{-4}$. 
\Cref{fig7:a} illustrates the parameters selected by \Cref{alg:sf} to construct $\subspace_\EigSpe$ and the resulting errors on $\TrainSet$. 
We see that in certain regions of $\prmtrSet$ the method selects many parameters close to each other, indicating that the eigenvectors that span the eigenspaces on these parameters are (almost) linearly independent despite the parameters being close to each other. 
This observation suggests a slow decay of $d_n(\calM)$, i.e. the Kolmorogov $n$-width associated to the eigenvalue problem as discussed in \Cref{subsec:problem setting}.
Despite the complexity of the model, we see that within the domain $\Xi$, the calculated error remains below the desired accuracy. 
{Moreover, the average computational time to compute the complete eigendecomposition of the the reduced matrix is only $0.11$ seconds, while, the one to compute the smallest eigenvalues and associated eigenvectors of the sparse full size matrix was around $3.54$ seconds, making the \ROM more than $30$ times faster than solving the original challenging bllq chain model. Note that, due to the relatively large size of the \ROM compared to the one of the previous examples, for the bblq chain model we found beneficial to use the iterative methods embedded in \EIGS also for the evaluation of the eigenvalues and associated eigenvectors of the \ROM.
}
\Cref{fig7:b} portrays the behavior of the error estimate and its components along the iterations $j$; the slow decay of the error estimator \eqref{eqn:new:err:est} confirms the slow decay of $d_n(\calM)$. As noted in \Cref{subsec:QSS:xxz}, initially, the error estimate is primarily influenced by the eigenvalue error component. However, as more interpolation points are chosen, the accuracy of the eigenvalue approximation improves significantly, resulting in the residual error component becoming predominant.
\begin{figure}[t]
	\centering{
	  \begin{subfigure}[t]{.48\linewidth}
			\input{Img_Tex/Example_3_Plot_3.tex}
    		\subcaption{Approximation error $\| \Proj^{\EigSpe^{\bot}_1(\prmtr)} \Vek_1^\basisV(\prmtr) \|$ over $\Xi$. Red crosses are the selected interpolation points by \Cref{alg:sf}{; $\log_{10}$ scale is used}.}  
   		\label{fig7:a}
   	\end{subfigure}
\hfill
    \begin{subfigure}[t]{.48\linewidth}
\tikzexternaldisable
\begin{tikzpicture}

\begin{axis}[%
	width=\imageWidth,
	height=\imageHeight,
	scale only axis,
	scaled ticks=false,
	grid=both,
	grid style={line width=.1pt, draw=gray!10},
	major grid style={line width=.2pt,draw=gray!50},
	axis lines*=left,
	axis line style={line width=\lineWidth},
xmin=0,
xmax=351,
xlabel style={font=\color{white!15!black}},
xlabel={$j$},
ymode=log,
ymin=1e-12,
ymax=5e9,
yminorticks=true,
ylabel style={font=\color{white!15!black}},
	axis background/.style={fill=white},
	legend style={%
		legend cell align=left, 
		align=left, 
		font=\tiny,
		draw=white!15!black,
		at={(1.0,1.0)},
		anchor=north east,},
]

\addplot [color=mycolor1, line width=1.5*\lineWidth]
 table [x=x1, y=y1, col sep=comma]{Img_Tex/DataCSV/Example_3_2.csv};
 \addlegendentry{$\Delta_{\subspace_\gamma,\varepsilon_\gamma}(\prmtr_{j},\ISubspace{j-1})$ (see \eqref{eqn:new:err:est})}
 \addplot [color=mycolor2, dashed, line width=1.5*\lineWidth]
 table [x=x1, y=y2, col sep=comma]{Img_Tex/DataCSV/Example_3_2.csv};
 
 \addlegendentry{$\|\ResSubsp{\ISubspace{j-1}}(\prmtr_{j})\|/\gamma^{\subspace_\gamma}(\prmtr_{j})$}
 \addplot [color=mycolor3, dashdotted, line width=1.5*\lineWidth]
 table [x=x1, y=y3, col sep=comma]{Img_Tex/DataCSV/Example_3_2.csv};

\addlegendentry{$H(\prmtr_{j},\ISubspace{j-1})/\gamma^{\subspace_\gamma}(\prmtr_{j})$ (see \eqref{eqn:eig:surr})}
\end{axis}
\end{tikzpicture}%
   		\subcaption{Decay of the estimator \eqref{eqn:new:err:est} with its residual and eigenvalue error components {along the iterations}.} 
   		\label{fig7:b}
   	\end{subfigure}

	}	
\caption{\QSS example, bblq chain model: Approximation of the eigenspace for $\Matrix(\prmtr)\in\R^{\fdim\times \fdim}$ sparse matrix as in \eqref{eqn:QSS:bblq:aff} with $\fdim=3^{10}=59049$.
The resulting reduced space $\subspace_{\EigSpe}$ has dimension $\rdim={409}$. }
\label{fig7}
\end{figure}

\Cref{fig8} shows why this problem leads to a large subspace dimension. 
Indeed, for the spectral gap approximation, $\ell(j)$ is at least $2$ since the eigenvector associated with the smallest and second smallest eigenvalue needs to be added to the subspace to guarantee the interpolation properties. 
For this example, $\ell(j)$ is greater than $2$ in almost all iterations, indicating that there is degeneracy in the smallest or second smallest eigenvalue (or in both); see \Cref{fig8:a} (blue line). 
Instead, for the eigenspace approximation, $\ell(j)$ is at least one, and in \Cref{fig8:a} (red dashed line), we observe $\ell(j)=2$ several times. 
\begin{figure}[t]
	\centering{
	  \begin{subfigure}[t]{.48\linewidth}
%
\tikzexternaldisable
\begin{tikzpicture}

\begin{axis}[%
	width=\imageWidth,
	height=0.95*\imageHeight,
	scale only axis,
	scaled ticks=false,
	grid=both,
	grid style={line width=.1pt, draw=gray!10},
	major grid style={line width=.2pt,draw=gray!50},
	axis lines*=left,
	axis line style={line width=\lineWidth},
xmin=2,
xmax=341,
xlabel style={font=\color{white!15!black}},
xlabel={$j$-th iteration},
ymin=0.5,
ymax=5,
yminorticks=true,
ylabel style={font=\color{white!15!black}},
	axis background/.style={fill=white},
	legend style={%
		legend cell align=left, 
		align=left, 
		font=\tiny,
		draw=white!15!black,
		at={(0.9,0.8)},
		anchor=south east,},
]

\addplot [color = mycolor1, line width=\lineWidth] table [x=x1, y=ne, col sep=comma]{Img_Tex/DataCSV/Example_3_1.csv};
 \addlegendentry{$\ell(j)$ for Alg.~\ref{alg:gap}}

\addplot [color = mycolor2, dashed, line width=\lineWidth] table [x=x1, y=ne, col sep=comma]{Img_Tex/DataCSV/Example_3_1_bis.csv};
\addlegendentry{$\ell(j)$ for Alg.~\ref{alg:sf}}
\end{axis}
\end{tikzpicture}%
    		\subcaption{$\ell(j)$ as a function of the iteration $j$ for \Cref{alg:gap} ($\ell(j)=m_1(\prmtr_j)+m_2(\prmtr_j)$) and \Cref{alg:sf} ($\ell(j)=m_1(\prmtr_j)$) for $\prmtr_j$ the selected parameters.}
   		\label{fig8:a}
   	\end{subfigure}
\hfill
	\begin{subfigure}[t]{.48\linewidth}
%
\tikzexternaldisable
\begin{tikzpicture}

\begin{axis}[%
	width=\imageWidth,
	height=0.95*\imageHeight,
	scale only axis,
	scaled ticks=false,
	grid=both,
	grid style={line width=.1pt, draw=gray!10},
	major grid style={line width=.2pt,draw=gray!50},
	axis lines*=left,
	axis line style={line width=\lineWidth},
xmin=1,
xmax=2500,
xlabel style={font=\color{white!15!black}},
xlabel={$i$},
ymode=log,
ymin=1e-6,
ymax=6e1,
yminorticks=true,
ylabel style={font=\color{white!15!black}},
	axis background/.style={fill=white},
	legend style={%
		legend cell align=left, 
		align=left, 
		font=\tiny,
		draw=white!15!black,
		at={(0.65,0.8)},
		anchor=south east,},
]

\addplot [color = mycolor1, line width=\lineWidth] table [x=x1, y=y1, col sep=comma]{Img_Tex/DataCSV/Example_3_GAP_ETA.csv};
 \addlegendentry{$\calF_{\gamma}(\prmtr_i)$ (see \eqref{eqn:plot:a})}
\addplot [color = mycolor2, dashed, line width=\lineWidth] table [x=x1, y=y1, col sep=comma]{Img_Tex/DataCSV/Example_3_EIG_ETA.csv};
 \addlegendentry{$\calF_{\EigSpe}(\prmtr_i)$ (see \eqref{eqn:plot:b})}
\end{axis}
\end{tikzpicture}%
   		\subcaption{Illustration of the terms \eqref{eqn:plot:a}, \eqref{eqn:plot:b} for the exact recovery of eigenspace dimensions for all the parameter points $\prmtr_i\in\TrainSet$ with $|\TrainSet|=2500$.}
   		\label{fig8:b}
   	\end{subfigure}
	}
	
	\caption{\QSS example, bblq chain model: Discrete grid $\Xi $ of $50\times 50$ Chebyschev nodes for sparse system matrix $\Matrix(\prmtr)\in\R^{\fdim\times \fdim}$ as in \eqref{eqn:QSS:bblq:aff} with $\fdim=3^{10}=59049$. 
    }
	\label{fig8}
\end{figure}

We recall that if \eqref{eqn:hyp} for $k=1$ holds, then our spectral gap approximation is certified, that is, $\calE^{\subspace_{\gamma}}\le\varepsilon_{\gamma}$ for all $\prmtr\in\Xi$, and if \eqref{eqn:hyp} for $k=1$ holds, then the dimension of $\EigSpe_1(\prmtr)$ is recovered exactly by $\EigSpe^{\subspace_{\EigSpe}}_1(\prmtr)$.
Since for the bblq chain model the ground state space has dimension larger than one for a number of parameters, we also plot
\begin{align}
     \calF_{\gamma}(\prmtr)\;
     \vcentcolon =&\;
     \eta_{\ast}^{(J_{\gamma})}(\prmtr, s(1))-\lambda^{\subspace_{\gamma}}_{s(1)}(\prmtr) -  \varepsilon^{(J_{\gamma})}(\prmtr,s(1)),\quad\quad \quad s(1)=m_1(\prmtr,\subspace_\gamma),
     \label{eqn:plot:a}\\
     \calF_{\EigSpe}(\prmtr)\;\vcentcolon=&\;\eta_{\ast}^{(J_{\EigSpe})}(\prmtr,s(1))-\lambda^{\subspace_{\EigSpe}}_{s(1)}(\prmtr)-\varepsilon^{(J_{\EigSpe})}(\prmtr,s(1)),\quad\; \quad s(1)=m_1(\prmtr,\subspace_\EigSpe),
     \label{eqn:plot:b}
\end{align}
after termination of \Cref{alg:gap} and of \Cref{alg:sf}, respectively,
in \Cref{fig8:b}.
We observe that both are always greater than zero for all $\prmtr\in\Xi$. 
This means that \eqref{eqn:hyp} for $k=1$ holds for both subspaces constructed. {We conclude by mentioning that for both $\subspace_\gamma$ and $\subspace_\EigSpe$ there were around $300$ and $30$ points, respectively, over $2500$ not satisfying the condition \eqref{eqn:hyp} for $k=1$; thus, the exact recovery of $m_1(\prmtr)$ was enforced through lines \ref{alg gap: second loop}-\ref{alg gap: end second loop} of \Cref{alg:gap} for the spectral gap and lines \ref{alg: second loop}-\ref{alg: end second loop} of \Cref{alg:sf} for the eigenspace approximation.}

In summary, regarding the challenging blbq-chain model, our method is again able to deliver reduced spaces for the approximation of the spectral gap and of the ground states with certified error control, despite resulting in reduced spaces of larger size with respect to the other test problems. 
However, we remark that this is caused by the complexity of the model, with many phase transitions causing several eigenvalues to have an algebraic multiplicity larger than one; this naturally promotes large reduced spaces.  
It is also worth mentioning that the refinement of the set $\Xi$ leads to a further increase in the dimensions of the reduced spaces and, in particular, to an even slower decrease in the error indicator function --- another circumstantial evidence for the complexity of the model. 
Finally, the blbq-chain model also shows that our two-stage procedure successfully recovers the correct dimension of the ground state, which is often greater than one for this example.
{The computational results demonstrate that our method provides rigorous certification even for problems where the Kolmogorov $n$-width decays slowly. 
The achieved speed-up reflects the true complexity of approximating the solution set, a challenge that any model reduction technique would face.}

\section{Conclusions}\label{conclusion}
In this work, we have considered the approximation problem of the smallest eigenpair of a large-scale parametric Hermitian matrix $\Matrix(\prmtr)$ using the subspace approach, which can be interpreted as projection-based \MOR or \RBM. 
For the construction of the subspaces, we have relied on greedy-type algorithms that ensure to computation efficiency with error certification for the approximation. 

After systematically discussing the connection between \RBM for source problems and eigenvalue problems in \Cref{subsec:problem setting} and reviewing some general error bounds for eigenvalue problems in \Cref{subsec:EigVal EigVec bounds}, we have introduced a novel error bound for the subspace approximation of the eigenspace $\EigSpe_1(\prmtr)$ associated with the smallest eigenvalue $\lambda_1(\prmtr)$ (see \Cref{thm:error bound EigVec}), from which an a posteriori error estimate has been derived; see \Cref{prop1}. 

As the spectral gap $\gamma(\prmtr)=\tilde \lambda_2(\prmtr) - \lambda_1(\prmtr) $ in \eqref{eq:gap full space} is crucial for the approximation of $\EigSpe_1(\prmtr)$, we have proposed a practically computable eigenvalue lower bound $\lambda_k^\SLB(\prmtr,\ISubspace{J})$ in \Cref{prop:EigVal SLB} and based on that upper and lower bounds for $\gamma(\prmtr)$ in \Cref{subsec:gap bounds}. 

To approximate $\EigSpe_1(\prmtr)$, we have then presented a two-stage greedy strategy in \Cref{sec4}: A first greedy algorithm is executed and results in a subspace for the approximation of $\gamma(\prmtr)$ with certified error control, see \Cref{alg:gap}, based on which a second greedy algorithm is executed and results in a subspace for the approximation of $\EigSpe_1(\prmtr)$ with certified error control, see \Cref{alg:sf}. 
A remarkable feature of the two-stage designed procedure, which we did not find in other methods described in the literature, is that the computed eigenspace approximation, under an efficient verifiable condition during the reduced space construction, has the same dimension of the eigenspace of the original problems for all parameters in the considered domain; see \Cref{thm:exact:dim:rec}. 
In addition, the developed theory could be extended to approximate the eigenspace associated with the $k$-th smallest eigenvalue.

Numerical experiments on randomly generated dense Hermitian matrices and on \QSS have validated our proposed method and associated theoretical findings.

\appendix
\section{An auxiliary property}
\begin{lemma} \label{lem:aux func monoton increasing ell bumps}
For $\rho,\lambda_1, \dots, \lambda_k \in \R$ satisfying $\rho \geq 0$ and $\lambda_1 \leq \cdots \leq \lambda_k$, we consider the function 
    \begin{align*}
    h_k: \R \rightarrow \R, \, \eta \mapsto \eta - \frac{2\rho^2}{ g_k(\eta) + \sqrt{ g_k(\eta)^2 + 4 \rho^2 }} 
\end{align*}
with $g_k(\eta) \coloneqq \min\{ |\eta-\lambda_j| : 1\leq j \leq k \}$.
Then, the function $h_k$ increases monotonically.
\end{lemma}
We note that \eqref{eq:aux func concrete k} is exactly $h_k$ with $\rho = \rho^{(J)} (\prmtr,s)$ and $\lambda_k = \RedEigVal{\ISubspace{J}}_k(\prmtr)$ for $1\leq j \leq k$. 
\begin{proof} 
Without loss of generality, we can assume $\lambda_1 < \cdots < \lambda_k$, since $g_k$ and $ h_k$ remain unchanged after removing duplicates among $\{\lambda_j\}_{j=1}^k$. 
Let us define two sets of auxiliary intervals 
\begin{align*}
    I_j \vcentcolon=
    \begin{cases}
        (-\infty, \lambda_1), & j=1, \\
        \left(\tfrac{1}{2}(\lambda_{j-1}+\lambda_j), \lambda_j \right), & 2\leq j \leq k,
    \end{cases}\quad\text{and}\quad  J_j \vcentcolon=
    \begin{cases}
        \left(\tfrac{1}{2}(\lambda_{j-1}+\lambda_j), \lambda_j \right), & 1\leq j \leq k-1, \\
        (\lambda_k,+\infty), & j = k.
    \end{cases}
\end{align*}
We observe that for $1\leq j \leq k$, the function $g_k$ is reduced to a simple expression, i.e. 
\begin{align*}
    g_k(\eta) = 
    \begin{cases}
        \lambda_j - \eta, & \eta \in I_j, \\
        \eta - \lambda_j, & \eta \in J_j. 
    \end{cases}
\end{align*}
Therefore, we can also simplify the expression of $h_k$ in each interval $I_j,J_j$, and obtain 
\begin{align*}
    h_k'(\eta) = 
    \begin{cases}
        \psi(\eta), & \eta \in I_j, \\
        \varphi(\eta), & \eta \in J_j, \\
    \end{cases}
\end{align*}
where 
\begin{align*}
    \psi_j : I_j \rightarrow \R, \, \eta &\mapsto 1 - \frac{2 \rho^2 \left( 1+ \frac{\lambda_j-\eta}{2\sqrt{(\lambda_j-\eta)^2+4\rho^2}} \right)}{\left(\sqrt{(\lambda_j-\eta)^2+4\rho^2}+\lambda_j-\eta\right)^2},
\end{align*}
and 
\begin{align*}
    \varphi_j : J_j \rightarrow \R, \, \eta &\mapsto 1 + \frac{2 \rho^2 \left( 1+ \frac{\eta-\lambda_j}{2\sqrt{(\eta-\lambda_j)^2+4\rho^2}} \right)}{\left(\sqrt{(\eta-\lambda_j)^2+4\rho^2}+\eta-\lambda_j\right)^2}.
\end{align*}
Clearly $\varphi_j$ is always positive in $J_j$. 
As for $\psi_j$, we observe that $\psi_j$ is monotonically decreasing with $\psi_j(\lim_{\eta \uparrow \lambda_j} \eta ) = 0$, and so $\psi_j$ is also always non-negative in $I_j$. 
Consequently, $h_k'$ is non-negative in all intervals $I_j,J_j$, meaning that $h_k$ is monotonically increasing in $I_j,J_j$. 
Since $h_k$ is a continuous function, we deduce that $h_k$ is indeed monotonically increasing in the whole real line.
\end{proof}

\subsection*{Acknowledgments}
Funded by Deutsche Forschungsgemeinschaft (DFG, German Research Foundation) under Germany's Excellence Strategy - EXC 2075 – 390740016. We acknowledge the support by the Stuttgart Center for Simulation Science (SimTech).
MM acknowledges funding from the BMBF (grant no.~05M22VSA).

\bibliographystyle{plain-doi}
\bibliography{journalabbr,Literature_MSZ24}    

\begin{thebibliography}{10}

\bibitem{AlgBB25}
\textsc{M.~Alghamdi, D.~Boffi, and F.~Bonizzoni}.
\newblock \href{https://doi.org/10.1016/j.cam.2024.116270}{A greedy mor method for the tracking of eigensolutions to parametric elliptic pdes}.
\newblock {\em J. Comput. Appl. Math.}, 457:116270, 2025.

\bibitem{AndS12}
\textsc{R.~Andreev and C.~Schwab}.
\newblock \href{https://doi.org/10.1007/978-3-642-22061-6_7}{Sparse {tensor} {approximation} of {parametric} {eigenvalue} {problems}}.
\newblock In \textsc{I.~G. Graham, T.~Y. Hou, O.~Lakkis, and R.~Scheichl}, editors, {\em Numerical {Analysis} of {Multiscale} {Problems}}, volume~83, pages 203--241. Springer Berlin Heidelberg, Berlin, Heidelberg, 2012.
\newblock Series Title: Lecture Notes in Computational Science and Engineering.

\bibitem{Atk89}
\textsc{K.~E. Atkinson}.
\newblock \href{http://www.worldcat.org/isbn/0471500232}{{\em {An Introduction to Numerical Analysis}}}.
\newblock John Wiley \& Sons, New York, second edition, 1989.

\bibitem{BaiDDRV00}
\textsc{Z.~Bai, J.~Demmel, J.~J. Dongarra, A.~Ruhe, and H.~van~der Vorst}.
\newblock \href{https://api.semanticscholar.org/CorpusID:117183209}{Templates for the solution of algebraic eigenvalue problems}.
\newblock In {\em Software, environments, tools}, 2000.

\bibitem{Bau84}
\textsc{H.~Baumgärtel}.
\newblock \href{https://doi.org/doi:10.1515/9783112721810}{{\em Analytic Perturbation Theory for Matrices and Operators}}.
\newblock De Gruyter, Berlin, Boston, 1984.

\bibitem{BenGW15}
\textsc{P.~Benner, S.~Gugercin, and K.~Willcox}.
\newblock \href{https://doi.org/10.1137/130932715}{A survey of projection-based model reduction methods for parametric dynamical systems}.
\newblock {\em {SIAM} Rev.}, 57:483--531, June 2015.

\bibitem{BinCDDPW11}
\textsc{P.~Binev, A.~Cohen, W.~Dahmen, R.~DeVore, G.~Petrova, and P.~Wojtaszczyk}.
\newblock \href{https://doi.org/10.1137/100795772}{Convergence rates for greedy algorithms in reduced basis methods}.
\newblock {\em SIAM J. Math. Anal.}, 43(3):1457--1472, January 2011.

\bibitem{Boffi_2010}
\textsc{D.~Boffi}.
\newblock \href{https://doi.org/10.1017/S0962492910000012}{Finite element approximation of eigenvalue problems}.
\newblock {\em Acta Numer.}, 19:1–120, 2010.

\bibitem{BHP24}
\textsc{D.~Boffi, A.~Halim, and G.~Priyadarshi}.
\newblock \href{https://doi.org/10.1007/s40314-024-02917-x}{Reduced basis approximation of parametric eigenvalue problems in presence of clusters and intersections}, 2024.

\bibitem{Brehmer2023}
\textsc{P.~Brehmer, M.~F. Herbst, S.~Wessel, M.~Rizzi, and B.~Stamm}.
\newblock \href{https://doi.org/10.1103/PhysRevE.108.025306}{Reduced basis surrogates for quantum spin systems based on tensor networks}.
\newblock {\em Phys. Rev. E}, 108:025306, Aug 2023.

\bibitem{BuhEOR14}
\textsc{A.~Buhr, C.~Engwer, M.~Ohlberger, and S.~Rave}.
\newblock \href{https://arxiv.org/abs/1407.8005}{A numerically stable a posteriori error estimator for reduced basis approximations of elliptic equations}, 2014.

\bibitem{DavK70}
\textsc{C.~Davis and W.~M. Kahan}.
\newblock \href{http://www.jstor.org/stable/2949580}{The rotation of eigenvectors by a perturbation. iii}.
\newblock {\em {SIAM} J. Numer. Anal.}, 7(1):1--46, 1970.

\bibitem{DecLS11}
\textsc{G.~De~Chiara, M.~Lewenstein, and A.~Sanpera}.
\newblock \href{https://doi.org/10.1103/PhysRevB.84.054451}{Bilinear-biquadratic spin-1 chain undergoing quadratic zeeman effect}.
\newblock {\em Phys. Rev. B}, 84:054451, Aug 2011.

\bibitem{VliM12}
\textsc{J.~De~Vlieger and K.~Meerbergen}.
\newblock \href{https://lirias.kuleuven.be/1947222&lang=en}{A subspace method for unimodal symmetric eigenvalue optimization problems involving large scale matrices}.
\newblock Technical report, KU Leuven, Leuven, Belgium, 2012.

\bibitem{DoeE24}
\textsc{J.~Dölz and D.~Ebert}.
\newblock \href{https://doi.org/10.1137/22m1529324}{On uncertainty quantification of eigenvalues and eigenspaces with higher multiplicity}.
\newblock {\em {SIAM} J. Numer. Anal.}, 62(1):422–451, February 2024.

\bibitem{GarS24}
\textsc{L.~Garrigue and B.~Stamm}.
\newblock \href{https://arxiv.org/abs/2408.11924}{On reduced basis methods for eigenvalue problems, and on its coupling with perturbation theory}, 2024.
\newblock arXiv:2408.11924.

\bibitem{GreU19}
\textsc{C.~Greif and K.~Urban}.
\newblock \href{https://doi.org/10.1016/j.aml.2019.05.013}{Decay of the {Kolmogorov} {N} -width for wave problems}.
\newblock {\em Appl. Math. Lett.}, 96:216--222, October 2019.

\bibitem{GSH23}
\textsc{L.~Grubišić, M.~Saarikangas, and H.~Hakula}.
\newblock \href{https://doi.org/10.1007/s00211-022-01339-3}{Stochastic collocation method for computing eigenspaces of parameter-dependent operators}.
\newblock {\em Numer. Math.}, 153(1):85--110, January 2023.

\bibitem{GugM23}
\textsc{N.~Guglielmi and M.~Manucci}.
\newblock \href{https://doi.org/10.1137/22M1520189}{Model order reduction in contour integral methods for parametric pdes}.
\newblock {\em {SIAM} J. Sci. Comput.}, 45(4):A1711--A1740, 2023.

\bibitem{HeiS22}
\textsc{C.~Heissenberg and A.~Sagnotti}.
\newblock \href{https://doi.org/10.1017/9781108952002}{{\em Classical and Quantum Statistical Physics: Fundamentals and Advanced Topics}}.
\newblock Cambridge University Press, 2022.

\bibitem{Herbst2022}
\textsc{M.~F. Herbst, B.~Stamm, S.~Wessel, and M.~Rizzi}.
\newblock \href{https://doi.org/10.1103/PhysRevE.105.045303}{Surrogate models for quantum spin systems based on reduced-order modeling}.
\newblock {\em Phys. Rev. E}, 105:045303, Apr 2022.

\bibitem{HesRS16}
\textsc{J.~Hesthaven, G.~Rozza, and B.~Stamm}.
\newblock \href{https://doi.org/10.1007/978-3-319-22470-1}{{\em Certified Reduced Basis Methods for Parametrized Partial Differential Equations}}.
\newblock Springer, January 2016.

\bibitem{HWD17}
\textsc{T.~Horger, B.~Wohlmuth, and T.~Dickopf}.
\newblock \href{https://doi.org/10.1051/m2an/2016025}{Simultaneous reduced basis approximation of parameterized elliptic eigenvalue problems}.
\newblock {\em ESAIM: Math. Model. Numer. Anal.}, 51(2):443--465, March 2017.

\bibitem{HuyRSP07}
\textsc{D.~B.~P. Huynh, G.~Rozza, S.~Sen, and P.~A. T.}
\newblock \href{https://doi.org/10.1016/j.crma.2007.09.019}{A successive constraint linear optimization method for lower bounds of parametric coercivity and inf-sup stability constants}.
\newblock {\em C. R. Math.}, 345(8):473--478, 2007.

\bibitem{KanMMM18}
\textsc{F.~Kangal, K.~Meerbergen, E.~Mengi, and W.~Michiels}.
\newblock \href{https://doi.org/10.1137/16M1070025}{A subspace method for large scale eigenvalue optimization}.
\newblock {\em {SIAM} J. Matrix Anal. Appl.}, 39(1):48--82, 2017.

\bibitem{Kato76}
\textsc{T.~Kato}.
\newblock {\em {Perturbation theory for linear operators; 2nd ed.}}
\newblock Grundlehren der mathematischen Wissenschaften: a series of comprehensive studies in mathematics. Springer, Berlin, 1976.

\bibitem{KreLV18}
\textsc{D.~Kressner, D.~Lu, and B.~Vanderycken}.
\newblock \href{https://doi.org/10.1137/17M1127545}{Subspace acceleration for the crawford number and related eigenvalue optimization problems}.
\newblock {\em {SIAM} J. Matrix Anal. Appl.}, 39(2):961---982, 2018.

\bibitem{KreV14}
\textsc{D.~Kressner and B.~Vandereycken}.
\newblock \href{https://doi.org/10.1137/120869432}{Subspace methods for computing the pseudospectral abscissa and the stability radius}.
\newblock {\em {SIAM} J. Matrix Anal. Appl.}, 35(1):292--313, 2014.

\bibitem{LiL05}
\textsc{C.~K. Li and R.~C. Li}.
\newblock \href{https://doi.org/https://doi.org/10.1016/j.laa.2004.08.026}{A note on eigenvalues of perturbed hermitian matrices}.
\newblock {\em Linear Algebra Appl.}, 395:183--190, 2005.

\bibitem{MacMOPR00}
\textsc{L.~Machiels, Y.~Maday, I.~B. Oliveira, A.~Patera, and D.~V. Rovas}.
\newblock \href{https://doi.org/https://doi.org/10.1016/S0764-4442(00)00270-6}{Output bounds for reduced-basis approximations of symmetric positive definite eigenvalue problems}.
\newblock {\em Comptes Rendus de l'Académie des Sciences - Series I - Mathematics}, 331(2):153--158, 2000.

\bibitem{ManMG24}
\textsc{M.~Manucci, E.~Mengi, and N.~Guglielmi}.
\newblock \href{https://doi.org/10.48550/arXiv.2409.05791}{Uniform approximation of eigenproblems of a large-scale parameter-dependent {H}ermitian matrix}.
\newblock {\em arXiv:2409.05791}, 2024.

\bibitem{ManN15}
\textsc{A.~Manzoni and F.~Negri}.
\newblock \href{https://doi.org/10.1007/s10444-015-9413-4}{Heuristic strategies for the approximation of stability factors in quadratically nonlinear parametrized {PDE}s}.
\newblock {\em Adv. Comput. Math.}, 41(5):1255--1288, 2015.

\bibitem{Meerbergen2017}
\textsc{K.~Meerbergen, E.~Mengi, W.~Michiels, and R.~Van~Beeumen}.
\newblock \href{https://doi.org/10.1093/imanum/drw065}{Computation of pseudospectral abscissa for large-scale nonlinear eigenvalue problems}.
\newblock {\em IMA J. Numer. Anal.}, 37(4):1831--1863, 2017.

\bibitem{MehM24}
\textsc{V.~Mehrmann and E.~Mengi}.
\newblock \href{https://doi.org/10.48550/arXiv.2409.04297}{Minimization of the pseudospectral abscissa of a quadratic matrix polynomial}.
\newblock {\em arXiv:2409.04297}, 2024.

\bibitem{MenYK14}
\textsc{E.~Mengi, E.~A. Yildirim, and M.~Kili\c{c}}.
\newblock \href{https://doi.org/10.1137/130933472}{Numerical optimization of eigenvalues of {H}ermitian matrix functions}.
\newblock {\em {SIAM} J. Matrix Anal. Appl.}, 35(2):699--724, 2014.

\bibitem{OhlR16}
\textsc{M.~Ohlberger and S.~Rave}.
\newblock \href{http://www.iam.fmph.uniba.sk/amuc/ojs/index.php/algoritmy/article/view/389}{Reduced basis methods: Success, limitations and future challenges}.
\newblock {\em Proceedings of the Conference Algoritmy}, pages 1--12

\bibitem{parlett98}
\textsc{B.~N. Parlett}.
\newblock \href{https://epubs.siam.org/doi/book/10.1137/1.9781611971163}{{\em The Symmetric Eigenvalue Problem}}.
\newblock Society for Industrial and Applied Mathematics, Philadelphia, 1998.

\bibitem{PraB24}
\textsc{D.~Pradovera and A.~Borghi}.
\newblock \href{https://doi.org/10.1016/j.jcp.2024.113384}{Match-based solution of general parametric eigenvalue problems}.
\newblock {\em J. Comput. Phys.}, 519:113384, 2024.

\bibitem{PruRVP02}
\textsc{C.~Prud'homme, D.~V. Rovas, K.~Veroy, and A.~T. Patera}.
\newblock \href{https://doi.org/10.1051/m2an:2002035}{A mathematical and computational framework for reliable real-time solution of parametrized partial differential equations}.
\newblock {\em ESAIM: Math. Model. Numer. Anal.}, 36(5):747--771, 2002.

\bibitem{SRFB2004}
\textsc{U.~Schollw{\"o}ck, R.~J., D.~J.~J. Farnell, and R.~F. Bishop}, editors.
\newblock \href{https://doi.org/10.1007/b96825}{{\em Quantum Magnetism}}.
\newblock Lecture Notes in Physics. Springer Nature, United States, 2004.

\bibitem{Sch11}
\textsc{U.~Schollwoeck}.
\newblock \href{https://doi.org/10.1016/j.aop.2010.09.012}{The density-matrix renormalization group in the age of matrix product states}.
\newblock {\em Ann. Phys.}, 326:96--192, 2011.

\bibitem{Axl15}
\textsc{A.~Sheldon}.
\newblock \href{https://doi.org/10.1007/978-3-319-11080-6}{{\em Linear Algebra Done Right}}.
\newblock Undergraduate Texts in Mathematics. Springer, 3rd edition, 2015.

\bibitem{SirK16}
\textsc{P.~Sirkovic and D.~Kressner}.
\newblock \href{https://doi.org/10.1137/15M1017181}{Subspace acceleration for large-scale parameter-dependent hermitian eigenproblems}.
\newblock {\em {SIAM} J. Matrix Anal. Appl.}, 37(2):1323--1353, 2016.

\bibitem{Wed83}
\textsc{P.~A. Wedin}.
\newblock \href{https://link.springer.com/content/pdf/10.1007/BFb0062089.pdf}{{On Angles between subspaces of a finite dimensional inner product space}}.
\newblock {\em Lect. Notes Math.}, 973:263--285, 1983.

\end{thebibliography}

\end{document}